\documentclass[11pt,oneside,reqno]{amsart}
\usepackage{geometry}
\usepackage[utf8]{inputenc}
\usepackage{amstext,latexsym,amsbsy,amsmath,amssymb,amsthm,mathtools,relsize,geometry,enumerate}
\usepackage{hyperref}
\hypersetup{
  colorlinks   = true, 
  urlcolor     = blue, 
  linkcolor    = red, 
  citecolor   = red 
}
\usepackage{mathtools,enumerate,enumitem}

\usepackage{mathrsfs}

\usepackage{comment}
\usepackage{graphicx}
\usepackage{epstopdf}
\setkeys{Gin}{width=\linewidth,totalheight=\textheight,keepaspectratio}
\graphicspath{{./images/}}

\usepackage{multicol}
\usepackage{booktabs}
\usepackage{mathrsfs}
\usepackage{color}

\newcommand{\rbb}{\mathbb{R}}

\renewcommand{\L}{\mathcal{L}}

\newcommand{\W}{\mathcal{W}}

\newcommand{\Pcal}{\mathcal{P}}
\newcommand{\B}{\mathcal{B}}
\newcommand{\la}{\langle}
\newcommand{\ra}{\rangle}

\newcommand{\X}{\mathbf{X}}

\newcommand{\alphat}{\tilde{\alpha}}
\newcommand{\betat}{\tilde{\beta}}

\newcommand{\D}{\mathcal{D}}

\renewcommand{\d}{\textup{d}}

\newcommand{\f}{\varphi}

\newcommand{\grad}{\nabla}

\newcommand{\Fcal}{\mathcal{F}}
\newcommand{\E}{\mathbb{E}}

\renewcommand{\P}{\mathbb{P}}

\newcommand{\TV}{\textup{TV}}

\theoremstyle{plain}
\newtheorem{theorem}{Theorem}[section]

\newtheorem{lemma}[theorem]{Lemma}
\newtheorem{assumption}[theorem]{Assumption}

\newtheorem{definition}[theorem]{Definition}

\theoremstyle{definition}

\newtheorem{remark}[theorem]{Remark}

\numberwithin{equation}{section}

\title{Lyapunov functions for the Langevin equation with dissipative damping}

\author{Barkat Mian$^1$ and Hung D.~Nguyen$^1$}

\address{$^1$ Department of Mathematics, University of Tennessee, Knoxville, Tennessee, USA}

\begin{document}

\begin{abstract}
    We study the well-posedness and the large time stability of the Langevin equation in the presence of potentials and dissipative damping together with additive white noise. While the confining potential is assumed to be smooth, the nonlinear damping is allowed to have either polynomial growth or singular behaviors. We construct novel Lyapunov functions and determine the corresponding parameter regimes that lead to global well-posedness or finite-time explosion. These Lyapunov functions also allow for deriving quantitative mixing rates toward equilibrium. Our results extend existing ergodic theory established for the instance of nonlinear damping satisfying suitable integrability conditions.       
\end{abstract}

\maketitle

\section{Introduction} \label{sec:intro}

We are interested in the Langevin equation modeling the motion of a particle suspended in a continuous medium and given by
\begin{align}\label{eq:Langevin-SDE}
    \d x & = v\d t,\notag\\
    \d v(t) & = - v(t)\d t -\grad U(x(t))\d t-\grad G(v(t))\d t+\sqrt{\sigma} \d W(t).
\end{align}
In \eqref{eq:Langevin-SDE}, the pair $(x(t),v(t))$ takes values in $\rbb^d\times \rbb^d$ representing the position and the velocity of the particle. The function $-\grad U(x)$ is the spatial derivative of a smooth confining potential $U:\rbb^d\to[0,\infty)$ with polynomial growth. The force term $-\grad G(v)$ describes a nonlinear effect on the damping where $G$ is assumed to exhibit either polynomial growth or singularity at the origin. The particle's trajectories are also posited to be randomly perturbed through $W(t)$, a standard $d$-dimensional Brownian motion with the parameter $\sigma>0$ denoting the noise intensity.

Historically, the Langevin equation with nonlinear potential in the presence of singular effect on the \textit{displacement} has been studied extensively in literature. See \cite{baudoin2021gamma,camrud2021weighted, cooke2017geometric, grothaus2015hypocoercivity,herzog2017ergodicity,lu2019geometric} and the references therein. More precisely, concerning the following equation
\begin{align} \label{eq:Langevin:singular-potential}
  \d x & = v\d t,\notag\\
    \d v(t) & = - v(t)\d t -\grad U(x(t))\d t-\grad G(x(t))\d t+ \d W(t),
\end{align}
where $G$ is a repulsive force acting on $x(t)$, the large-time stability of \eqref{eq:Langevin:singular-potential} is well understood that the system is attractive toward equilibrium exponentially fast. However, compared with the setting of singular potentials, the appearance of nonlinear damping $G(v)$ as in \eqref{eq:Langevin-SDE} seems to receive less attention. Recently, the well-posedness and ergodicity were explored in \cite{song2020well} under the assumption that the damping effect is a singular function satisfying certain integrability conditions. It remains an open problem in \cite{song2020well} on how to apply the technique and the results therein to the case when $G(v)$ is either a polynomial or as singular as the Coulomb and Lennard-Jones functions. The main goal of the present paper is thus to address the stability issue of \eqref{eq:Langevin-SDE} for the aforementioned instances of nonlinearities on the velocity field.

Our main results can be summarized in the following meta theorem.
\begin{theorem}\label{thm:meta:well-posed-mixing} 
    Let $\lambda\ge 2$ and $\alpha>2$ and $\beta\ge 1$ be given constants. Suppose that $U(x)$ satisfies
    \begin{align} \label{cond:meta:U}
        U(x)&\propto |x|^{\lambda}+\textup{lower order terms},\quad x\in\rbb^d,
    \end{align}
and that $G(v)$ satisfies 
    \begin{align} 
        \text{either}\qquad G(v)&\propto |v|^{\alpha}+\textup{lower order terms},\quad v\in\rbb^d;\label{cond:meta:G:poly}\\
        \text{or}\qquad \grad G(v)&\propto -\frac{v}{|v|^{\beta+1}}+\textup{lower order terms},\quad v\in\rbb^d\setminus\{0\}.\label{cond:meta:G:singular}
    \end{align}
Then, under suitable assumptions on the lower order terms, the following hold: 

1. Well-posedness: equation \eqref{eq:Langevin-SDE} admits a unique global solution except for the case $d=1$, $\beta=1$ and $\sigma$ is sufficiently large, where finite time explosions occur a.s.

2. Ergodicity and mixing: 

\hspace{0.5cm} \textup{(i)} If $G(v)$ satisfies \eqref{cond:meta:G:poly}, then \eqref{eq:Langevin-SDE} converges toward equilibrium with either exponential rates or at least algebraic rates of arbitrary orders corresponding to either $\lambda\ge \alpha$ or $\lambda<\alpha$, respectively.

\hspace{0.5cm} \textup{(ii)} If $G$ satisfies \eqref{cond:meta:G:singular}, then in dimension $d=1$, there does not exist invariant probability measures whereas in dimension $d\ge 2$, under further relations between $\lambda$ and $\beta$, namely,
\begin{align*}
     \beta\ge 1,\quad 2\le \lambda\le 4 \quad \text{or}\quad  1\le \beta <\frac{6(\lambda-1)}{\lambda-4},\quad \lambda>4 ,
\end{align*}
\eqref{eq:Langevin-SDE} converges toward equilibrium with at least algebraic rates of arbitrary orders.
\end{theorem}

We refer the reader to Theorem \ref{thm:well-posed} for the rigorous statement of Theorem \ref{thm:meta:well-posed-mixing}, part 1 on the well-posedness of \eqref{eq:Langevin-SDE} and to Theorem \ref{thm:mixing} for Theorem \ref{thm:meta:well-posed-mixing}, part 2 on the ergodic behaviors. We note that Theorem \ref{thm:meta:well-posed-mixing}, part 2 (ii) covers for example when $U(x)$ is a quartic polynomial and $G$ is either Coulomb or Lennard-Jones functions in arbitrary dimension $d\ge 2$. 

\subsection{Literature background} Research on the large-time stability of the Langevin equation first appeared in the 1970s in \cite{tropper1977ergodic} and was later popularized by the seminal work of \cite{mattingly2002ergodicity, talay2002stochastic} dealing with exponential mixing for the setting with smooth confining potential $U$, i.e., $G\equiv 0$ in \eqref{eq:Langevin-SDE}. The issue of ergodicity becomes considerably more challenging in \eqref{eq:Langevin:singular-potential}, which incorporates the repulsive force $G(x)$ on the displacement. Efforts in this direction started as early as in \cite{conrad2010construction,grothaus2015hypocoercivity} using the hypocoercivity method to derive convergence rates for Lennard-Jones type functions in the setting of multiple interacting particles. Alternatively, the approach of Lyapunov function was explored in \cite{cooke2017geometric, herzog2017ergodicity,lu2019geometric}, which have proven to be applicable to a broad class of singular potentials including both the Coulomb and Lennard-Jones functions. While qualitative results on mixing were established in  \cite{grothaus2015hypocoercivity,herzog2017ergodicity,lu2019geometric}, more recent work of \cite{baudoin2021gamma,camrud2021weighted} has focused on obtaining explicit quantitative estimates on the convergence rates by making use of the Gamma calculus. For a study on ergodicity of related Hamiltonian systems, we refer the reader to \cite{cuneo2018non, hairer2009hot,hairer2009slow,mertz2026controllability, nguyen2025geometric}

Turning back to equation \eqref{eq:Langevin-SDE}, while the presence of polynomial damping in second order systems arises naturally in a variety of applications, e.g., the study of stochastic non-linear oscillators \cite{markus1988stochastic,markus1989oscillation,markus1993stochastic}, earthquake-induced and structural vibrations \cite{chatterjee2020identification,fan2004surface,phillips2009damping}, acoustic wave
propagation in biomedical applications \cite{khokhlova2006effects}, and nonlinear optics \cite{hasegawa1973}, the consideration of singular damping is motivated by singular structures such as point vortices and dipoles in ideal incompressible fluid dynamics \cite{kulik2010interaction,smith2011singularities,yanovsky2009singularities}.

With regards to the well-posedness and large-time stability of Langevin with nonlinear damping, the existence of global solutions of \eqref{eq:Langevin-SDE} for smooth $G(v)$ was investigated in \cite{markus1993stochastic} under a dissipation condition involving $LH$ where $L$ is the generator of \eqref{eq:Langevin-SDE} and $H$ is the usual Hamiltonian structure $H(x,v)=U(x)+\frac{1}{2}|v|^2$. However, to the best of the authors' knowledge, there does not seem to be an existing work for the study of ergodicity under smooth nonlinear damping. In this regard, the result of Theorem \ref{thm:meta:well-posed-mixing}, part 2(i) seems to be the first established in this direction for the Langevin systems. We also would like to mention that related studies have explored analogous questions in infinite-dimensional settings, particularly for stochastic wave equations, where polynomial mixing rates have been established; see, for instance \cite{barbu2007stochastic,gao2026ergodicity,nguyen2024polynomial}. In contrast, the consideration of singular-typed damping appeared in \cite{song2020well}, which investigates the equation
\begin{align}\label{eq:Langevin-SDE:Song}
    \d x & = v\d t,\notag\\
    \d v(t) & = - v(t)\d t -\grad U(x(t))\d t+F(v(t))\d t+\sqrt{\sigma} \d W(t).
\end{align}
Under the assumption that $F\in L^p(\rbb^d)$ for suitable $p\ge 1$ determined by the growth rate of the potential $U(x)$, namely, 
\begin{align*}
     \max\{2,d\}<p\le \infty \quad\text{when}\quad  \lambda=2, \quad\text{and}\quad   2d<p\le \infty \quad\text{when}\quad   \lambda>2,
\end{align*}
it is established in \cite[Theorem 1.2]{song2020well} that \eqref{eq:Langevin-SDE:Song} possesses unique global solution and furthermore is exponential attractive toward the unique invariant probability measure. Specific instances of singular $F$ within the scope of \cite{song2020well} include truncated functions of the form 
\begin{align*}
    F(v)=\frac{v}{|v|^{\beta+1}}\mathbf{1}\{|v|\le K\},
\end{align*} for arbitrary constant $K>0$ and appropriate constant $\beta>0$. The approach employed in \cite{song2020well} for the well-posedness makes use of the Zvonkin's transform to regularize the singularity in $F$ whereas the Krylov's estimate for parabolic equations together with a Lyapunov function was employed to study the problem of mixing rates. It is worth to mention that while the phase space considered in \cite{song2020well} is the whole space $\rbb^d\times\rbb^d$ as an implication of the singular $F$ in \eqref{eq:Langevin-SDE:Song}, in our setting of \eqref{eq:Langevin-SDE} when $G$ is singular as in condition \eqref{cond:meta:G:singular}, the phase space is restricted to $\rbb^d\times(\rbb^d\setminus\{0\})$ in dimension $d\ge 2$ and $\rbb\times (0,\infty)$ in dimension $d=1$ where the velocity is confined in the positive real line. This is a manifestation of the fact that the repulsive force $-\grad G(v)$ prevents the origin from being reachable by the velocity, except when finite time explosion occurs. This distinction highlights a fundamental difference between the singular perturbations treated in \cite{song2020well} and the singular damping considered in the present work.

\subsection{{Methodology of the proofs}} With regard to the proof of Theorem \ref{thm:meta:well-posed-mixing}, part 1 on the well-posedness of \eqref{eq:Langevin-SDE}, since the nonlinear forcing is locally Lipschitz, it is a classical result that local solutions always exist up to a random explosion time. In order to achieve the global well-posedness, we construct suitable Lyapunov functions ensuring the total energy is appropriately controlled without escaping to infinity. While such estimate for the smooth damping given by condition \eqref{cond:meta:G:poly} is relatively simple making use of the Hamiltonian structure as discussed in \cite{markus1993stochastic}, the argument for singular damping corresponding to condition \eqref{cond:meta:G:singular} is more delicate. Unlike the approach of \cite{song2020well} involving a regularization procedure, we opt for treating the issue of global solutions by directly constructing different Lyapunov functions, depending on the dimension $d$, the leading singular parameter $\beta$ and the noise intensity $\sigma$. Our analysis also reveals a regime in which well-posedness fails. That is in the case when $d=1$, $\beta=1$ and $\sigma$ is sufficiently large, the ill-posedness occurs as stated in Theorem \ref{thm:meta:well-posed-mixing}, part 1, in the sense that the velocity can reach the origin in finite time a.s. This result actually is not surprising resembling the behaviors of the Bessel process which exhibits an analogous property of hitting the origin under suitable parameter conditions. In fact, the proof of finite time explosion in this specific scenario involves a comparison argument where we will relate the trajectories of the local solution of \eqref{eq:Langevin-SDE} to those of a Bessel-like process. We refer the reader to Appendix \ref{sec:illposed:auxiliary-equation} for a detailed discussion of this auxiliary process. The complete well-posedness study will be given in Section \ref{sec:wellposed}, where we supply the proof of Theorem \ref{thm:well-posed}, the rigorous statement of Theorem \ref{thm:meta:well-posed-mixing}, part 1.

Turning to the second main topic of the paper on the ergodicity, we proceed to establish Theorem \ref{thm:meta:well-posed-mixing}, part 2 by combining the Lyapunov technique with hypoellipticity and the solvability of an associated control problem. While the latter two properties are quite standard following the approach of \cite{duong2024asymptotic,hairer2009hot,herzog2017ergodicity,lu2019geometric,mattingly2002ergodicity} established for related Langevin systems, the main crucial ingredient of the proof is the existence of strong Lyapunov functions. We note that for the well-posedness, it is sufficient to derive a weak Lyapunov bound of the form $LV\le cV+C$, whereas for ergodicity, we need to construct a function $V$ satisfying 
\begin{align} \label{ineq:meta:Lyapunov}
    LV\le -c V^p+C,
\end{align}
where $p>0$ is a constant ensuring the contraction of the dynamics toward equilibrium.

In the setting of polynomial damping given by condition \eqref{cond:meta:G:poly}, as it turns out, such a function $V$ is determined by the relation between two leading parameters $\lambda$ and $\alpha$ from \eqref{cond:meta:U} and \eqref{cond:meta:G:poly}, respectively. On the one hand, when $\lambda\ge \alpha$, we follow the usual perturbation $H(x,v)+\varepsilon\la x,v\ra$ with $\varepsilon$ small enough to show that $LH\le -CH+C$, which produces exponential convergence rate. This result is a reminiscence of analogous behaviors of the standard Langevin equation established in \cite{mattingly2002ergodicity}, since in this regime, the potential is regarded as the leading force and dominates the nonlinear damping. On the other hand, when $\lambda<\alpha$, since neither the potential nor the damping can be treated as a lower-order contribution, we resort to the approach of \cite{barbu2007stochastic,gao2026ergodicity,nguyen2024polynomial} explored for the stochastic wave equations. More specifically, we augment the Hamiltonian by the term $G(v)+\la v,\grad U(x)\ra$, which produces an extra dissipation through $-|\grad G(v)+\grad U(x)|^2$ after applying the generator $L$. This allows us to successfully obtain a Lyapunov function $V$ satisfying \eqref{ineq:meta:Lyapunov} for a suitable exponent $p\in (0,1)$ and thus extract polynomial rates of any order. The specific construction will be presented in Section \ref{sec:ergodicity:poly-damping}.

With regard to the singular damping as in condition \eqref{cond:meta:G:singular}, the Lyapunov technique becomes substantially challenging since we have to treat simultaneously the regimes of large displacement and small velocity as $v$ approaches the origin. Nevertheless, in dimension $d=1$, since the velocity is posed on the positive real line, the position trajectory is always strictly increasing. Heuristically, this implies that the system is not able to return to bounded sets, leading to the non-existence of invariant measures. In higher dimension $d\ge 2$, when $\beta=1$, since $G$ is proportional to $-\log|v|$, we employ the construction for the smooth damping discussed above and modify the Hamiltonian by the correction $-\log|v|+\la v,\grad U(x)\ra$. When $\beta>1$, we draw upon the technique of Lyapunov functions from \cite{duong2024asymptotic, lu2019geometric} dealing with similar issue for singular damping described by equation \eqref{eq:Langevin:singular-potential}. More specifically, for \eqref{eq:Langevin:singular-potential}, we recall that a Lyapunov function has the form $H(x,v)+\la x,v\ra-\la x,v\ra/|x|$ where the last term incorporates the singular effect in the displacement direction. Motivated by this trick, for \eqref{eq:Langevin-SDE}, we replace the term $G(v)+\la v,\grad U(x)\ra$ by a weighted version $|v|^{-q}[G(v)+\la v,\grad U(x)\ra]$ for a suitable constant $q\in (0,1)$. In turn, this enables us to achieve the Lyapunov bound \eqref{ineq:meta:Lyapunov} with $p\in(0,1)$, yielding algebraic rates of arbitrary order. It is worthwhile to mention that a fundamental difference of \eqref{eq:Langevin-SDE} from the \eqref{eq:Langevin:singular-potential} is the interaction between noise and singularity. In \eqref{eq:Langevin-SDE}, we need to facilitate the noise to produce a dissipation in the $v$ direction when $v$ is small. In contrast, the singular force in \eqref{eq:Langevin:singular-potential} only acts on the $x$-direction whereas the noise appears directly in the velocity equation and thus does not induce a strong effect on small $x$. Finally, we remark that our results only apply to smooth functions $U(x)$ and do not extend to the regime of repulsive potentials, nor the setting of multi-particle with interaction potentials. We opt for leaving these problems as future work. The detailed construction will be carried out in Section \ref{sec:ergodicity:singular-damping} while the proof of Theorem \ref{thm:mixing}, the rigorous statement of Theorem \ref{thm:meta:well-posed-mixing} part 2, will be given in Section \ref{sec:ergodicity-mixing}.

\subsection{Organization of the paper}
The rest of the paper is organized as follows: in Section \ref{sec:main-results}, we introduce the notations as well as the assumptions that we make on the nonlinearities. We
also state the main results of the paper, including Theorem \ref{thm:well-posed} on the well-posedness of \eqref{eq:Langevin-SDE} and Theorem \ref{thm:mixing} on its ergodic properties. In Section \ref{sec:wellposed}, we detail the energy estimates for the existence of global solutions and a comparison argument for the finite time explosion. The proof of Theorem \ref{thm:well-posed} is also supplied in this section. In
Section \ref{sec:ergodicity}, we present the construction of Lyapunov functions allowing for extracting the mixing rates. In Section \ref{sec:ergodicity-mixing}, we prove Theorem \ref{thm:mixing} by combining the Lyapunov functions of Section \ref{sec:ergodicity} with auxiliary results on hypoellipticity, the associated control problem and the non-existence of invariant probability measures in dimension $d=1$ under singular damping. 
The paper concludes with Appendices \ref{sec:illposed:auxiliary-equation} and \ref{sec:LVFour-Weighted-Correction-Derivation}. In Appendix \ref{sec:illposed:auxiliary-equation}, we study a Bessel-like process in dimension $d=1$ and establish the illposedness property, the result of which was employed in the proof of Theorem \ref{thm:well-posed}. In Appendix \ref{sec:LVFour-Weighted-Correction-Derivation}, we supply a series of calculations, which appeared in the derivation of Lyapunov bounds in Section \ref{sec:ergodicity}.

\section{Assumptions and main results} \label{sec:main-results}

Throughout, we fix a stochastic basis $(\Omega,\Fcal,(\Fcal_t)_{t\ge 0}, \P)$ satisfying the usual conditions \cite{karatzas2012brownian} and let $W(t)$ be a standard $d$-dimensional Brownian motion on $(\Omega,\Fcal,\P)$ adapted to the filtration $\Fcal_t$. 

In this section, we state the main assumptions that we will employ for the analysis, as well as the main results on the well-posedness and large time stability of system \eqref{eq:Langevin-SDE}. More precisely, in Section \ref{sec:main-results:well-posed}, we detail a general set of conditions that we impose on the potential $U(x)$ and the nonlinear damping $G(v)$ allowing us to establish the well-posedness of system \eqref{eq:Langevin-SDE} through Theorem \ref{thm:well-posed}. In particular, while the global solutions are guaranteed to uniquely exist for a majority of the hypotheses, the finite time explosion occurs in dimension $d=1$ when the noise parameter dominates the singular effect from the damping. In Section \ref{sec:main-results:ergodicity}, we establish the ergodic properties of \eqref{eq:Langevin-SDE} together characterizing mixing behaviors in the parameter regimes corresponding to the global well-posedness. Notably, while the unique ergodicity is guaranteed for a large class of nonlinearities, the non-existence of invariant probability measures is exhibited in the case of singular damping in dimension $d=1$.   

\subsection{Well-posedness} \label{sec:main-results:well-posed} To derive the well-posedness of equation \eqref{eq:Langevin-SDE}, we start with the potential $U$ and impose the following popular conditions \cite{mattingly2002ergodicity,pavliotis2014stochastic}.

\begin{assumption} \label{cond:U}

\textup{(i)} There exist constants $a_1>0$ and $\lambda\ge 2$ such that $U\in C^\infty(\rbb^d;[0,\infty))$ satisfies for  all $x\in\rbb^d$
\begin{align}
\frac{1}{a_1}|x|^{\lambda}-a_1 &\le  U(x)\le a_1(1+|x|^{\lambda}),\label{cond:U:U(x)=O(x^lambda+1)}  \\
 |\grad U(x)| & \le a_1(1+|x|^{\lambda-1}).\label{cond:U:U'(x)=O(x^lambda)} \\
\text{and}\quad|\nabla^2U(x)|
&\le
a_1\left(1+|x|^{\lambda-2}\right).\label{cond:U:U''(x)=O(x^lambda)}
\end{align}

\textup{(ii)} Furthermore,
\begin{equation} \label{cond:U:x.U'(x)>-x^(lambda+1)}
\la \grad U(x), x\ra \ge a_2|x|^{\lambda}-a_3,\quad x\in\rbb^d,
\end{equation}
for some positive constants $a_2$ and $a_3$.

\end{assumption}

Turning to the dissipative damping $G(v)$, in Assumption \ref{cond:G:poly} below, we impose that $G$ is a smooth function satisfying a polynomial growth.

\begin{assumption}[Polynomial damping] \label{cond:G:poly}
    $G\in C^\infty(\rbb^d;[0,\infty))$ satisfies for all $v \in \rbb^d$
    \begin{align}
        \big|\grad G(v)- a_4|v|^{\alpha-2} v\big| & \le a_5\big(|v|^{\alphat-1}+1\big), \label{cond:G:poly:grad.G}\\
        and\quad \big|\grad^2 G(v)\big|& \le a_5(|v|^{\alpha-2}+1),\label{cond:G:poly:grad^2.G}
    \end{align}
    for some constants $\alpha > 2$, $1\le \alphat<\alpha$, $a_4>0$ and $a_5\ge 0$.
\end{assumption}

We note that while condition \eqref{cond:G:poly:grad.G} is slightly stronger than the more standard coercive condition \eqref{cond:U:x.U'(x)>-x^(lambda+1)} imposed on $U(x)$, condition \eqref{cond:G:poly:grad.G} essentially requires that $G$ behave as a polynomial function of the form $|v|^{\alpha}$ with $\alpha>2$. 
Alternatively, as mentioned in the introduction, we will also consider a class of damping $G(v)$ exhibiting singular behavior at the origin. This is stated precisely below in Assumption \ref{cond:G:singular}. 

\begin{assumption}[Singular damping] \label{cond:G:singular}
    $G\in C^\infty(\rbb^d\setminus\{0\}; \rbb)$ satisfies
    \begin{align} \label{cond:G:singular:|gradG|<1/|v|^beta}
    \bigg|\grad G(v) +a_4\frac{v}{|v|^{\beta+1}} \bigg|\le |a_6|\bigg(\frac{1}{|v|^{\betat}}+1\bigg).
\end{align}
    for some constants $\beta\ge 1$, $a_4> 0$, $a_6\in\rbb$, and $0\le \betat <\beta$.
    
\end{assumption}

\begin{remark}
1. Unlike condition \eqref{cond:G:poly:grad.G} imposing polynomial-like behaviors, condition \eqref{cond:G:singular:|gradG|<1/|v|^beta} states that $G$ can be expressed as the form $\grad G(v) = -a_4 v/|v|^{\beta+1}+$ lower order terms. In particular, from both Assumption \ref{cond:G:poly} and Assumption \ref{cond:G:singular}, $G$ is not $L^p$-integrable on the respective domains $\rbb^d$ and $\rbb^d\setminus\{0\}$. 

2. As mentioned elsewhere in \cite{duong2024asymptotic,duong2026trend}, a routine calculation shows that condition \eqref{cond:G:singular:|gradG|<1/|v|^beta} covers the repulsive Coulomb functions
\begin{align*} 
G(v)=\begin{cases}
-\log|v|,& d=2,\\
\frac{1}{|v|^{d-2}},& d\ge 3,
\end{cases}
\end{align*}
with $\beta=d-1$ and $a_6=0$ (for which $\betat$ is arbitrary), as well as the Lennard-Jones functions 
\begin{align*}
G(v) = \frac{c_0}{|v|^{12}}-\frac{c_1}{|v|^6},
\end{align*} 
with $\beta = 13$ and $\betat=7$.

\end{remark}

Having stated the sufficient conditions on the nonlinearities, we address the first topic of the paper concerning the well-posedness of \eqref{eq:Langevin-SDE}. In order to do so, we introduce the phase space for the pair $(x,v)$ given by
$$\X = \rbb^d\times\D,$$ 
where $\D$ is the phase space for the velocity depending on the choice of $G$ and defined as
\begin{align*}
    \D = \begin{cases}
\rbb^d,& \text{polynomial damping},\\
\rbb^d \setminus \{0\},& \text{singular damping},\quad d\ge 2,\\
(0,\infty),& \text{singular damping},\quad d= 1.
    \end{cases}
\end{align*}
In Theorem \ref{thm:well-posed} stated below, we characterize the parameter regimes that lead to either the existence and uniqueness of the global solutions or finite time explosion.

\begin{theorem} \label{thm:well-posed} Let $U$ satisfy Assumption \ref{cond:U}.

1. Suppose that $G$ satisfies either Assumption \ref{cond:G:poly} or Assumption \ref{cond:G:singular}. Moreover, for the latter assumption, at least one of the following conditions hold:
\begin{align} \label{cond:d-beta-sigma:well-posed}
    d>1,\quad \textup{or}\quad \beta>1,\quad \textup{or}\quad \sigma\le 2a_4,
\end{align}
where $\beta,\, a_4$ are the constants as in Assumption \ref{cond:G:singular} and $\sigma$ is the noise parameter in \eqref{eq:Langevin-SDE}. Then, for all initial conditions $(x_0,v_0)\in \X$, there exists a unique strong solution for equation \eqref{eq:Langevin-SDE}.

2. Suppose that $G$ satisfies Assumption \ref{cond:G:singular} and that
\begin{align} \label{cond:d-beta-sigma:ill-posed}
    d=1,\qquad \beta=1,\qquad \sigma> 2a_4.
\end{align}
Then, finite time explosion occurs a.s. That is there exists a random variable $T>0$ such that for all $(x_0,v_0)\in \X=\rbb\times (0,\infty)$, almost surely,
\begin{align*}
    T<\infty, \quad\textup{and}\quad \lim_{t\to T^-}v(t)=0.
\end{align*}
    
\end{theorem}

The proof of Theorem \ref{thm:well-posed}, part 1, relies on the standard technique of Lyapunov bounds ensuring the energy of the system does not escape to infinity. On the one hand, we note that the estimates for the case of polynomial damping, i.e., when Assumption \ref{cond:G:poly} holds, is relatively straightforward following the Hamiltonian structure, cf. Lemma \ref{lem:well-posed:poly-damping}. On the other hand, the construction of the Lyapunov functions in the case of singular damping as in Assumption \ref{cond:G:singular} is non trivial, requiring the validity of at least one of the conditions in \eqref{cond:d-beta-sigma:well-posed}. We refer the reader to Lemma \ref{lem:well-posed:singular-damping} for a detailed derivation. Altogether, they allow us to invoke a stopping time argument to extend local solutions to global well-posedness, of which the explicit explanation will be provided in Section \ref{sec:wellposed}.  

While the well-posedness draws upon the usual approach for SDEs, the finite time explosion in the case of singular damping satisfying condition \eqref{cond:d-beta-sigma:ill-posed} is much more delicate. In particular, the threshold $\sigma >2a_4$ stems from the classical fact that the one-dimensional Bessel process 
\begin{align*}
    \d v(t) = \frac{a}{v(t)}\d t + \sqrt{\sigma}\d W(t)
\end{align*}
reaches the origin in finite time when $\sigma>2a$; see \cite[Chapter XI, Section 1, p.~442]{revuz2013continuous}, in the paragraph preceding Proposition~1.5. To demonstrate the finite time explosion for the solutions of \eqref{eq:Langevin-SDE}, we will use a direct comparison argument and show that when \eqref{cond:d-beta-sigma:ill-posed} holds, the behaviors of \eqref{eq:Langevin-SDE} can be related to a Bessel-like process that also exhibits the origin-reachability property. In Appendix \ref{sec:illposed:auxiliary-equation}, we discuss this latter process in details, whereas in Section \ref{sec:wellposed:singular-damping}, Lemma \ref{lem:well-posed:singular-damping:sigma>2:illposed}, we carry out the comparison argument, allowing us to establish Theorem \ref{thm:well-posed}, part 2.

\subsection{Ergodicity and mixing} \label{sec:main-results:ergodicity}

As a consequence of the global well-posedness formulated in Theorem \ref{thm:well-posed}, part (i), under the same hypotheses, we can thus introduce the Markov transition probabilities of the solution $\big(x(t),v(t)\big)$ by
\begin{align*}
P_t\big((x_0,v_0),A\big) :=\P\big( (x(t),v(t))\in A|(x(t),v(t))=(x_0,v_0)\big),
\end{align*} 
which are well-defined for $t\ge 0$, initial condition $(x_0,v_0)\in \X$ and Borel sets $A\subset \X$.  Letting $\B_b(\X)$ denote the set of bounded Borel measurable functions $f:\X \rightarrow \rbb$, the associated Markov semigroup $P_t:\B_b(\X)\to\B_b(\X)$ is defined as
\begin{align*}
P_t f(x_0,v_0)=\E_{(x_0,v_0)}\big[f(x(t),v(t)))\big], \,\, f\in \B_b(\X).
\end{align*}
Letting $\Pcal r(\X)$ be the space of probability measures on $\X$, we denote by $P_t\nu$ the push-forward measure of $\nu\in\Pcal r(\X) $ under the action of $P_t$. That is
\begin{align*}
P_t\nu(A) = \int_{\X}P_t\big((x_0,v_0),A\big)\nu(\d x_0,\d v_0),
\end{align*}
for all Borel sets $A\subset \X$. Recall that $\nu\in \Pcal r(\X)$ is called invariant for $P_t$ if $P_t\nu = \nu$, $t\ge 0$. We will measure the convergence rate toward equilibrium through the total variation distance between two probability measures $\nu_1,\nu_2\in \Pcal r(\X)$ defined as
\begin{align*}
    \W_{\TV}(\nu_1,\nu_2) = \inf \E[\mathbf{1}\{X_1\neq X_2\}],
\end{align*}
where the infimum runs over all bivariate random variables $(X_1,X_2)\in \X\times \X$ such that $X_1\sim \nu_1$ and $X_2\sim\nu_2$.

Owing to the interaction between the singular damping $G(v)$, the potential $U(x)$ and the noise term, when Assumption \ref{cond:G:singular} holds with $\beta>1$, we will impose further restriction on $G$ as well as the parameters representing the dominating behaviors of $G$ and $U$. This is precisely stated in the following assumption. 

\begin{assumption} \label{cond:G:singular:beta>1}
    Let $\lambda$ be the constant as in Assumption \ref{cond:U} and $\beta,\betat,a_4,a_6$ be the constants from Assumption \ref{cond:G:singular}. When $\beta>1$, we assume further that $G$ satisfies
    \begin{align} \label{cond:G:singular:gradG=1/|v|^beta:beta>1:mixing}
    \grad G(v) = -a_4\frac{v}{|v|^{\beta+1}} - a_6\frac{v}{|v|^{\betat+1}},
\end{align}
and that
\begin{align} \label{cond:G:singular:gradG=1/|v|^beta:beta>1:mixing:lambda-beta}
\begin{cases}
\beta>1, & 2\le\lambda\le4,\\[1mm]
1<\beta<\dfrac{6(\lambda-1)}{\lambda-4}, & \lambda>4.
\end{cases}
\end{align}

\end{assumption}

\begin{remark} \label{rem:G:singular:mixing}
     It is clear that on the one hand, the extra conditions \eqref{cond:G:singular:gradG=1/|v|^beta:beta>1:mixing}-\eqref{cond:G:singular:gradG=1/|v|^beta:beta>1:mixing:lambda-beta} still cover the instances of singular damping as Coulomb and Lennard-Jones functions, when the potential $U$ grows at most like a quartic polynomial, i.e., $\lambda\le 4$. On the other hand, as $\lambda$ gets large, the damping $G(v)$ is allowed to be as singular as $1/|v|^5$. In particular, for all $\lambda$ large enough, condition \eqref{cond:G:singular:gradG=1/|v|^beta:beta>1:mixing:lambda-beta} excludes the the Coulomb function in dimension $d\ge 8$ and the Lennard-Jones function.
\end{remark}

We are now in a position to state the main result of the paper characterizing the ergodicity of system \eqref{eq:Langevin-SDE}, as well as the mixing rates toward equilibrium.

\begin{theorem} \label{thm:mixing}
Let $U$ satisfy Assumption \ref{cond:U}.

1.  Suppose that $G$ satisfies Assumption \ref{cond:G:poly}. Then, there exists a unique invariant probability measure $\nu$ for $P_t$. Moreover, let $\lambda$ and $\alpha$ be the constants as in Assumption \ref{cond:U} and Assumption \ref{cond:G:poly}, respectively. Then, 

\hspace{0.5cm}\textup{(i)} if $\alpha\le \lambda$, there exists a function $V:\X\to [1,\infty)$ such that the following holds
\begin{align*}
    \W_{\TV}\big(P_t((x_0,v_0),\,\cdot\,),\nu\big) \le Ce^{-ct}V(x_0,v_0),\quad t\ge 0,\quad (x_0,v_0)\in\X,
\end{align*}
for some positive constants $c,C$ independent of $t$ and $(x_0,v_0)$;

\hspace{0.5cm}\textup{(ii)} otherwise, if $\alpha> \lambda$, for all $n\ge 1$, there exists a function $V_n:\X\to [1,\infty)$ such that 
\begin{align*}
    \W_{\TV}\big(P_t((x_0,v_0),\,\cdot\,),\nu\big) \le  \frac{C}{(1+t)^n}V_n(x_0,v_0),\quad t\ge 0,\quad (x_0,v_0)\in\X.
\end{align*}

2. Suppose that $G$ satisfies Assumption \ref{cond:G:singular} in dimension $d=1$ and that at least one of the following conditions hold:
\begin{align*}
    \beta>1 \quad \textup{or}\quad \sigma\le 2a_4,
\end{align*}
where $\beta,\, a_4$ are the constants as in Assumption \ref{cond:G:singular} and $\sigma$ is the noise parameter in \eqref{eq:Langevin-SDE}. Then, there does not exist an invariant probability measure for $P_t$.

3. Suppose that in dimension $d\ge 2$, $G$ satisfies Assumption \ref{cond:G:singular} when $\beta=1$ and Assumption \ref{cond:G:singular:beta>1} when $\beta>1$. Then, there exists a unique invariant probability measure $\nu$ for $P_t$. Moreover, for all $n\ge 1$, there exists a function $V_n:\X\to [1,\infty)$ such that 
\begin{align*}
    \W_{\TV}\big(P_t((x_0,v_0),\,\cdot\,),\nu\big) \le  \frac{C}{(1+t)^n}V_n(x_0,v_0),\quad t\ge 0,\quad (x_0,v_0)\in\X.
\end{align*}

\end{theorem}

Since system \eqref{eq:Langevin-SDE} is clearly hypoelliptic, cf. Lemma \ref{lem:Hormander}, the unique ergodicity and mixing rates of Theorem \ref{thm:mixing} are mainly a consequence of the existence of Lyapunov functions, which are much more delicate than the corresponding analysis derived for the well-posedness in Theorem \ref{thm:well-posed}. In particular, while the Lyapunov functions used to handle the polynomial damping in Theorem \ref{thm:mixing}, part 1, are motivated by the related work of \cite{barbu2007stochastic,nguyen2024polynomial} on the stochastic wave equations, the construction for the singular damping as stated in Theorem \ref{thm:mixing}, part 3, is drawn upon the previous literature established for the setting Langevin equation with repulsive potentials in \cite{duong2024asymptotic,lu2019geometric}. Both results together with a study of the associated control problem, in the spirit of the Strook-Varadhan Support Theorem, will allow us to deduce the corresponding mixing rates. In contrast, the one-dimensional singular damping induces the positivity of the velocity process, forcing the displacement trajectories strictly increasing over time. The system therefore is not able to relax to equilibrium, hence the non-existence of invariant probability measures as captured in Theorem \ref{thm:mixing}, part 2. In Section \ref{sec:ergodicity}, we present the construction of Lyapunov functions in detail whereas the proof of Theorem \ref{thm:mixing} will be given in Section \ref{sec:ergodicity-mixing}. This section also develops the auxiliary results required for the ergodicity, including the verification of H\"ormander's condition and the study of the associated control problem.

\section{Well-posedness} \label{sec:wellposed}

In this section, we prove Theorem \ref{thm:well-posed} by establishing the well-posedness and the finite time explosion for \eqref{eq:Langevin-SDE}. In Section \ref{sec:wellposed:explosion-time}, we follow the standard theory of SDEs to obtain the local solutions up to an explosion time. In Section \ref{sec:wellposed:poly-damping} where the damping $G$ is assumed to have polynomial growth, we obtain energy estimate by means of Hamiltonian structures to deduce the explosion time is infinite a.s., thereby yielding the global solutions. The same strategy through more complicated Lyapunov functions is also employed in Section \ref{sec:wellposed:singular-damping} under the assumption of singular dampings to deduce the well-posedness. At the end of Section \ref{sec:wellposed:singular-damping}, we carry out a comparison argument tailored to dimension $d=1$ to show finite time exlosion occurs a.s. when condition \eqref{cond:d-beta-sigma:ill-posed} holds.

\subsection{Local solutions}
\label{sec:wellposed:explosion-time}

For each choice of the damping term, the state space $\X = \rbb^d\times\D$ is an open subset of a Euclidean space and therefore a locally compact metric space. Following~\cite[Def.~1.1.4, p.~11]{Hsu}, let $\X^c:=\X\cup\{\partial_{\X}\}$ denote the one-point compactification of $\X$. An $\X$-valued path $y$ with explosion time $T\in(0,\infty]$ is a continuous map $y: [0, \infty) \to\X^c$ such that $y(t)\in\X$ for $0\le t<T$ and, if $T<\infty$, $y(t)=\partial_{\X}$ for all $t\ge T$. Thus, when $T=\infty$, the path remains in $\X$ for all $t\ge0$. The explosion time can equivalently be characterized through successive exit times associated with an exhaustion of the state space. More precisely, let $\{O_n\}_{n\ge1}$ be an exhaustion of $\X$ by relatively compact open sets, that is,
\[
\overline{O_n}\subset O_{n+1},
\qquad
\bigcup_{n\ge1}O_n=\X.
\]
For an $\X$-valued path $y$ with explosion time $T$, define
\begin{align}\label{eq:General-Exit-Times}
T_{O_n}
:=
\inf\left\{ t \ge 0: y(t)\notin O_n \right\}.
\end{align}
Then, by~\cite[Prop.~1.1.5\textup{(1)}, p.~11]{Hsu},
\begin{align}\label{eq:General-Explosion-Time}
T_{O_n}\uparrow T \qquad \text{as }n\to\infty.
\end{align}
We use the notion of a solution up to a stopping time in the sense of~\cite[Def.~1.1.7, p.~12]{Hsu} for the polynomial damping case, where $\X = \mathbb R^{2d}$; in the singular damping case, where $\X$ is regarded as a smooth manifold, we use the corresponding notion in~\cite[Def.~1.2.3, p.~20]{Hsu}. Under either damping assumption, the coefficients of~\eqref{eq:Langevin-SDE} are smooth on $\X$. With these conventions, we establish the existence and uniqueness of solutions to~\eqref{eq:Langevin-SDE} up to the explosion time in the following lemma.

\begin{lemma}\label{Lem:wellposed:existence-up-to-explosion}
Suppose that $U$ satisfies Assumption~\ref{cond:U} and that $G$ is given either by Assumption~\ref{cond:G:poly} or by Assumption~\ref{cond:G:singular}. For every $(x_0,v_0)\in\X$, there exists a unique solution $\{(x(t),v(t))\}_{t\in[0,T)}$ to~\eqref{eq:Langevin-SDE} taking values in $\X$ up to its explosion time $T$.
\end{lemma}

\begin{proof}[Proof of Lemma~\ref{Lem:wellposed:existence-up-to-explosion} ] We treat the two choices of the damping term separately. Suppose first that $G$ is given by Assumption~\ref{cond:G:poly}. In this case, $\X=\mathbb R^{2d}$. Set $Y(t):=(x(t),v(t))^{\mathsf T}\in\mathbb R^{2d}$ and $Z(t):=(W(t),t)^{\mathsf T}\in\mathbb R^{d+1}$, where $ x(t), v(t)\in\mathbb R^d$, $t\in[0,\infty)$, and $W$ is a standard $d$-dimensional Brownian motion. Then~\eqref{eq:Langevin-SDE} can be written in the form $\d Y(t)=\Sigma(Y(t))\,\d Z(t)$, where $\Sigma:\mathbb R^{2d} \longrightarrow M(2d,d+1)$ is defined, for $y=(x,v)^{\mathsf T} \in \mathbb R^{2d}$, by
\[
\Sigma(y)
=
\begin{pmatrix}
0_{d\times d} & v\\
\sqrt{\sigma}\,I_d & -v-\grad U(x)-\grad G(v)
\end{pmatrix}.
\]
Since the map $\grad U$ is smooth on $\mathbb R^d$ by Assumption~\ref{cond:U} and $\grad G$ is smooth, and hence locally Lipschitz on $\mathbb R^d$ by Assumption~\ref{cond:G:poly}, $\Sigma$ is locally Lipschitz on $\mathbb R^{2d}$. Therefore,~\cite[Thm.~1.1.8, p.~12]{Hsu} yields a unique solution to~\eqref{eq:Langevin-SDE} up to its explosion time.

Suppose now that $G$ is given by Assumption~\ref{cond:G:singular}. In this case,
\[
\X
=
\begin{cases}
\mathbb R^d\times(\mathbb R^d\setminus\{0\}), & d\ge2,\\
\mathbb R\times(0,\infty), & d=1,
\end{cases}
\]
which is an open subset of $\mathbb R^{2d}$ and hence a smooth manifold without boundary. Set $Y(t) := (x(t), v(t))^{\mathsf T} \in \X$ and $Z(t) := (W^1(t),\ldots,W^d(t) , t)^{\mathsf T} \in \mathbb R^{d+1}$, where $W^1,\ldots,W^d$ are the components of the $d$-dimensional Brownian motion $W$. Let $e_1, \ldots, e_d$ denote the standard basis of $\mathbb R^d$. For $y = (x, v)^{\mathsf T} \in \X$, define the vector fields $V_i : \X \to \mathbb R^{2d}$ by
\[
V_i(y)
:=
\begin{cases}
(0,\sqrt{\sigma}\,e_i)^{\mathsf T}, & i=1,\ldots,d,\\[1mm]
\left( v, -v-\grad U(x)-\grad G(v) \right)^{\mathsf T}, & i=d+1.
\end{cases}
\]
By Assumptions~\ref{cond:U} and~\ref{cond:G:singular}, these vector fields are smooth on $\X$. With these definitions,~\eqref{eq:Langevin-SDE} can be written as
\begin{align}
\d Y(t)
&=
\begin{pmatrix}
v(t)\,\d t\\
\left( -v(t)-\grad U(x(t))-\grad G(v(t)) \right)\d t + \sqrt{\sigma}\,\d W(t)
\end{pmatrix} \nonumber\\
&=
\sum_{i=1}^d
\begin{pmatrix}
0\\
\sqrt{\sigma}\,e_i
\end{pmatrix}
\d W^i(t)
+
\begin{pmatrix}
v(t)\\
-v(t)-\grad U(x(t))-\grad G(v(t))
\end{pmatrix}
\d t \nonumber\\
&=
\sum_{i=1}^d V_i(Y(t))\,\d W^i(t) + V_{d+1}(Y(t))\,\d t. \nonumber
\end{align}
The second equality follows from $\d W(t)=\sum_{i=1}^d e_i\,\d W^i(t)$, and the last equality follows from the definitions of $V_1,\ldots,V_{d+1}$. By the It\^o--Stratonovich conversion formula in~\cite[p.~19]{Hsu},
\begin{align}
\sum_{i=1}^d V_i(Y(t))\circ\d W^i(t)
&=
\sum_{i=1}^d V_i(Y(t))\,\d W^i(t) + \frac{1}{2} \sum_{i=1}^d D V_i(Y(t))V_i(Y(t))\,\d t \nonumber\\
&=
\sum_{i=1}^d V_i(Y(t))\, \d W^i(t), \nonumber
\end{align}
where $D V_i(y)$ denotes the Jacobian matrix of the map $V_i : \X \to \mathbb R^{2d}$ at $y\in\X$, and the second equality follows because each $V_i$, $i=1,\ldots,d$, is constant on $ \X $, and hence $D V_i=0$. Consequently,
\begin{align}
\d Y(t)
=
\sum_{i=1}^d V_i(Y(t))\circ\d W^i(t) + V_{d+1}(Y(t))\,\d t
=
\sum_{i=1}^{d+1} V_i(Y(t))\circ\d Z^i(t), \nonumber
\end{align}
where the second equality uses the fact that integration with respect to the finite-variation process $t$ is the same in the It\^o and Stratonovich senses. Thus, \eqref{eq:Langevin-SDE} is precisely a Stratonovich stochastic differential equation on the smooth manifold $\X$ driven by the $\mathbb R^{d+1}$-valued semimartingale $Z$. Using the notion of a solution of a stochastic differential equation on a manifold up to a stopping time in~\cite[Def.~1.2.3, p.~20]{Hsu},~\cite[Thm.~1.2.9, p.~23]{Hsu} yields a unique solution taking values in $\X$ up to its explosion time.
\end{proof}

Observe that, for any fixed $(x_0,v_0)\in\X$, if $\{(x(t), v(t))\}_{t \in [0,T)}$ denotes the unique solution given by Lemma~\ref{Lem:wellposed:existence-up-to-explosion}, with explosion time $T$, then, with $y(t)=(x(t), v(t))^{\mathsf T}$, the explosion time is given by $T = \lim_{n \to \infty}T_{O_n}$ almost surely for any exhaustion $\{O_n\}_{n\ge1}$ of $\X$, as in~\eqref{eq:General-Explosion-Time}. The remainder of this section is devoted to determining whether $T$ is finite or infinite in the different damping regimes.

\subsection{Polynomial damping} \label{sec:wellposed:poly-damping}

First, we consider the case of polynomial damping, i.e., when $G$ satisfies Assumption \ref{cond:G:poly}. For notational convenience, we denote by $L$ the generator associated with \eqref{eq:Langevin-SDE} given by
\begin{align} \label{form:L}
    L V = \la \grad_xV, v\ra+ \la \grad_v V, -v-\grad U(x)-\grad G(v)\ra+\frac{1}{2}\sigma\triangle_v V,
\end{align}
that is defined for $V\in C^2(\X;\rbb)$. In order to show that the explosion time $T$ obtained in Lemma \ref{Lem:wellposed:existence-up-to-explosion} is infinite a.s., it suffices to demonstrate the existence of a function $V:\X\to[0,\infty)$ satisfying that $V(x,v)\to\infty$ whenever $(x,v)\to \partial \X$ and that
\begin{align} \label{ineq:LV<cV+C}
    LV \le cV+V.
\end{align}
In Lemma \ref{lem:well-posed:poly-damping}, stated and proven next, we show that such a candidate for \eqref{eq:Langevin-SDE} turns out to be the well-known function $U(x)+\frac{1}{2}|v|^2$, which is the Hamiltonian structure for the Langevin equation in the absence of nonlinear damping, i.e., $G\equiv 0$ in \eqref{eq:Langevin-SDE}.

\begin{lemma} \label{lem:well-posed:poly-damping}
Suppose that $U$ satisfies Assumption \ref{cond:U} and $G$ is given by Assumption \ref{cond:G:poly}. Let $V_1$ be defined as 
\begin{align}\label{Def:V-One}
    V_1(x,v) = U(x)+\frac{1}{2}|v|^2.
\end{align}
Then, $V_1$ satisfies
\begin{align} \label{ineq:well-posed:poly-damping}
    LV_1(x,v) \le cV_1(x,v)+C,\quad (x,v)\in \X=\rbb^d\times \D.
\end{align}

\end{lemma}

    We remark that the result of Lemma \ref{lem:well-posed:singular-damping} is not new as it is essentially a consequence of \cite[Theorem 2.1]{markus1993stochastic}. Nevertheless, we opt for presenting its statement and proof here for the sake of completeness.

\begin{proof}[Proof of Lemma~\ref{lem:well-posed:poly-damping}]
By the definition of $V_1$ in~\eqref{Def:V-One}, we have $\grad_xV_1=\grad U(x)$, $\grad_vV_1=v$, and $\triangle_vV_1=d$.
Substituting these identities into~\eqref{form:L}, we obtain the first equality below.
\begin{align*}
L V_1(x,v)
=
\la \grad U(x),v\ra + \la v, -v-\grad U(x)-\grad G(v) \ra + \frac{\sigma d}{2}.
\end{align*}
Recalling Assumption \ref{cond:G:poly}, we have
\begin{align*}
    \la v,-\grad G(v)\ra  = -a_4|v|^\alpha-\la v,\grad G(v)-a_4v|v|^{\alpha-2}\ra \le -a_4|v|^\alpha+c|v|^{\alphat}+c|v|.
\end{align*}
Since $1\le \alphat<\alpha$ by virtue of Assumption \ref{cond:G:poly}, we may absorb $|v|^{\alphat}+|v|$ into $-a_4|v|^\alpha$, implying
\begin{align*}
 L V_1(x,v)\le - |v|^2 - \frac{a_4}{2}|v|^{\alpha} +C.
\end{align*}
This proves estimate \eqref{ineq:well-posed:poly-damping}, as claimed.
\end{proof}

\subsection{Singular damping} \label{sec:wellposed:singular-damping}

Next, we turn to the case of singular damping, i.e., when Assumption \ref{cond:G:singular} holds for $G$. As it turns out, the choice of a Lyapunov function satisfying \eqref{ineq:LV<cV+C} is much more non trivial than the usual function $U(x)+\frac{1}{2}|v|^2$. In fact, we will construct different energy estimates depending on the relation between the dimension $d$ and the parameters $\beta,a_4$, where $\beta,a_4$ are as in Assumption \ref{cond:G:singular}. This result is sumarized below in Lemma \ref{lem:well-posed:singular-damping}.

\begin{lemma} \label{lem:well-posed:singular-damping}
Suppose that $U$ satisfies Assumption \ref{cond:U} and $G$ is given by Assumption \ref{cond:G:singular}. Let $V_2$ be defined as
\begin{align} \label{form:V_2:well-posed}
    V_2 (x,v) =
\begin{cases}
J(x,v)^2+ |v|^{-\beta_1},
& \beta>1,\\
J(x,v)^2-\log|v|+c_0,
& \beta=1,\ d\ge2,\\
J(x,v)^2-\log v+c_0,
& \beta=1,\ d=1,\ 0<\sigma<2a_4,\\
J(x,v)^2-\log v- v^{\beta_2}+c_0,
& \beta=1,\ d=1,\ \sigma=2a_4,
\end{cases}
\end{align}
where $J(x,v):=U(x)+\frac{1}{2}|v|^2$, $c_0$ is a sufficiently large constant, and $\beta_1$ and $\beta_2$ are constants satisfying
\[
\max\{\beta-3,0\} < \beta_1 <\beta\Big(\frac{\lambda+1}{\lambda-1}\Big)-1, \quad \textup{and} \quad 0 < \beta_2 <\min\Big\{\frac{2}{\lambda+1},1-\betat\Big\}.
\]
Then, $V_2$ satisfies
\begin{align}\label{ineq:well-posed:singular-damping}
    LV_2(x,v) \le cV_2(x,v)+C,\quad (x,v)\in \X=\rbb^d\times \D.
\end{align}

\end{lemma}

\begin{proof} 
Since $J(x,v)=U(x)+\frac{1}{2}|v|^2$, it follows from~\eqref{form:L} that
\begin{align} \label{eq:LJ:singular}
LJ(x,v)
&=
-|v|^2
-
\la v,\grad G(v)\ra
+
\frac{\sigma d}{2} \notag 
\\
&=
-|v|^2
+
\frac{a_4}{|v|^{\beta-1}}- \Big\la v,\grad G(v)+a_4\frac{v}{|v|^{\beta+1}} \Big\ra
+
\frac{\sigma d}{2}, 
\end{align}
where the second equality is derived in accordance with
Assumption~\ref{cond:G:singular}, condition \eqref{cond:G:singular:|gradG|<1/|v|^beta}. Using this identity and
$\grad_vJ=v$, the generator chain rule yields
\begin{align} \label{eq:LJsqure}
&LJ(x,v)^2  \nonumber\\
&=
2J(x,v)LJ(x,v)+\sigma|v|^2
\nonumber\\
&=
-2J(x,v)|v|^2
+
2a_4\frac{ J(x,v)}{|v|^{\beta-1}} -2J(x,v) \Big\la v,\grad G(v)+a_4\frac{v}{|v|^{\beta+1}} \Big\ra
+
\sigma d\, J(x,v)
+
\sigma|v|^2.
\end{align}
Now there are three cases depending on the parameters $\beta$, $\sigma$ and dimension $d$. \vspace{.2cm}

Case 1: $\beta>1$. Since
\begin{align*}
\grad_v \bigg(\frac{1}{|v|^{\beta_1}}\bigg)
&=
-\frac{\beta_1v}{|v|^{\beta_1+2}},\\
\text{and }\quad \triangle_v \bigg(\frac{1}{|v|^{\beta_1}} \bigg)
&=
\operatorname{Tr}\left(\frac{\beta_1(\beta_1+2)\,v\otimes v}{|v|^{\beta_1+4}}-\frac{\beta_1 I_d}{|v|^{\beta_1+2}}\right)
=
\frac{\beta_1(\beta_1+2-d)}{|v|^{\beta_1+2}}, \nonumber
\end{align*}
using~\eqref{form:L}, we obtain
\begin{align}
&L\bigg(\frac{1}{|v|^{\beta_1}}\bigg)\\
&=
\Big\la -\frac{\beta_1v}{|v|^{\beta_1+2}}, -v-\grad U(x)-\grad G(v) \Big\ra
+
\frac{\sigma\beta_1(\beta_1+2-d)}{2|v|^{\beta_1+2}} \nonumber\\
&=
\frac{\beta_1}{|v|^{\beta_1}}
+
\frac{\beta_1\la v,\grad U(x)\ra}{|v|^{\beta_1+2}}
-
\frac{\beta_1a_4}{|v|^{\beta_1+\beta+1}} + \Big\la \frac{\beta_1v}{|v|^{\beta_1+2}}, \grad G(v) +a_4\frac{v}{|v|^{\beta+1}} \Big\ra
+
\frac{\sigma\beta_1(\beta_1+2-d)}{2|v|^{\beta_1+2}}, \nonumber
\end{align}
where the second equality is derived in accordance with
Assumption~\ref{cond:G:singular}. Since $V_2(x,v)=J(x,v)^2+ |v|^{-\beta_1}$, combining the preceding identity with~\eqref{eq:LJsqure} gives
\begin{align}\label{eq:LVTwo-BetaGreaterOne}
LV_2(x,v)
&=
-2J(x,v)|v|^2
+
2a_4\frac{J(x,v)}{|v|^{\beta-1}} -2J(x,v) \Big\la v,\grad G(v)+a_4\frac{v}{|v|^{\beta+1}} \Big\ra
+
\sigma d\,J(x,v)
+
\sigma|v|^2 \nonumber\\
&\quad
+
\frac{\beta_1}{|v|^{\beta_1}}
+
\frac{\beta_1\la v,\grad U(x)\ra}{|v|^{\beta_1+2}}
-
\frac{\beta_1a_4}{|v|^{\beta_1+\beta+1}}
 \nonumber\\
& \quad+
\frac{\sigma\beta_1(\beta_1+2-d)}{2|v|^{\beta_1+2}}+  \Big\la \frac{\beta_1v}{|v|^{\beta_1+2}}, \grad G(v) +a_4\frac{v}{|v|^{\beta+1}} \Big\ra.
\end{align}
Since $\beta_1>\beta-3$, Young's inequality and condition \eqref{cond:G:singular:|gradG|<1/|v|^beta} imply that, for every
$\delta>0$,
\begin{align}\label{eq:YoungJV-Singular}
a_4\frac{J(x,v)}{|v|^{\beta-1}}-J(x,v) \Big\la v,\grad G(v)+a_4\frac{v}{|v|^{\beta+1}} \Big\ra &\le  
a_4\frac{J(x,v)}{|v|^{\beta-1}}+ c\frac{J(x,v)}{|v|^{\betat-1}}+cJ(x,v)|v|
\notag \\
&\le
\delta J(x,v)^2
+
\frac{\delta}{|v|^{\beta_1+\beta+1}}
+
C_\delta.
\end{align}
Furthermore, since $\beta_1<\beta\big(\frac{\lambda+1}{\lambda-1}\big)-1$, Young's inequality yields the second inequality below:
\begin{align} \label{eq:YoungGradientU-Singular}
\frac{\la v,\grad U(x)\ra}{|v|^{\beta_1+2}}
\le
\frac{|\grad U(x)|}{|v|^{\beta_1+1}}
&\le
\frac{\delta}{|v|^{\beta_1+\beta+1}}
+
C_\delta \left(1+|\grad U(x)|^{\frac{2\lambda}{\lambda-1}}\right) \nonumber\\
&\le
\frac{\delta}{|v|^{\beta_1+\beta+1}}
+
C_\delta \left(1+J(x,v)^2\right),
\end{align}
where the final inequality follows from~\eqref{cond:U:U(x)=O(x^lambda+1)}--\eqref{cond:U:U'(x)=O(x^lambda)} and $U(x)^2\le J(x,v)^2$, after enlarging $C_\delta$ if necessary. For the last two terms in~\eqref{eq:LVTwo-BetaGreaterOne}, since $\beta>1$ and $0\le \betat<\beta$, Young's inequality and condition \eqref{cond:G:singular:|gradG|<1/|v|^beta} once again imply that, for every $\delta>0$,
\begin{align} \label{eq:YoungLaplacian-Singular}
&\frac{1}{|v|^{\beta_1+2}} + \Big\la \frac{\beta_1v}{|v|^{\beta_1+2}}, \grad G(v) +a_4\frac{v}{|v|^{\beta+1}} \Big\ra \notag \\
& \le \frac{C}{|v|^{\beta_1+2}}+ \frac{C}{|v|^{\beta_1+\betat+1}} + \frac{C}{|v|^{\beta_1+1}} \le 
\frac{\delta}{|v|^{\beta_1+\beta+1}}
+
C_\delta.
\end{align}
The two nonsingular positive terms in~\eqref{eq:LVTwo-BetaGreaterOne} are bounded above by $C(1+J^2)$ since $J\le1+J^2$ and $|v|^2\le2J$, while $\frac{\beta_1}{|v|^{\beta_1}}$ is controlled directly by $V_2$. Applying this observation and~\eqref{eq:YoungJV-Singular}--\eqref{eq:YoungLaplacian-Singular} in~\eqref{eq:LVTwo-BetaGreaterOne}, we obtain
\begin{align}
LV_2(x,v)
&\le
C_\delta\left(1+J(x,v)^2\right)
+
\frac{\beta_1}{|v|^{\beta_1}}
-
\frac{\beta_1a_4-\widetilde C\delta}{|v|^{\beta_1+\beta+1}}, \nonumber
\end{align}
where $\widetilde C>0$ does not depend on $\delta$ and absorbs all fixed coefficients, including
$\frac{1}{2}\sigma\beta_1|\beta_1+2-d|$. Choosing $\delta>0$ sufficiently small so that $\widetilde C\delta=\frac{1}{2}\beta_1a_4$ and discarding the resulting negative term, the desired estimate \eqref{ineq:well-posed:singular-damping} follows for suitable constants $c,C>0$. \vspace{.2cm}

Case 2: $\beta=1$ and either $d\ge2$ or $d=1$ with $0<\sigma<2a_4$. In this case, we recall that
\begin{align*}
   V_2(x,v)=J(x,v)^2-\log|v|+c_0. 
\end{align*}
A direct computation gives
\begin{align}
\grad_v\bigl(-\log|v|\bigr)
=
-\frac{v}{|v|^2},
\qquad
\triangle_v\bigl(-\log|v|\bigr)
=
\operatorname{Tr}\left(\frac{2\,v\otimes v}{|v|^4} -\frac{I_d}{|v|^2} \right)
=
-\frac{d-2}{|v|^2}. \nonumber
\end{align}
Thus, using~\eqref{form:L} and Assumption~\ref{cond:G:singular} with $\beta=1$, we obtain
\begin{align}  \label{eq:L(-log|v|)}
L\bigl(-\log|v|\bigr)
&=
\bigg\la -\frac{v}{|v|^2}, -v-\grad U(x)-\grad G(v) \bigg\ra
-
\frac{\sigma(d-2)}{2|v|^2} \\
&=
1 + \frac{\la v,\grad U(x)\ra}{|v|^2} - \frac{ \left(a_4+ \frac{\sigma}{2}(d-2) \right)}{|v|^2}+\Big\la \frac{v}{|v|^2}, \grad G(v)+a_4\frac{v}{|v|^2} \Big\ra. \nonumber
\end{align}
The above together with~\eqref{eq:LJsqure} for $\beta=1$ and $Lc_0=0$, gives
\begin{align}\label{eq:LVTwo-BetaOne}
LV_2(x,v)
&=
-2J(x,v)|v|^2 + (2a_4+\sigma d)J(x,v) -2J(x,v) \Big\la v,\grad G(v)+a_4\frac{v}{|v|^{2}} \Big\ra + \sigma|v|^2 
\notag \\
&\qquad+ 1 + \frac{\la v,\grad U(x)\ra}{|v|^2} - \frac{\left(a_4+\frac{\sigma}{2}(d-2)\right)}{|v|^2}+\Big\la \frac{v}{|v|^2}, \grad G(v)+a_4\frac{v}{|v|^2} \Big\ra.
\end{align}
As in~\eqref{eq:YoungGradientU-Singular}, Young's inequality and~\eqref{cond:U:U(x)=O(x^lambda+1)}--\eqref{cond:U:U'(x)=O(x^lambda)}
imply that, for every $\delta>0$,
\begin{align}
\frac{\la v,\grad U(x)\ra}{|v|^2}
\le
\frac{|\grad U(x)|}{|v|}
\le
\frac{\delta}{100}\frac{1}{|v|^2}
+
C_\delta|\grad U(x)|^2
\le
\frac{\delta}{100}\frac{1}{|v|^2}
+
C_\delta\left(1+J(x,v)^2\right), \nonumber
\end{align}
whereas Assumption \ref{cond:G:singular} yields
\begin{align*}
  &-2J(x,v) \Big\la v,\grad G(v)+a_4\frac{v}{|v|^2} \Big\ra+  \Big\la \frac{v}{|v|^2}, \grad G(v)+a_4\frac{v}{|v|^2} \Big\ra\\
  & \le C \frac{J(x,v)}{|v|^{\betat-1}}+CJ(x,v)|v|+ C\frac{1}{|v|^{\betat+1}}+C\frac{1}{|v|}\\
  &\le C(1+J(x,v)^2)+\frac{\delta}{100}\frac{1}{|v|^2}.
\end{align*}
Using the preceding estimate in~\eqref{eq:LVTwo-BetaOne} and bounding the nonsingular positive terms by $C(1+J^2)$, we obtain
\begin{align}
LV_2(x,v)
&\le
C_\delta\left(1+J(x,v)^2\right)
-
\frac{\left(a_4+\frac{\sigma}{2}(d-2)-\delta\right)}{|v|^2}
\le
C\left(1+J(x,v)^2\right), \nonumber
\end{align}
where the second inequality follows by choosing $\delta=\frac{1}{2}\left(a_4+\frac{\sigma}{2}(d-2)\right)$, which is
positive when either $d\ge2$ or $d=1$ with $0<\sigma < 2a_4$. Finally, since the function $f(r) :=\log r-\frac{1}{8}r^4$ is continuous on $(0, \infty)$ and $f(r) \to -\infty$ as $r\to0$ or $r\to\infty$, there exists a constant $K>0$ such that the first inequality below holds:
\[
\log|v| \le \frac{1}{8}|v|^4 + K \le \frac12J(x,v)^2 + K.
\]
The second inequality follows from $|v|^2\le 2J(x,v)$. Therefore, from the definition of $V_2$ given by~\eqref{form:V_2:well-posed},
\[
V_2(x,v) \ge \frac12J(x,v)^2+c_0-K.
\]
Thus, choosing $c_0\ge K+1$, we obtain $V_2(x,v)\ge1$ and $J(x,v)^2 \le 2 V_2(x,v)$. Substituting this into the preceding generator estimate yields the desired estimate for suitable constants $c, C > 0 $. \vspace{1cm}

Case 3: $\beta=1$, $d=1$, and $\sigma=2a_4$.
In this case, we note that $v>0$ and recall from \eqref{form:V_2:well-posed} that
\begin{align*}
   V_2(x,v)= J(x,v)^2-\log v- v^{\beta_2}+c_0.
\end{align*}
Using~\eqref{form:L} and
Assumption~\ref{cond:G:singular} with $\beta=1$, we obtain
\begin{align}
L\bigl(-v^{\beta_2}\bigr)
&=
-\beta_2v^{\beta_2-1} \left( -v-\grad U(x)-\grad G(v) \right)
+
\beta_2(1-\beta_2)a_4v^{\beta_2-2} \nonumber\\
&=
\beta_2v^{\beta_2}
+
\frac{\beta_2\grad U(x)}{v^{1-\beta_2}}
-
\frac{\beta_2^2a_4}{v^{2-\beta_2}}+\beta_2v^{\beta_2-1}\Big(\grad G(v)+a_4\frac{1}{v^{\beta}}\Big). \nonumber
\end{align}
Combining~\eqref{eq:LVTwo-BetaOne} for $\beta=1$, $d=1$, and $\sigma=2a_4$ with the preceding identity gives the second equality below.
\begin{align}\label{eq:LVTwo-BetaOne-dOne-SigmaTwo}
LV_2(x,v)
&=
L\bigl(J(x,v)^2-\log v\bigr)
+
L\bigl(-v^{\beta_2}\bigr) \nonumber\\
&=
-2J(x,v)v^2
+
4a_4J(x,v) -2J(x,v) v\Big(\grad G(v)+\frac{a_4}{v}\Big)
+
2a_4v^2 \notag \\
&\qquad
+
1
+
\frac{\grad U(x)}{v} + \frac{1}{v}\Big(\grad G(v)+\frac{a_4}{v}\Big) \notag \\
&\qquad
+
\beta_2v^{\beta_2}
+
\frac{\beta_2\grad U(x)}{v^{1-\beta_2}}
-
\frac{\beta_2^2a_4}{v^{2-\beta_2}}+\beta_2v^{\beta_2-1}\Big(\grad G(v)+a_4\frac{1}{v^{\beta}}\Big).
\end{align}
For the singular terms involving $\grad U(x)$, observe that $2-\beta_2$ and $\frac{2-\beta_2}{1-\beta_2}$ are conjugate exponents, since $0<\beta_2<\frac{2}{\lambda+1}< 1$. Therefore, Young's inequality gives, for every $\delta>0$,
\begin{align*}
\frac{|\grad U(x)|}{v}+ \beta_2\frac{|\grad U(x)|}{v^{1-\beta_2}} &\le C\frac{|\grad U(x)|}{v} +C|\grad U(x)| \\
&\le
\frac{1}{100}\frac{\delta}{v^{2-\beta_2}}
+
C_\delta |\grad U(x)|^{\frac{2-\beta_2}{1-\beta_2}}+C_\delta. \nonumber
\end{align*}
Next, using~\eqref{cond:U:U'(x)=O(x^lambda)} and enlarging the
constant if necessary, we obtain the first inequality below.
\begin{align*}
|\grad U(x)|^{\frac{2-\beta_2}{1-\beta_2}}
&\le
C\left( 1+|x|^{(\lambda-1)\frac{2-\beta_2}{1-\beta_2}} \right)
\le
C\left( 1+|x|^{2\lambda} \right)
\le
C\left( 1+J(x,v)^2 \right)
\end{align*}
The second to last inequality uses $\frac{2-\beta_2}{1-\beta_2}<\frac{2\lambda}{\lambda-1}$, which follows from $\beta_2<\frac{2}{\lambda+1}$, whereas the last inequality follows from the lower bound in~\eqref{cond:U:U(x)=O(x^lambda+1)} and the definition of $J$. Combining the preceding estimates, we obtain
\begin{align} \label{eq:YoungGradientU-SigmaTwo}
\frac{|\grad U(x)|}{v}+ \beta_2\frac{|\grad U(x)|}{v^{1-\beta_2}}
\le
\frac{1}{100}\frac{\delta}{v^{2-\beta_2}}
+
C_\delta\left( 1+J(x,v)^2 \right).
\end{align}
Next, we employ condition \eqref{cond:G:singular:|gradG|<1/|v|^beta} to see that
\begin{align*}
    -2J(x,v) v\Big(\grad G(v)+\frac{a_4}{v}\Big) \le 2J(x,v) v\cdot |a_5|\Big(\frac{1}{v^{\betat}}+1\Big)&\le 2|a_5| J(x,v)\big(v^{1-\betat}+v\big)\\
    &\le CJ(x,v)(v+1)\\
    &\le C(1+J(x,v)^2).
\end{align*}
Also, we employ the condition $\beta_2<1-\betat$ to infer
\begin{align*}
     & \frac{1}{v}\Big(\grad G(v)+\frac{a_4}{v}\Big) +\beta_2v^{\beta_2-1}\Big(\grad G(v)+a_4\frac{1}{v^{\beta}}\Big) \\
    &\le |a_5|\Big( \frac{1}{v^{1+\betat}}+\frac{1}{v}+\beta_2\frac{1}{v^{1+\betat-\beta_2}}+\frac{1}{v^{1-\beta_2}}  \Big) \le \frac{1}{100} \frac{\delta}{v^{2-\beta_2}}+C_\delta.
\end{align*}
Using the above with~\eqref{eq:YoungGradientU-SigmaTwo}
in~\eqref{eq:LVTwo-BetaOne-dOne-SigmaTwo}, and bounding the nonsingular positive terms by $C(1+J^2)$, we obtain
\begin{align}
LV_2(x,v)
\le
C_\delta\left( 1+J(x,v)^2 \right)
-
\frac{  \beta_2^2a_4-\delta  }{ v^{2-\beta_2} }
\le
C\left( 1+J(x,v)^2 \right), \nonumber
\end{align}
where the second inequality follows by choosing $\delta = \frac{\beta_2^2 a_4}{2}$. Since $0<\beta_2<1$, we have the first inequality below:
\[
\log v+v^{\beta_2} \le \frac18v^4+K \le \frac12J(x,v)^2+K,
\]
where, as in Case~2, $K>0$ is a sufficiently large constant, since
$f(r):=\log r+r^{\beta_2}-\frac18r^4$ is continuous and bounded above on $(0,\infty)$. Repeating the final argument of Case~2, we obtain $V_2\ge1$ and the desired estimate
\[
LV_2(x,v)\le cV_2(x,v)+C
\]
for suitable constants $c,C>0$.
\end{proof}

\begin{lemma}\label{lem:well-posed:singular-damping:sigma>2:illposed}
Suppose that $U$ satisfies Assumption~\ref{cond:U} and $G$ is defined in Assumption~\ref{cond:G:singular}. Let $d=1$, $\beta=1$, and $\sigma > 2a_4$, and fix $(x_0,v_0)\in\X=\mathbb R \times (0, \infty)$. Let $\{(x(t), v(t))\}_{t \in [0,T)}$ be the unique solution given by Lemma~\ref{Lem:wellposed:existence-up-to-explosion}, with explosion time $T$. Then the explosion time is finite almost surely and the explosion occurs through the velocity boundary $0$; more precisely,
\[
T<\infty,
\qquad
\lim_{t\uparrow T}v(t)=0,
\qquad
\sup_{0\le t<T}\bigl(|x(t)|+v(t)\bigr)<\infty
\qquad
\textup{a.s.}
\]
\end{lemma}
\begin{proof}[Proof of Lemma~\ref{lem:well-posed:singular-damping:sigma>2:illposed}]
Fix $(x_0,v_0)\in\X=\mathbb R\times(0,\infty)$, and let
$\{(x(t),v(t))\}_{t\in[0,T)}$ be the unique solution given by
Lemma~\ref{Lem:wellposed:existence-up-to-explosion}, with explosion time $T$. Since $d=1$ and $\beta=1$, by Assumption~\ref{cond:G:singular}, condition \eqref{cond:G:singular:|gradG|<1/|v|^beta}, there exists a function
$r:(0,\infty)\to\mathbb R$ such that
\[
G'(v) = -\frac{a_4}{v}+r(v), \qquad |r(v)|\le |a_5|\left(v^{-\betat } + 1 \right), \qquad v>0,
\]
where $0 \le \betat < 1$. Define the Lyapunov function $J(x,v) := U(x)+\frac{1}{2} v^2$ for $(x,v) \in \mathbb R \times(0,\infty)$. Then, using the generator in~\eqref{form:L}, we obtain\begin{align}\label{eq:Generator-Energy-Singular}
LJ(x,v)
&= vU'(x) + v \left(-v- U'(x) + \frac{a_4}{v}-r(v) \right)+ \frac{\sigma}{2} \nonumber\\
&=-v^2+a_4-vr(v) + \frac{\sigma}{2}\nonumber\\
&\le-v^2+a_4+|a_5| \left(v^{1 -\betat } +
v\right)+
\frac{\sigma}{2} \le C.
\end{align}
Here, the final inequality follows from Young's inequality, since
$0< 1- \betat \le1$, which gives, for every $\varepsilon>0$,
$v^{1-\betat} \le \varepsilon v^2+C_\varepsilon$ and $v\le \varepsilon v^2+C_\varepsilon$; the constants $a_4$ and $ \sigma/2$ are absorbed into $C_\varepsilon$. Choosing $ \varepsilon >0$ sufficiently small so that $2 |a_5| \varepsilon \le 1/2$ yields the final inequality for some constant $C>0$ independent of $(x,v)$. Choose $m_0>\max\{|x_0|,v_0,v_0^{-1},1\}$ and, for each $m\ge m_0$, set
\[
O_m
:=
(-m,m)\times\left(\frac1m,m\right).
\]
Then $\{O_m\}_{m\ge m_0}$ is an exhaustion of $\X = \mathbb R \times (0,\infty)$ and, with $y(t) := (x(t),v(t))^{\mathsf T}$, the corresponding exit times $T_m:= T_{O_m}$ defined in~\eqref{eq:General-Exit-Times} are given by
\begin{align}\label{eq:Singular-Damping-Exit-Times}
T_m
=
\inf\left\{ 0\le t<T: |x(t)|\ge m \ \text{or}\ v(t) \notin \left( \frac1m,m \right) \right\} \wedge T.
\end{align}
Moreover, $T_m\uparrow T$ almost surely as $m\to\infty$ by~\eqref{eq:General-Explosion-Time}. Also, $T_m<T$ almost surely on $\{T<\infty\}$ for every $m\ge m_0$. Indeed, since $\overline{O_m}$ is a compact subset of $\X$, the set $\X^c \setminus \overline{O_m}$ is an open neighborhood of
$\partial_{\X}$ in the one-point compactification $\X^c$. On $\{T<\infty\}$, we have $y(T)=\partial_{\X}$, and the continuity of the path $y$ therefore implies that there exists some $s<T$ such that $y(s) \notin \overline{O_m}$ and hence $y(s)\notin O_m$. By the definition of $T_m$, it follows that $T_m\le s<T$.

Choose $n_0>\max\{v_0^{-1},1\}$. For $n\ge n_0$, define the stopping times corresponding to exit through the lower velocity levels by
\[
T_n^-
:=
\inf\left\{ 0\le t<T: v(t)\le\frac1n \right\}, \qquad 
T^-
:=
\lim_{n\to\infty}T_n^-,
\]
with the convention that $\inf\varnothing=\infty$. Since the levels $1/n$ decrease with $n$, the sequence $\{T_n^-\}_{n\ge n_0}$ is nondecreasing, and hence the limit $T^-$ is well defined. By continuity of $v$, whenever $T_n^-<\infty$, we have $v(T_n^-)=1/n$.

We first show that, for every $n\ge n_0$, if the explosion time $T$ is finite, then the solution reaches the level $1/n$ strictly before $T$. Fix $n\ge n_0$ and, for $m \ge \max\{m_0,n\}$, set
\[
\mathcal T_{m,n}:=T_m\wedge T_n^-.
\]
For every $t>0$, the solution is defined up to $t \wedge \mathcal T_{m,n}$, since $\mathcal T_{m,n}\le T_m<T$ almost surely on $\{T<\infty\}$, while this is immediate on $\{T = \infty\}$. Applying It\^o's formula to $J(x(t), v(t))$ up to this stopping time and using the second equality in~\eqref{eq:Generator-Energy-Singular}, we obtain
\begin{align}
J\bigl(x(t \wedge \mathcal T_{m,n}),v(t \wedge \mathcal T_{m,n})\bigr)
&= J(x_0,v_0) + \int_0^{t\wedge\mathcal T_{m,n}} \left(-v(s)^2+a_4-v(s)r(v(s)) +
\frac{\sigma}{2}\right)\d s  \nonumber\\ 
& \quad+ \sqrt{\sigma} \int_0^{t \wedge \mathcal T_{m,n}} v(s)\,\d W(s). \nonumber
\end{align}
Also, by It\^o's isometry,
\[
\mathbb E\left[ \left| \int_0^{t\wedge\mathcal T_{m,n}} v(s)\,\d W(s) \right|^2 \right] = \mathbb E\left[ \int_0^t \mathbf 1_{\{s<\mathcal T_{m,n}\}} v(s)^2\,\d s \right]
\le m^2t < \infty,
\]
where the inequality follows since $s<\mathcal T_{m,n}\le T_m$ implies $s<T_m$; by the definition of $T_m$ as the first exit time from $O_m=(-m,m)\times(1/m,m)$, we have $(x(s),v(s))\in O_m$ for every $s<\mathcal T_{m,n}$, and hence $v(s)<m$. Thus,
$\int_0^{t\wedge\mathcal T_{m,n}}v(s)\,\d W(s)$ is a square-integrable martingale and has expectation zero. Taking expectations and using~\eqref{eq:Generator-Energy-Singular}, we obtain
\begin{align}\label{eq:Stopped-Energy-Estimate}
\mathbb E\left[ J\bigl(x(t\wedge\mathcal T_{m,n}) ,v(t\wedge \mathcal T_{m,n})\bigr) \right] \le J(x_0,v_0) + C\,t.
\end{align}

Fix $t>0$. On the event $\{T_m\le t,\ T_m<T_n^-\}$, we cannot have $v(T_m) =1/m$. Indeed, since $m\ge n$, $1/m\le1/n<v_0$, and the continuity of $v$ would then imply that $v(r)=1/n$ for some $r\le T_m$, contradicting $T_m<T_n^-$. Therefore, by the definition~\eqref{eq:Singular-Damping-Exit-Times} of $T_m$ and
the continuity of $x$ and $v$, either $|x(T_m)|=m$ or $v(T_m)=m$. In the first case, we have the second inequality below by~\eqref{cond:U:U(x)=O(x^lambda+1)}
\[
J\bigl(x(T_m),v(T_m)\bigr)
=
U\bigl(x(T_m)\bigr) + \frac{1}{2} v(T_m)^2
\ge
U\bigl(x(T_m)\bigr)
\ge
\frac1{a_1}m^\lambda-a_1,
\]
whereas in the second case, since $U\ge0$,
\[
J\bigl(x(T_m),v(T_m)\bigr)
\ge
\frac{1}{2}v(T_m)^2
=
\frac{1}{2}m^2.
\]
Hence, on $\{T_m\le t,\ T_m<T_n^-\}$, we have
\[
J\bigl(x(T_m),v(T_m)\bigr)
\ge
\min\left\{ \frac1{a_1}m^\lambda-a_1, \frac{1}{2}m^2 \right\}
=:c_m.
\]
Since $\lambda\ge2$, we have $c_m\to\infty$ as $m\to\infty$, and
$c_m>0$ for all sufficiently large $m$. Moreover, on $\{T_m\le t,\ T_m<T_n^-\}$, $t\wedge\mathcal T_{m,n}=T_m$. Consequently,
\[
c_m
\mathbf 1_{\{T_m\le t,\ T_m<T_n^-\}}
\le
J\bigl( x(t\wedge\mathcal T_{m,n}), v(t\wedge\mathcal T_{m,n}) \bigr).
\]
Taking expectations and using~\eqref{eq:Stopped-Energy-Estimate}, we obtain
\begin{align}\label{eq:Exit-Probability-Before-Tn-Minus}
\mathbb P\left( T_m\le t,\  T_m<T_n^- \right)
\le
\frac{ J(x_0,v_0) + C\,t }{ c_m }.
\end{align}
On the event $\{T\le t,\ T\le T_n^-\}$, we have
$T_m<T\le t$ and $T_m<T\le T_n^-$ almost surely. Hence,
\[
\{T\le t,\ T\le T_n^-\}
\subseteq
\{T_m\le t,\ T_m<T_n^-\}
\qquad
\textup{a.s.}
\]
Therefore, by~\eqref{eq:Exit-Probability-Before-Tn-Minus},
\[
\mathbb P\left( T\le t,\  T\le T_n^- \right)
\le
\frac{ J(x_0,v_0) + C\, t }{ c_m }.
\]
Letting $m\to\infty$, we obtain $\mathbb P\left( T\le t,\  T\le T_n^- \right) = 0$. Since
\[
\{T<\infty,\ T\le T_n^-\}
=
\bigcup_{k=1}^{\infty} \{T\le k,\ T\le T_n^-\},
\]
we conclude that
\begin{align}\label{eq:No-Explosion-Before-Tn-Minus}
\mathbb P\left( T<\infty,\  T\le T_n^- \right)
=
0, \qquad n\ge n_0.
\end{align}

Since $v(t)>0$ for every $0\le t<T$, the first equation
in~\eqref{eq:Langevin-SDE} gives
\[
x(t)
=
x_0+\int_0^t v(s)\,\d s \ge x_0,
\qquad 0\le t<T.
\]
Also, since $d=1$, Assumption~\ref{cond:U} gives $xU'(x)\ge a_2|x|^\lambda-a_3$. Choose $R>\max\{x_0,1\}$ sufficiently large that $a_2R^\lambda-a_3\ge0$. Then, for every $x\ge R$,
\[
xU'(x)
\ge a_2x^\lambda-a_3 \ge 0,
\]
and hence $U'(x)\ge0$. Since $U'$ is continuous on the compact interval $[x_0,R]$, the constant
\[
C_{x_0}
:=
\max\Big\{ 0,\, \sup_{x_0\le x\le R}\bigl(-U'(x)\bigr) \Big\}
\]
is finite. Consequently,
\[
-U'(x)\le C_{x_0}, \qquad x\ge x_0.
\]

Since $\sigma>2a_4$, choose $\eta>0$ such that $k:=a_4+\eta<\sigma/2$. Since $0\le\betat<1$, there exists $K_\eta\ge0$ such that the second inequality below holds:
\begin{align}\label{eq:Singular-Perturbation-Comparison}
|r(v)| \le |a_5| \left(v^{-\betat} + 1\right) \le \frac{\eta}{v} + K_\eta, \qquad v>0.
\end{align}
Indeed, for $0<\betat<1$, this follows from Young's inequality applied to $(v^{-1})^{\betat}$; for $\betat=0$, the remainder is bounded by $2|a_5|$. Consequently, for $x\ge x_0$ and $v>0$,
\[
\begin{aligned}
-v-U'(x)+\frac{a_4}{v}-r(v)
&\le -v+C_{x_0} + K_\eta+\frac{a_4+\eta}{v}\\
&= -v+C_{x_0}+K_\eta+\frac{k}{v},
\end{aligned}
\]
since $k=a_4+\eta$. We therefore consider the one-dimensional stochastic differential equation with state space $(0,\infty)$, driven by the same Brownian motion $W$,
\begin{align}\label{eq:Auxiliary-Comparison-Diffusion}
\d u(t)
=
\left( -u(t)+C_{x_0}+K_\eta+\frac{k}{u(t)} \right)\d t + \sqrt{\sigma}\,\d W(t), \qquad u(0)=v_0.
\end{align}
By Lemma~\ref{Lem:Comparison-Diffusion}, applied with $C=C_{x_0}+K_\eta$ and $\sigma>2k$, there exists a unique solution $\{u(t)\}_{t\in[0,S)}$ of~\eqref{eq:Auxiliary-Comparison-Diffusion} with explosion time $S$ such that
\begin{align} \label{Eq:Lem:Comparison-Diffusion-Deductions}
S<\infty, \qquad \lim_{t\uparrow S}u(t)=0, \qquad \sup_{0\le t<S}u(t)<\infty \qquad \textup{a.s.}
\end{align}
For $n\ge n_0$, define the stopping times corresponding to the lower levels by
\[
S_n^- := \inf\left\{ 0\le t<S: u(t)\le\frac1n \right\}, \qquad S^-:=\lim_{n\to\infty}S_n^-.
\]
with the convention that $\inf\varnothing=\infty$. Then, by the continuity of $u$ and the fact that $\lim_{t\uparrow S}u(t)=0$, we have $S_n^-\uparrow S^-=S$ almost surely. Fix $n\ge n_0$ and set
\[
\eta_n:=T_n^-\wedge S_n^-.
\]
By~\eqref{eq:No-Explosion-Before-Tn-Minus}, we have $T_n^-<T$ almost surely on $\{T<\infty\}$, while $S_n^-<S$ almost surely. Hence, both $v$ and $u$ are defined on $[0,\eta_n]$ almost surely. Set
\[
D(t):=v(t)-u(t).
\]
Since $v$ and $u$ are driven by the same Brownian motion and
$D(0)=0$, subtracting~\eqref{eq:Auxiliary-Comparison-Diffusion}
from the velocity equation in~\eqref{eq:Langevin-SDE} and using $ G'(v) =-\frac{a_4}{v}+r(v)$ yields
\[
D(t\wedge\eta_n) = \int_0^{t\wedge\eta_n} \left[ -D(s)-U'(x(s))-C_{x_0}-K_\eta+\frac{a_4}{v(s)}-r(v(s))-\frac{k}{u(s)} \right]\d s.
\]
Since $D$ is absolutely continuous, the chain rule for its positive part, $(D^+)'(t)=\mathbf 1_{\{D(t)>0\}}D'(t)$ for almost every $t$, gives
\[
D^+(t\wedge\eta_n)
=
\int_0^{t\wedge\eta_n} \mathbf 1_{\{D(s)>0\}}
\left[ -D(s)-U'(x(s))-C_{x_0}-K_\eta+\frac{a_4}{v(s)}-r(v(s))-\frac{k}{u(s)} \right]\d s.
\]
By~\eqref{eq:Singular-Perturbation-Comparison}, $-U'(x(s)) \le C_{x_0}$, and recalling that $k:=a_4 + \eta$, the expression in brackets is at most
\[
-D(s) + k\left(\frac{1}{v(s)} - \frac1{u(s)} \right).
\]
This is nonpositive on $\{D(s)>0\}$, since $v(s)>u(s)$ there. Hence $D^+(t\wedge\eta_n)\le0$. Since $D^+\ge0$, we conclude that $D^+(t\wedge\eta_n)=0$, and therefore
\begin{align}\label{eq:Velocity-Comparison}
v(t)\le u(t), \qquad 0\le t\le\eta_n.
\end{align}

Also, observe that if $S_n^-<T_n^-$, then $\eta_n=S_n^-$, and~\eqref{eq:Velocity-Comparison} gives $v(S_n^-)\le u(S_n^-)=1/n$. Since $T_n^-$ is the first time that $v$ reaches the level $1/n$, this implies $T_n^-\le S_n^-$, contradicting $S_n^-<T_n^-$. Therefore, $T_n^-\le S_n^-$ almost surely for every $n\ge n_0$, which gives the first inequality below:
\[
T^- = \lim_{n\to\infty}T_n^- \le \lim_{n\to\infty}S_n^- = S < \infty \qquad \textup{a.s.},
\]
where the finiteness holds by~\eqref{Eq:Lem:Comparison-Diffusion-Deductions}. In particular, $T_n^-<\infty$ almost surely for every $n\ge n_0$. On the event $\{T<\infty\}$,~\eqref{eq:No-Explosion-Before-Tn-Minus} gives $T_n^-<T$ almost surely, while on $\{T=\infty\}$ this inequality is immediate. Hence, $T_n^-<T$ almost surely for every $n\ge n_0$, and therefore $T^-\le T$ almost surely. On the event $\{T^-<T\}$, the continuity of $v$ on $[0,T)$ and
$v(T_n^-)=1/n$ imply that $v(T^-)=\lim_{n\to\infty}v(T_n^-)=0$, which is not possible since $T^-<T$ and the solution takes values in $\X=\mathbb R\times(0,\infty)$ on $[0,T)$. Therefore,
\[
T=T^-\le S<\infty \qquad \textup{a.s.}
\]

Also, using $T_n^-\le S_n^-$ almost surely for every $n\ge n_0$, we have $\eta_n=T_n^-$. Thus,~\eqref{eq:Velocity-Comparison} implies that, for every $t<T$, choosing $n$ sufficiently large that $t<T_n^-$ gives $v(t)\le u(t)$. Therefore,
\[
\sup_{0\le t<T}v(t) \le \sup_{0\le t<S}u(t) < \infty \qquad \textup{a.s.},
\]
where the finiteness follows from~\eqref{Eq:Lem:Comparison-Diffusion-Deductions}. Since $x(t)=x_0+\int_0^t v(s)\,\d s$ and $T<\infty$ almost surely, it follows that
\[
\sup_{0\le t<T}|x(t)|
\le |x_0| + T\sup_{0\le t<T}v(t)
<
\infty \qquad \textup{a.s.}
\]
Combining the above estimates, we have
\[
\sup_{0\le t<T}\bigl(|x(t)|+v(t)\bigr)
<
\infty \qquad \textup{a.s.}
\]
The above, together with the definition of the explosion time and the fact that $T<\infty$ almost surely, implies that
\[
\lim_{t\uparrow T}v(t)=0 \qquad \textup{a.s.}
\]
\end{proof}

Finally, we provide the proof of Theorem \ref{thm:well-posed}.

\begin{proof}[Proof of Theorem \ref{thm:well-posed}]
In light of Lemma \ref{lem:well-posed:poly-damping} and Lemma \ref{lem:ergodicity:singular-damping} giving energy estimates on \eqref{eq:Langevin-SDE}, following the framework of \cite{veretennikov2024lyapunov}, we obtain the global well-posedness of \eqref{eq:Langevin-SDE} and establish Theorem \ref{thm:well-posed}, part 1 and 3. On the other hand, Lemma \ref{lem:well-posed:singular-damping:sigma>2:illposed} proves that in dimension $d=1$, when condition \eqref{cond:d-beta-sigma:ill-posed} holds, the explosion time $T$ constructed in Lemma \ref{Lem:wellposed:existence-up-to-explosion} is finite a.s. and that 
\begin{align*}
    \lim _{t\to T^{-}}v(t)=0.
\end{align*}
This concludes Theorem \ref{thm:well-posed}, part 2, as claimed.

\end{proof}

\section{Lyapunov functions for ergodicity} \label{sec:ergodicity}

In this section, we construct suitable Lyapunov functions for the analysis of ergodicity and mixing of $P_t$. In Section \ref{sec:ergodicity:poly-damping}, we focus on the derivation for the polynomial damping under Assumption \ref{cond:G:poly}. In Section \ref{sec:ergodicity:singular-damping}, we present the construction for the repulsive damping $G$ satisfying Assumption \ref{cond:G:singular} when the leading singular parameter $\beta=1$ and the additional hypothesis of Assumption \ref{cond:G:singular:beta>1} when $\beta>1$. The results of this section will be employed to establish Theorem \ref{thm:mixing} for ergodicity later in Section \ref{sec:ergodicity-mixing}.

\subsection{Polynomial damping} \label{sec:ergodicity:poly-damping}

With regard to the setting of the polynomial damping, before diving into details, we provide a heuristic derivation on the construction. For this purpose, we recall the parameters $\lambda$ and $\alpha$, respectively from Assumption \ref{cond:U} and Assumption \ref{cond:G:poly} representing the leading order terms of $U$ and $G$. We note that our construction will depend on the competition nature between these two constants. 

We first consider the case $\lambda\ge \alpha$. Heuristically, in this case, the potential effect on the velocity equation of \eqref{eq:Langevin-SDE} dominates the nonlinear effect from the damping. As a first ansatz, we may consider the classical Lyapunov approach from \cite{mattingly2002ergodicity,pavliotis2014stochastic} for the Langevin system when $G\equiv 0$. Such a function is well known and given by
\begin{align*}
    V= U(x)+\frac{1}{2}|v|^2+\varepsilon\la x,v\ra,
\end{align*}
where $\varepsilon$ is sufficiently small. As it turns out, the above function will produce the desired Lyapunov effect of the form
\begin{align*}
    LV\le -cV+C,
\end{align*}
which leads to the exponential convergence rate. This result can be found below in Lemma \ref{lem:ergodicity:poly-damping}, part (i). 

Turning to the converse scenario when $\lambda<\alpha$, we note that in this case, while $G$ has a larger magnitude compared with $U$, neither $U$ nor $G$ can be considered the leading term in \eqref{eq:Langevin-SDE} and thus in the formula of the generator $L$ as in \eqref{form:L}. So, heuristically, as $|(x,v)|\to \infty$, 
\begin{align*}
    L \sim \la \grad_v,-\grad U(x)-\grad G(v)\ra.
\end{align*}
Following the framework of the stochastic wave equations in \cite{barbu2007stochastic,nguyen2024polynomial}, we look for a function $V$ satisfying 
\begin{align*}
    \la \grad_v V,-\grad U(x)-\grad G(v)\ra \le 0.
\end{align*}
Such a solution is given by $G(v)+\la v,\grad U(x)\ra$ which produces a dissipation
\begin{align*}
    \la \grad_v (G(v)+\la v,\grad U(x)\ra),-\grad U(x)-\grad G(v)\ra = -|\grad G(v)+\grad U(x)|^2.
\end{align*}
In order to control the term $\grad U$ in $\la v,\grad U\ra$, it suffices to consider 
\begin{align*}
    V = \Big(U(x)+\frac{1}{2}|v|^2\Big)^2+C_0\big(G(v)+\la v,\grad U(x)\ra\big)+c_0,
\end{align*}
where $C_0$ and $c_0$ are suitably chosen constants. In Lemma \ref{lem:ergodicity:poly-damping}, part (ii) below, we show that the above $V$ satisfies for all $n$
\begin{align*}
    LV^n\le -cV^{n-\frac{1}{2}}+C,
\end{align*}
leading to polynomial mixing rates of arbitrary orders.

Now, we state Lemma \ref{lem:ergodicity:poly-damping} rigorously establishing the Lyapunov functions for the setting of polynomial damping.

\begin{lemma} \label{lem:ergodicity:poly-damping}
Suppose that $U$ satisfies Assumptions~\ref{cond:U} and $G$ is given by Assumption~\ref{cond:G:poly}. Let $V_3$ be defined for $(x,v)\in\X=\rbb^{d}\times \rbb^d$
\begin{align}\label{Def:V-Three}
V_3(x,v)
:=\begin{cases}
J(x,v)+\varepsilon\la x,v\ra+c_0,&\alpha\le \lambda,\\
   J(x,v)^2 + C_0 G(v) + C_0\la v,\grad U(x)\ra + c_0,& \alpha> \lambda,\\
\end{cases}
\end{align}
where $J(x,v) := U(x) + \frac{1}{2}|v|^2$, and $\varepsilon, C_0, c_0$ are positive constants. Then, for all $\varepsilon$ sufficiently small and $c_0,C_0$ sufficiently large, $V_3\ge1$ and the followings hold:
\begin{itemize}
    \item[\textup{(i)}] When $\alpha\le \lambda$, 
    \begin{align} \label{ineq:Lyapunov:poly-damping:exponential}
        L V_3(x,v) \le -cV_3(x,v)+C,\quad (x,v)\in\X,
    \end{align}
    for some positive constant $c,C$ independent of $(x,v)\in\X$.
    \item[\textup{(ii)}] When $\alpha>\lambda$, for every real $n\ge1$, there exist constants
$c_n,C_n>0$, independent of $(x,v)\in\X$, such that
\begin{align}\label{ineq:Lyapunov:poly-damping:polynomial}
    L \bigl( V_3(x,v)^n \bigr) \le - c_n V_3(x,v)^{n-\frac12} + C_n, \qquad (x,v) \in \X.
\end{align}

\end{itemize}

\end{lemma}

\begin{proof} 
(i) Considering the case $\alpha\le \lambda$, recall that $J(x,v) = U(x) + \frac{1}{2}|v|^2$. Since $\grad_xJ = \grad U(x)$, $\grad_vJ=v$, and $\triangle_vJ=d$, we have by~\eqref{form:L} and Assumption~\ref{cond:G:poly},
\begin{align}
LJ(x,v)
&=
\la \grad U(x),v\ra + \la v, -v-\grad U(x)-\grad G(v) \ra + \frac{\sigma d}{2} \notag \\
&=
-|v|^2-a_4|v|^{\alpha}-\la v,\grad G(v)-a_4v|v|^{\alpha-2}\ra+\frac{\sigma d}{2}. 
\end{align}
In view of Assumption \ref{cond:G:poly}, condition \eqref{cond:G:poly:grad.G}, we readily have
\begin{align*}
    \la v,\grad G(v)-a_4v|v|^{\alpha-2}\ra \le a_5(|v|^{\alphat}+|v|) \le \frac{1}{100}a_4|v|^\alpha+C,
\end{align*}
whence
\begin{align} \label{ineq:LJ}
    LJ(x,v)\le -|v|^2-\frac{3}{4}a_4|v|^{\alpha} +C.
\end{align}
Next, concerning $\varepsilon\la x,v\ra$, it holds that
\begin{align*}
    \varepsilon L\la x,v\ra =\varepsilon |v|^2-\varepsilon \la x,v\ra-\varepsilon\la x,\grad G(v)\ra-\varepsilon\la x,\grad U(x)\ra.
\end{align*}
In view of \eqref{cond:U:x.U'(x)>-x^(lambda+1)}, we readily have
\begin{align*}
    -\varepsilon\la x,\grad U(x)\ra \le -\varepsilon a_2|x|^\lambda+\varepsilon a_3,
\end{align*}
whereas by Young's inequality and Assumption \ref{cond:G:poly}, condition \eqref{cond:G:poly:grad.G},
\begin{align*}
    -\varepsilon\la x,\grad G(v)\ra  & = -\varepsilon \la x,v\ra|v|^{\alpha-2}-\varepsilon\la x,\grad G(v)-v|v|^{\alpha-2}\ra \\
    &\le \varepsilon |x| |v|^{\alpha-1}+\varepsilon a_5 |x|(|v|^{\alphat-1}+1)\\
    &\le \varepsilon C |x|(|v|^{\alpha-1}+1)\\
    &\le \frac{1}{100}a_4|v|^{\alpha}+C\varepsilon^\alpha|x|^\alpha+C.
\end{align*}
In the above, $C=C(\alpha)$ is a positive constant independent of $\varepsilon$. Likewise,
\begin{align*}
    -\varepsilon\la x,v\ra \le \frac{1}{2}|v|^2 +C \varepsilon^2|x|^2.
\end{align*}
Altogether, we get
\begin{align*}
    LV_3 & = L\big(J(x,v)+\varepsilon\la x,v\ra\big)  \\
     &\le -\big(\frac{1}{2}-\varepsilon\big)|v|^2-\frac{1}{2}a_4|v|^\alpha-\varepsilon a_2|x|^\lambda+C\varepsilon^\alpha|x|^\alpha+C\varepsilon^2 |x|^2+C.
\end{align*}
Since $2<\alpha \le \lambda$, observe that by taking $\varepsilon$ sufficiently small, the term $C\varepsilon^\alpha|x|^\alpha+C\varepsilon^2 |x|^2$ can be subsumed into $-\varepsilon a_2|x|^\lambda$. It follows that
\begin{align*}
    LV_3 
     &\le -\big(\frac{1}{2}-\varepsilon\big)|v|^2-\frac{1}{2}|v|^\alpha-\frac{1}{2}a_4\varepsilon a_2|x|^\lambda+C.
\end{align*}
We invoke condition \eqref{cond:U:U(x)=O(x^lambda+1)}, namely, $U(x)\le a_1(|x|^\lambda+1)$, to deduce that 
\begin{align*}
    -\big(\frac{1}{2}-\varepsilon\big)|v|^2-\frac{1}{2}a_4|v|^\alpha-\frac{1}{2}\varepsilon a_2|x|^\lambda \le - c\big(U(x)+\frac{1}{2}|v|^2+\varepsilon\la x,v\ra\big)+C=-c V_3+C.
\end{align*}
As a consequence,
\begin{align*}
     LV_3 
     \le  -cV_3+C.
\end{align*}
This finishes the Lyapunov inequality \eqref{ineq:Lyapunov:poly-damping:exponential} in the setting $\alpha\le \lambda$. \vspace{.2cm}

(ii) With regard to the scenario $\alpha>\lambda$, we first establish \eqref{ineq:Lyapunov:poly-damping:polynomial} for the case $n=1$. By the chain rule (for the generator $L$ defined in~\eqref{form:L}), we have $$L(\phi(f))=\phi'(f)Lf+\frac{\sigma}{2}\phi''(f)|\grad_v f|^2,$$
producing
\begin{align}\label{eq:LVThree-JSquare}
LJ(x,v)^2
&=
2J(x,v)LJ(x,v)+\sigma|v|^2 \nonumber\\
&\le 
-2J(x,v)|v|^2 -a_4 J(x,v)|v|^\alpha +CJ(x,v) +\sigma|v|^2,
\end{align}
where the last estimate above follows from \eqref{ineq:LJ}. Next, we compute
\begin{align} \label{eq:LG}
LG(v)
&= \la\grad G(v),-v-\grad U(x)-\grad G(v)\ra + \frac{\sigma}{2} \triangle_v G(v) \nonumber\\ 
&= -\la v, \grad G(v) \ra -\la\grad G(v),\grad U(x)\ra -|\grad G(v)|^2 + \frac{\sigma}{2}\triangle_v G(v).
\end{align}
Also, since $\grad_x \la v, \grad U(x) \ra = \nabla^2U(x)v$, $\grad_v \la v, \grad U(x)\ra=\grad U(x)$, and $\triangle_v \la v, \grad U(x)\ra = 0$, we have by~\eqref{form:L}
\begin{align}\label{eq:LVThree-Cross-Term}
L\la v,\grad U(x)\ra
&= \la \nabla^2U(x)v,v\ra + \la \grad U(x),-v-\grad U(x) - \grad G(v)\ra \nonumber\\
&= \la v,\nabla^2U(x)v\ra -\la v,\grad U(x)\ra -|\grad U(x)|^2 -\la \grad U(x),\grad G(v)\ra.
\end{align}
Set $F(x,v):=\grad U(x)+\grad G(v)$. Then combining~\eqref{eq:LG} and~\eqref{eq:LVThree-Cross-Term}, we obtain
\begin{align}\label{eq:LVThree-Correction}
L\big(G(v)+\la v,\grad U(x)\ra\big)
= -|F(x,v)|^2 -\la v,F(x,v) \ra + \la v,\nabla^2 U(x) v \ra +
\frac{\sigma}{2} \triangle_v G(v).
\end{align}
Using Young's inequality and $J(x,v)\ge\frac{1}{2}|v|^2$, we have, for every $\delta>0$,
\begin{align} \label{eq:LVThree-Correction-One}
-\la v,F(x,v)\ra
\le \frac14|F(x,v)|^2+|v|^2
\le \frac14|F(x,v)|^2+\delta J(x,v)|v|^2+C_\delta.
\end{align}
Also, by Assumption~\ref{cond:G:poly}, condition \eqref{cond:G:poly:grad^2.G}, it holds that
\begin{align} \label{eq:LVThree-Correction-Two}
\triangle_vG(v) = \operatorname{Tr}\bigl(\nabla_v^2G(v)\bigr) \le c|v|^{\alpha-2}
\le
\delta J(x,v)|v|^\alpha+C_\delta. 
\end{align}
Moreover, by Cauchy-Schwarz and~\eqref{cond:U:U''(x)=O(x^lambda)}, we have
\begin{align} \label{eq:preHessian-Term-Bound}
\big|\la v,\nabla^2U(x)v\ra\big|
\le
a_1\left(1+|x|^{\lambda-2}\right)|v|^2.
\end{align}
Since $\lambda \ge 2$, Young's inequality and the lower bound in~\eqref{cond:U:U(x)=O(x^lambda+1)} imply for every $\delta > 0$,
\[
a_1(1+|x|^{\lambda-2})
\le
\frac{\delta}{a_1} |x |^\lambda + C
\le
\delta  U(x) + C
\le
\delta J(x,v)  +C,
\]
Since $ J(x, v) \ge \frac{1}{2} |v|^2$, Young's inequality gives
\[
C|v|^2
\le
\delta \frac{1}{2}|v|^4+C
\le
\delta J(x,v)|v|^2+C.
\]
Substituting this into~\eqref{eq:preHessian-Term-Bound} yields 
\begin{align} \label{eq:LVThree-Correction-Three}
\big|\la v,\nabla^2U(x)v\ra\big|
\le
2\delta J(x,v)|v|^2 + C_\delta,
\end{align}
for every $\delta>0$. Using the estimates~\eqref{eq:LVThree-Correction-One},~\eqref{eq:LVThree-Correction-Two}, and~\eqref{eq:LVThree-Correction-Three} in~\eqref{eq:LVThree-Correction}, and relabeling $\delta$ and $C_\delta $ if necessary, we obtain
\begin{align}\label{eq:LVThree-Correction-Bound}
L\big(G(v)+\la v,\grad U(x)\ra\big)
\le
-\frac{3}{4} |F(x, v)|^2 + \delta J(x,v) |v|^2 + \delta J(x,v) |v|^\alpha + C_\delta.
\end{align}

Next, we aim to establish an estimate for $J$. Choose $C>0$ sufficiently large that $a_3 \le \frac{1}{2} a_2 C^\lambda$. On the one hand, for every $|x| \ge C$, by condition \eqref{cond:U:x.U'(x)>-x^(lambda+1)} and the Cauchy-Schwarz inequality,
\begin{align} \label{cond:U:U'>x^(lambda-1)}
    |\grad U(x)|
\ge
\frac{1}{|x|}\la \grad U(x),x\ra
\ge
a_2|x|^{\lambda-1}-\frac{a_3}{|x|}
\ge
\frac{a_2}{2}|x|^{\lambda-1}.
\end{align}
Since $\lambda\ge2$, we have $2 \lambda-2 \ge \lambda$. Hence, for $|x| \ge \max\{C,1 \}$, the preceding estimate gives
\[
|\grad U(x)|^2 \ge c|x|^{2\lambda-2} \ge c|x|^\lambda.
\]
Using the upper bound in~\eqref{cond:U:U(x)=O(x^lambda+1)}, we therefore obtain
\[
U(x)
\le
a_1\left(1+|x|^\lambda\right)
\le
a_1\left(1+\frac{1}{c}|\grad U(x)|^2\right), \qquad |x| \ge \max\{C,1\}.
\]
On the other hand, since $U$ is continuous, it is bounded on the compact set $\{x \in \rbb^d : |x| \le \max\{C,1\} \}$. It follows that there exists $C_1>0$ sufficiently large such that
\begin{align} \label{eq:U(x)-upper-bound}
U(x)
\le
C_1\left(1+|\grad U(x)|^2\right), \qquad x\in\rbb^d.
\end{align}
Consequently, for every $\delta>0$, we claim that
\begin{align}\label{eq:J-Coercivity-Ergodicity}
J(x,v)
\le
\delta J(x,v)|v|^2 + C_\delta\left(1+|F(x,v)|^2\right).
\end{align}
Indeed, when $|v|^2\ge\delta^{-1}$, we have $\delta |v|^2 \ge 1 $, and hence $J(x,v) \le \delta J(x,v) |v|^2$, which gives~\eqref{eq:J-Coercivity-Ergodicity} in this case. Suppose now that $|v |^2 < \delta^{-1}$. Then, by Assumption~\ref{cond:G:poly}, we have
\[
|\grad G(v)| \le  a_4|v|^{\alpha-1}+a_5(|v|^{\alphat-1}+c) \le C.
\]
Since $F(x,v) = \grad U(x) + \grad G(v)$, it follows that
\[
|\grad U(x)|^2
=
|F(x,v)-\grad G(v)|^2
\le
2|F(x,v)|^2+2|\grad G(v)|^2
\le
2|F(x,v)|^2+C.
\]
The above together with~\eqref{eq:U(x)-upper-bound} yields $U(x) \le C_\delta \left(1 + |F(x, v) |^2 \right)$. Moreover, $|v |^2 < \delta^{-1}$ implies $ \frac{1 }{2} |v|^2 \le \frac{ 1}{ 2 \delta} $. Therefore, by the definition of $J$, we get
\[
J(x,v)
=
U(x)+\frac{1}{2}|v|^2 \le C_\delta \left(1+ |F(x,v) |^2 \right),
\]
after enlarging $C_\delta$ if necessary. Combining the two cases
yields~\eqref{eq:J-Coercivity-Ergodicity}.

From~\eqref{Def:V-Three},~\eqref{eq:LVThree-JSquare}, and~\eqref{eq:LVThree-Correction-Bound}, we get
\begin{align*}
LV_3(x,v)
&\le
-2J(x,v)|v|^2 -a_4 J (x,v)|v|^\alpha +CJ(x,v) +\sigma|v|^2 \nonumber\\
&\quad
-\frac{3 C_0}{4}| F(x,v)|^2 + C_0 \delta J(x,v)|v|^2 + C_0 \delta J(x,v )|v|^\alpha + C_0 C_\delta.
\end{align*}
By~\eqref{eq:J-Coercivity-Ergodicity}, we have, for every $\delta_1>0$,
\[
CJ(x,v)+\sigma|v|^2 \le CJ(x,v) \le C\delta_1 J(x,v)|v|^2+CC_{\delta_1}\big(1+F(x,v)^2\big),\]
whence
\begin{align*}
LV_3(x,v)
&\le
-\left(2-C_0 \delta-C\delta_1\right)J(x,v)|v|^2 -\left(a_4 - C_0 \delta \right) J(x, v) |v |^\alpha \\
&\quad
-\left( \frac{3 C_0}{4}-C_{\delta_1}C\right)|F(x,v)|^2 + C_0 C_\delta + CC_{\delta_1}. 
\end{align*}
In the above, we emphasize that the constant $C$ is independent of $\delta$ and $\delta_1$. Thus, by first taking $\delta_1$ sufficiently small, $C_0$ sufficiently large and then $\delta$ small enough, we arrive at
\begin{align}\label{eq:LVThree-Coercive}
LV_3(x,v)
&\le
-\frac{3}{2} J(x,v)|v|^2 -\frac{a_4}{2} J(x,v)|v|^\alpha - \frac{C_0}{2} |F(x, v)|^2 + C \nonumber\\
&\le
-c\left( J(x,v)|v|^2 + J(x,v) |v|^\alpha + |F(x, v)|^2 \right) + C,
\end{align}
for suitable constants $c,C$ independent of $(x,v)$.

It remains to compare the right-hand side of~\eqref{eq:LVThree-Coercive} with $V_3$. By Cauchy-Schwarz and Assumption~\ref{cond:U}, conditions \eqref{cond:U:U(x)=O(x^lambda+1)}--\eqref{cond:U:U'(x)=O(x^lambda)}, we infer that
\begin{align} \label{eq:LVThree-Correction-Three-v-GradU}
\big|\la v,\grad U(x)\ra\big|
&\le
a_1|v|\left(1+|x|^{\lambda-1}\right) \nonumber\\
&\le
C\sqrt{J(x,v)} \left[ 1+ J(x,v)^{\frac{\lambda - 1}{\lambda}} \right] \nonumber\\
&=
C\left( J(x,v)^{1/2} + J(x,v)^{ \frac{3}{2} - \frac1 \lambda} \right)
\le \delta J(x,v)^2+C_\delta,
\end{align}
for every $\delta>0$. Next, thanks to Assumption~\ref{cond:G:poly}, condition \eqref{cond:G:poly:grad.G}, we recast $G$ as follows:
\begin{align*}
    G(v)-G(0) = \int_0^1\la \grad G(sv),v\ra \d s & = \frac{a_4}{\alpha}|v|^\alpha +\int_0^1\la \grad G(sv) -sv|sv|^{\alpha-2},v\ra \d s\\
    &\ge \frac{a_4}{\alpha}|v|^\alpha - a_5\int_0^1(|sv|^{\alphat-1}+1)|v|\d s\\
    &\ge \frac{a_4}{2\alpha}|v|^\alpha-C.
\end{align*}
In other words, we have
\begin{align*}
    G(v)\ge \frac{a_4}{2\alpha}|v|^\alpha-C.
\end{align*}
Likewise, 
\begin{align*}
    G(v)\le 2\frac{a_4}{\alpha}|v|^\alpha+C.
\end{align*}
Hence, on the one hand, by~\eqref{Def:V-Three},~\eqref{eq:LVThree-Correction-Three-v-GradU}, and taking $\delta$ small enough, it holds that
\begin{align}\label{eq:VThree-Lower-Bound}
V_3(x,v) \ge \left(1-C_0\delta\right)J(x,v)^2 +C_0\frac{a_4}{2\alpha}|v|^\alpha +c_0-C.
\end{align}
It is clear that by choosing $c_0>0$ sufficiently large, we obtain $V_3(x,v)\ge1$ for all $(x,v) \in \rbb^{2d}$. On the other hand, we have
\begin{align}
V_3(x,v)
&\le
\left(1+C_0\delta\right)J(x,v)^2 + C_0\frac{2a_4}{\alpha}|v|^\alpha + C \nonumber\\
&\le
C\left( 1+J(x,v)^2+|v|^\alpha \right), \nonumber
\end{align}
whence,
\begin{align*}
    V_3(x,v)^{1/2}
&\le
C\left(1+J(x,v)+|v|^{\alpha/2}\right).
\end{align*}
We invoke \eqref{eq:J-Coercivity-Ergodicity} together with the estimate $|v |^{ \alpha/2} \le C\left( 1 + J(x,v) |v|^\alpha \right)$ to deduce further that
\begin{align}\label{eq:VThree-Powers-Comparison}
V_3(x,v)^{1/2}
&\le
C\left(
1+J(x,v)|v|^2+J(x,v)|v|^\alpha+|F(x,v)|^2
\right).
\end{align}
Combining \eqref{eq:VThree-Powers-Comparison} with~\eqref{eq:LVThree-Coercive} yields
\[
LV_3(x,v) \le -cV_3(x,v)^{1/2}+C.
\]
This establishes \eqref{ineq:Lyapunov:poly-damping:polynomial} for the base case $n=1$.

Turning to the case $n > 1$, for notational convenience, we set
\[
D(x,v):=J(x,v)|v|^2+J(x,v)|v|^\alpha+|F(x,v)|^2.
\]
Since $\nabla_vV_3=2Jv+C_0F$ by~\eqref{Def:V-Three}, it follows that
\[
|\nabla_vV_3|^2 \le 8J^2|v|^2+2C_0^2|F|^2 \le C V_3^{1/2}D.
\]
The last inequality holds since $J|v|^2\le D$ and $|F|^2\le D$ by the definition of $D$, and by~\eqref{eq:VThree-Lower-Bound}, with $c_0$ sufficiently large, we have $V_3\ge1+\frac12J^2+\frac{C_0}{\alpha}|v|^\alpha$, which in particular implies $J\le C V_3^{1/2}$ and $V_3^{1/2}\ge1$. Using the chain rule, we obtain the equality below.
\begin{align}\label{eq:LVThree-Powers-Drift}
L\big(V_3^n\big)
&=nV_3^{n-1}LV_3+\frac{\sigma}{2}n(n-1)V_3^{n-2}|\nabla_vV_3|^2 \nonumber \\
&\le -ncV_3^{n-1}D+nCV_3^{n-1}+C_nV_3^{n-\frac32}D \nonumber\\
&\le -\frac{nc}{2}V_3^{n-1}D+C_nV_3^{n-1}.
\end{align}
The first inequality follows from the definition of $D$,~\eqref{eq:LVThree-Coercive}, and the preceding estimate for $|\nabla_vV_3|^2$. For the final inequality in~\eqref{eq:LVThree-Powers-Drift}, observe that
\[
C_nV_3^{n-\frac32}D
=
C_nV_3^{-1/2}V_3^{n-1}D.
\]
Since $V_3^{-1/2}\to0$ as $V_3\to\infty$, we may choose $R_n\ge1$ sufficiently large so that $C_nV_3^{-1/2}\le\frac{1}{2}nc$ whenever $V_3\ge R_n$. Hence, on $\{V_3\ge R_n\}$, we have $C_nV_3^{n-\frac32}D\le\frac{nc}{2}V_3^{n-1}D$, and therefore the final inequality in~\eqref{eq:LVThree-Powers-Drift} holds in this region after enlarging $C_n$ if necessary. On the other hand, on the compact set $\{V_3\le R_n\}$, since $L(V_3^n)$, $V_3^{n-1}$, and $D$ are continuous, enlarging $C_n$ if necessary, we have
\[
L(V_3^n)+\frac{nc}{2}V_3^{n-1}D
\le C_n
\le C_nV_3^{n-1},
\]
where the second inequality uses $V_3^{n-1}\ge1$, and hence the final inequality in~\eqref{eq:LVThree-Powers-Drift} holds in this region as well.

Since~\eqref{eq:VThree-Powers-Comparison} gives
$D\ge \frac{1}{C}V_3^{1/2}-1$, combining this with~\eqref{eq:LVThree-Powers-Drift} produces
\[
L\big(V_3(x,v)^n\big) \le -cV_3(x,v)^{n-\frac12} + C V_3(x,v)^{n-1} \le -c V_3(x,v)^{n-\frac12} + C.
\]
This establishes the Lyapunov estimate \eqref{ineq:Lyapunov:poly-damping:polynomial} for the general case $n>1$. The proof is thus complete.

\end{proof}

\subsection{Singular damping} \label{sec:ergodicity:singular-damping}

We now turn to the scenario of singular damping under Assumption \ref{cond:G:singular} when the leading singular parameter $\beta=1$ and under the extra condition of Assumption \ref{cond:G:singular:beta>1} when $\beta>1$. Similar to the presentation of Section \ref{sec:ergodicity:poly-damping}, we first discuss the heuristic construction before giving the explicit argument.

Since the difficulty in this case is to find a function $V$ such that $LV$ can produce a dissipation effect on the singularity, there are two cases to be considered depending on the leading singular parameter $\beta$. 

Case 1: $\beta=1$. In this case, we note that $G$ satisfying Assumption \ref{cond:G:singular} behaves like a logarithm function. Following the same approach as in the polynomial damping setting, the candidate Lyapunov function should contain the term 
\begin{align*}
 G(v)+\la v,U(x)\ra\sim -\log |v|+ \la v,U(x)\ra.
\end{align*}
Moreover, to control $\la v,\grad U(x)\ra$, it is sufficient to incorporate the function $U(x)+\frac{1}{2}|v|^2$. However, if we choose $(U(x)+\frac{1}{2}|v|^2)^2$ as in the setting of polynomial damping, we will need to control the bad term $U(x)$ as in  
\begin{align*}
    L\Big[\Big(U(x)+\frac{1}{2}|v|^2\Big)^2\Big] = \dots+ U(x)+\cdots.
\end{align*}
This turns out to be challenging owing to the presence of the singular $\grad G(v)$. To circumvent the issue, we pick $(U(x)+\frac{1}{2}|v|^2)^p$ for suitable $p\in(1,2)$, which together with $-\log |v|+ \la v,U(x)\ra$ suggests a candidate of the form
\begin{align*}
    V=\Big(U(x)+\frac{1}{2}|v|^2\Big)^{p}-\log |v|+\la v,\grad U(x)\ra.
\end{align*}
On the one hand, recalling $L$ from \eqref{form:L}, it is important to observe that $\triangle_v (-\log |v|) = -\frac{d-2}{|v|^{2}}$. So that $L(-\log|v|)$ already produces a dissipation effect for the singularity in dimension $d\ge 3$. On the other hand, in dimension $d=2$, the dissipation effect from $-\log|v|$ is canceled, requiring a modification by adding $-|v|^\theta$ for suitable $\theta\in (0,1)$ to the above $V$. We note that such a modification also appeared in Lemma \ref{lem:well-posed:singular-damping} for the well-posedness. The exact form of $V$ in the case $\beta=1$ can be found in expression \eqref{Def:V-Four}.

Case 2: $\beta>1$. In this case, since $G$ behaves like $1/|v|^{\beta-1}$, a first ansatz is
\begin{align*}
    V = \Big(U(x)+\frac{1}{2}|v|^2\Big)^{p} +\frac{1}{|v|^{\beta-1}}+\la v,\grad U(x)\ra.
\end{align*}
However, a routine calculation shows that $\triangle_v (|v|^{-\beta+1})$ has the form $-(d-\beta+1)|v|^{-(\beta+1})$, which is a bad term if $\beta>d-1$. To circumnavigate the issue, we draw upon the technique in the setting of singular potentials appearing in the work of \cite{duong2024asymptotic,lu2019geometric}. Recall that for the singular potential setting \eqref{eq:Langevin:singular-potential}, a correction to the usual Lyapunov function $U(x)+\frac{1}{2}|v|-\varepsilon\la x,v\ra$ for \eqref{eq:Langevin:singular-potential} is given by
\begin{align*}
    U(x)+\frac{1}{2}|v|^2-\varepsilon\frac{\la x,v\ra}{|x|}+c_0.
\end{align*}
Motivated by the above expression, our candidate for \eqref{eq:Langevin-SDE} can be modified as
\begin{align*}
    V = \Big(U(x)+\frac{1}{2}|v|^2\Big)^{p} +\frac{1}{|v|^q}\Big(\frac{1}{|v|^{\beta-1}}+\la v,\grad U(x)\ra\Big),
\end{align*}
for a suitable constant $q\in (0,1]$. As it turns out, all important parameters $\beta,\lambda,p,q$ and the dimension $d$ are closely related, hence the extra assumption \ref{cond:G:singular:beta>1}. The exact form of $V$ when $\beta>1$ can also be found below in \eqref{Def:V-Four}.

\begin{lemma} \label{lem:ergodicity:singular-damping}
In dimension $d\ge2$, let $U$ satisfy Assumptions~\ref{cond:U} and $G$ satisfy Assumption~\ref{cond:G:singular}. Define
\begin{align}\label{Def:V-Four}
V_4(x,v) :=
\begin{cases}
J(x,v)^p-a_4\log|v|+\la v,\grad U(x)\ra+c_0, & d\ge3,\ \beta=1,\\[1mm]
J(x,v)^p-|v|^\theta-a_4\log|v|+\la v,\grad U(x)\ra+c_0, & d=2,\ \beta=1,\\[1mm]
J(x,v)^p+\tfrac{a_4}{\beta+q-1}|v|^{-(\beta+q-1)}+\tfrac{a_6}{\betat+q-1}|v|^{-(\betat+q-1)}
\\[1mm]
\hspace{1cm}+\frac{1}{1-q}|v|^{-q}\la v,\grad U(x)\ra+c_0, & \beta>1,
\end{cases}
\end{align}
for $(x,v)\in\X$, where $J(x,v):=U(x)+\frac12|v|^2$, $\beta,\betat,a_4,a_6$ are constants from Assumption \ref{cond:G:singular}, and $c_0$ is a positive constant. Suppose that the parameters $p,q,\theta$ satisfy the following conditions depending on dimension $d$, $\beta$ and $\lambda$ where $\lambda$ is the constant from Assumption \ref{cond:U}\textup{:}
\begin{itemize}
\item[\textup{(i)}] $d\ge3$, $\beta=1$\textup{:}
\[
2-\frac{2}{\lambda}<p<2.
\]

\item[\textup{(ii)}] $d=2$, $\beta=1$\textup{:}
\[
0<\theta<\min\Big\{\frac{\lambda}{\lambda-1},2(1-\betat)\Big\},
\qquad
2-\frac{2}{\lambda}<p<\min\Big\{2,\,
1+\frac{(2-\theta)(\lambda-1)}{\lambda}\Big\}.
\]

\item[\textup{(iii)}] $d\ge2$ and $\beta>1$\textup{:} Suppose further that $G$ is given by Assumption \ref{cond:G:singular:beta>1} and that
\[
\begin{cases}
\beta>1, & 2\le\lambda\le4,\\[1mm]
1<\beta<\dfrac{6(\lambda-1)}{\lambda-4}, & \lambda>4.
\end{cases}
\]
The parameters $q$ and $p$ satisfy
\begin{align}\label{eq:LVFour-q-range}
\max\left\{0,\,1-\frac{d-1}{\beta},\,
\frac{2\beta(\lambda-2)-4(\lambda-1)}{\beta\lambda+2(\lambda-1)}\right\}<q<1,
\qquad q\neq1-\betat,
\end{align}
and
\begin{align}\label{eq:LVFour-p-range}
\max\left\{1,\,2-\frac{2}{\lambda}-\frac q2\right\}
<p<\min\left\{3-\frac{2}{\lambda},\,
1+\frac{(q+2)(\lambda-1)}{\beta\lambda}\right\}.
\end{align}
\end{itemize}
Then for all $c_0$ sufficiently large, $V_4\ge1$ and, for every real $n\ge1$, there exist positive constants $c_n,C_n$ such that
\begin{align} \label{ineq:Lyapunov:singular-damping:polynomial}
    L\bigl(V_4(x,v)^n\bigr)\le-c_nV_4(x,v)^{n-\frac1p}+C_n,
\qquad (x,v)\in\X.
\end{align}

\end{lemma}

\begin{remark} \label{rem:G:singular:beta_n}
    Following closely the proof of Lemma \ref{lem:ergodicity:singular-damping} for the case $d\ge 2$ and $\beta>1$, we see that the argument therein can be generalized to cover function $G$ satisfying
    \begin{align*}
        \grad G(v) =-a\frac{ v}{|v|^{\beta+1}}-\sum_{i=1}^n c_i\frac{v}{|v|^{\beta_i+1}}, 
    \end{align*}
    where $0< \beta_1<\beta_2<\dots <\beta_n<\beta$, $c_i\in\rbb$. A suitable choice for the Lyapunov function is given by
    \begin{align*}
        V(x,v) &= J(x,v)^p+ \frac{a}{\beta+q-1}|v|^{-(\beta+q-1)}+\sum_{i=1}^n \frac{c_i}{\beta_i+q-1} | v|^{-(\beta_i+q-1)}\\
       &\qquad +\frac{1}{1-q}|v|^{-q}\la v,\grad U(x)\ra +c_0.
    \end{align*}
\end{remark}

\begin{proof}[Proof of Lemma~\ref{lem:ergodicity:singular-damping}]

1. In the first part of the proof, we proceed to establish the Lyapunov inequality \eqref{ineq:Lyapunov:singular-damping:polynomial} when $n=1$. For notational convenience, we set
\begin{align} \label{form:R(v)=gradG+v/v^beta}
R (v) =  \grad G(v)+a_4\frac{v}{|v|^{\beta+1}}.
\end{align}
In view of Assumption \ref{cond:G:singular}, condition  \eqref{cond:G:singular:|gradG|<1/|v|^beta}, $R$ satisfies
\begin{align} \label{ineq:R(v)}
    |R(v)| \le |a_6|\Big(\frac{1}{|v|^{\betat}}+1\Big).
\end{align}

(i) Suppose that $d\ge3$ and $\beta=1$.
Recall from \eqref{eq:LJ:singular} that,
\begin{align*}
    LJ(x,v)=-|v|^2+\frac{1}{2}(2a_4+\sigma d)-\la v,R(v) \ra.
\end{align*}
Consequently, using the chain rule for the generator $L$, we obtain
\begin{align}\label{eq:LVFour-Jp}
LJ(x,v)^p
&=
pJ(x,v)^{p-1}LJ(x,v) +\frac{\sigma}{2}p(p-1)J(x,v)^{p-2}|v|^2 \nonumber\\
&= -pJ(x,v)^{p-1}|v|^2
+p\frac{2a_4+\sigma d}{2}J(x,v)^{p-1} -\la v,R(v)\ra J(x,v)^{p-1} \notag \\
&\qquad+\frac{\sigma p(p-1)}{2} J(x,v)^{p-2}|v|^2.
\end{align}
Next, from \eqref{eq:L(-log|v|)}, we have
\begin{align*}
    &L(-a_4\log|v|)\\
    & =\Big\la -a_4\frac{v}{|v|^2}, -v-\grad U(x)-\grad G(v) \Big\ra
-
\frac{a_4\sigma(d-2)}{2|v|^2}\\
& =\Big\la \grad G(v)-R(v), -v-\grad U(x)-\grad G(v) \Big\ra
-
\frac{a_4\sigma(d-2)}{2|v|^2}\\
&= - \la \grad G(v),\grad U(x)+\grad G(v)\ra-\la \grad G(v),v\ra + \la R(v),v+\grad U(x)+\grad G(v)\ra -
\frac{a_4\sigma(d-2)}{2|v|^2}.
\end{align*}
Combining the above with~\eqref{eq:LVThree-Cross-Term}, we obtain
\begin{align}
&L\big(-a_4\log|v|+\la v,\grad U(x)\ra\big) \\
&= \la v,\nabla^2U(x)v\ra -\la v,\grad U(x)+\grad G(v)\ra+\la R(v),v+\grad U(x)+\grad G(v)\ra \notag \\
&\qquad-|\grad U(x)+\grad G(v)|^2
 -\frac{a_4\sigma}{2}\frac{d-2}{|v|^2}. \nonumber
\end{align}
Therefore, by~\eqref{Def:V-Four},~\eqref{eq:LVFour-Jp}, and the above identity, we have
\begin{align}\label{eq:LVFour-d-geq-3}
&LV_4(x,v) \notag \\
&= -pJ(x,v)^{p-1}|v|^2
+p\frac{2a_4+\sigma d}{2}J(x,v)^{p-1} -\la v,R(v)\ra J(x,v)^{p-1}+\frac{\sigma p(p-1)}{2} J(x,v)^{p-2}|v|^2 \nonumber\\
&\quad +\la v,\nabla^2U(x)v\ra 
+\la R(v)-v,\grad U(x)+\grad G(v)\ra+\la R(v),v\ra \notag \\
&\qquad-|\grad U(x)+\grad G(v)|^2  -a_4\frac{\sigma}{2} \frac{d-2}{|v|^2}.
\end{align}
By Young's inequality,
\[
\big| \la R(v)- v,\grad U(x)+\grad G(v)\ra \big| \le
\frac{1}{4}|\grad U(x)+\grad G(v)|^2 +|R(v)-v|^2.
\]
Since $R(v)$ satisfies \eqref{ineq:R(v)}, we get
\begin{align*}
    |R(v)-v|^2+|\la v,R(v)\ra| \le C(|R(v)|^2+|v|^2) \le C\Big(\frac{1}{|v|^{2\betat}}+|v|^2 +1\Big).
\end{align*}
Consequently,
\begin{align*} 
&-|\grad U(x)+\grad G(v)|^2 +\la R(v)-v,\grad U(x)+\grad G(v)\ra +\la R(v),v\ra\notag \\
&\qquad\le
-\frac{3}{4}|\grad U(x)+\grad G(v)|^2 +C\Big(\frac{1}{|v|^{2\betat}}+|v|^2 +1\Big) .
\end{align*}
Using $|x|^\lambda\le a_1J(x,v)+a_1^2$, we obtain $|x|^{\lambda-2} \le C\big(1+J(x,v)^{\frac{\lambda-2}{\lambda}}\big)$. Also, since $\frac{\lambda-2}{\lambda}<p-1$, Young's inequality gives, for every $\delta>0$, $J(x,v)^{\frac{\lambda-2}{\lambda}} \le \delta J(x,v)^{p-1}+C_\delta$. Therefore, after enlarging $C_\delta$ and choosing the coefficient of
$J(x,v)^{p-1}$ sufficiently small, we may write
\[
C_0a_1\bigl(1+|x|^{\lambda-2}\bigr)
\le \delta J(x,v)^{p-1}+C_\delta.
\]
By Cauchy-Schwarz inequality and condition~\eqref{cond:U:U''(x)=O(x^lambda)}, and using the preceding estimate, we obtain
\[
\big|\la v,\nabla^2U(x)v\ra\big|
\le a_1\bigl(1+|x|^{\lambda-2}\bigr)|v|^2
\le \frac{1}{100}\delta J(x,v)^{p-1}|v|^2+C_\delta|v|^2.
\]
Also, (recalling \eqref{ineq:R(v)} and $\betat\in[0,1)$ when $\beta=1$)
\begin{align*}
    -\la v,R(v)\ra J(x,v)^{p-1} \le |a_6|(|v|^{1-\betat}+|v|)J(x,v)^{p-1} &\le C(1+|v|)J(x,v) \\
    &\le C_\delta J(x,v)^{p-1}+\frac{1}{100}\delta J(x,v)^{p-1} |v|^2.
\end{align*}
Moreover, since $J(x,v)\ge\frac12|v|^2$ and $p>1$, Young's inequality gives, for every $\delta>0$,
\begin{align} \label{ineq:|v|^2<J^(p-1)|v|^2}
    \bigl(C+C_\delta\bigr)|v|^2
\le \frac{1}{100}\delta J(x,v)^{p-1}|v|^2 + C.
\end{align}
Therefore, combining the above observations, we obtain
\begin{align}\label{eq:LVFour-Energy-Estimate}
&-pJ(x,v)^{p-1}|v|^2- \la v,R(v)\ra J(x,v)^{p-1} + \la v,\nabla^2U(x)v\ra \notag \\
&\qquad-|\grad U(x)+\grad G(v)|^2 +\la R(v)-v,\grad U(x)+\grad G(v)\ra+\la v,R(v )\ra \notag \\
&\le -\bigl(p-\delta \bigr) 
J(x,v)^{p-1}|v|^2 -\frac{3}{4}|\grad U(x)+\grad G(v)|^2+ C_{\delta}\Big(\frac{1}{|v|^{2\betat}} +J(x,v)^{p-1}+1\Big)\nonumber\\
&\le -\frac{p}{2}J(x,v)^{p-1}|v|^2-\frac{3}{4}|\grad U(x)+\grad G(v)|^2
+  C_{\delta}\Big(\frac{1}{|v|^{2\betat}} +J(x,v)^{p-1}+1\Big),
\end{align}
where the second inequality holds by choosing $\delta$ sufficiently small.

Next, we claim that, for every $\eta>0$,
\begin{align}\label{eq:LVFour-Remaining-Estimate}
J(x,v)^{p-1} \le
\eta J(x,v)^{p-1}|v|^2 + \eta|\grad U(x)+\grad G(v)|^2 +
\frac{\eta}{|v|^2} + C_\eta.
\end{align}
Indeed, fix $\eta>0$ and choose $R_\eta:=\eta^{-1/2}$. For $|v|\ge R_\eta$, we have
\[
J(x,v)^{p-1}\le\eta J(x,v)^{p-1}|v|^2.
\]
It remains to consider $|v|<R_\eta$. From the estimate \eqref{cond:U:U'>x^(lambda-1)}, i.e., $|x|^{\lambda-1}\le C(1+|\grad U(x)|)$, and condition \eqref{cond:U:U(x)=O(x^lambda+1)}, i.e., $U(x)\le C\bigl(1+|x|^\lambda\bigr)$, for $|v|<R_\eta$, it holds that
\[
J(x,v)^{p-1}
=
\Big(U(x)+\frac12|v|^2 \Big)^{p-1}
\le C\left(1+|\grad U(x)|^{\frac{\lambda(p-1)}{\lambda-1}}\right)
\le\eta|\grad U(x)|^2+C_{\eta},
\]
for every $\eta>0$. The final inequality follows from Young's inequality, since
$0<\lambda(p-1)/(\lambda-1)<2$ by $1<p<2$ and $\lambda\ge2$. Furthermore,
since 
\begin{align*}
  |\grad G(v)|\le a_4|v|^{-1}+|a_6|(|v|^{-\betat}+1)\le C(|v|^{-1}+1)
\end{align*}
by virtue of Assumption \ref{cond:G:singular}, we have
\[
|\grad U(x)|^2\le 2|\grad U(x)+\grad G(v)|^2+ 2|\grad G(v)|^2\le2|\grad U(x)+\grad G(v)|^2+C \Big(\frac{1}{|v|^{2}}+1\Big).
\]
Taking $\eta$ small enough yields
\[
J(x,v)^{p-1}\le\eta|\grad U(x)+\grad G(v)|^2
+\frac{\eta}{|v|^2}+C_\eta,
\qquad |v|<R_\eta.
\]
Combining the two regions proves~\eqref{eq:LVFour-Remaining-Estimate}.

Turning back to \eqref{eq:LVFour-d-geq-3}, since $|v|^2\le2J(x,v)$, we have
\begin{align}\label{eq:LVFour-Remaining-Terms}
p\frac{2a_4+\sigma d}{2}J(x,v)^{p-1}
&+ \frac{\sigma p(p-1)}{2}J(x,v)^{p-2}|v|^2 +C_\delta J(x,v)^{p-1} \notag \\
&\le C_\delta J(x,v)^{p-1} \nonumber \\
&\le \eta J(x,v)^{p-1}|v|^2
+ \eta|\grad U(x)+\grad G(v)|^2
+ \frac{\eta}{|v|^2}
+ C.
\end{align}
The last inequality holds after applying~\eqref{eq:LVFour-Remaining-Estimate} to $J(x,v)^{p-1}$ for a positive constant $\eta/C_\delta$. Combining~\eqref{eq:LVFour-Energy-Estimate} and
\eqref{eq:LVFour-Remaining-Terms} with~\eqref{eq:LVFour-d-geq-3}, we obtain
\begin{align}
LV_4(x,v)
&\le -\left(\frac{p}{2}-\eta\right)
J(x,v)^{p-1}|v|^2 -\left(\frac{3}{4}-\eta\right)
|\grad U(x)+\grad G(v)|^2 \nonumber\\ 
&\quad -a_4\frac{\sigma}{2} \frac{d-2}{|v|^2}+\frac{\eta}{|v|^2}+\frac{C}{|v|^{2\betat}} + C, \nonumber
\end{align}
for every $\eta>0$. Since $\betat\in[0,1)$, we note that $|v|^{-2\betat}$ can be subsumed into $\eta|v|^{-2}$. This together with choosing
$\eta>0$ sufficiently small and the fact that $d\ge 3$ allows us to obtain constants $c,C>0$ such that
\begin{align}\label{eq:LVFour-Coercive}
LV_4(x,v)
\le -cJ(x,v)^{p-1}|v|^2
-c|\grad U(x)+\grad G(v)|^2
-\frac{c}{|v|^2} +C.
\end{align}

It remains to establish that $V_4\ge 1$ and compare the right-hand side of~\eqref{eq:LVFour-Coercive} with $V_4^{(p-1)/p}$. By Cauchy-Schwarz inequality and Assumption~\ref{cond:U}, we obtain the first inequality below.
\[
\big|\la v,\grad U(x)\ra\big|
\le C|v|\bigl(1+|x|^{\lambda-1}\bigr)
\le C\left( 1+
J(x,v)^{\frac12+\frac{\lambda-1}{\lambda}}
\right)
\le \delta J(x,v)^p+C_\delta,
\]
for every $\delta>0$. The second inequality uses that
$|v|\le \sqrt{2}J(x,v)^{1/2}$ and $|x|^\lambda\le a_1J(x,v)+a_1^2$. Since $\frac12+\frac{\lambda-1}{\lambda}<p$, Young's inequality gives the final inequality. Also, when $|v|\ge 1$, 
\begin{align*}
    -a_4\log|v| \ge -a_4|v| \ge -\delta J(x,v)^{p}-C_\delta.
\end{align*}
Therefore, on the one hand, by~\eqref{Def:V-Four} and the preceding estimates, for every $\delta>0$,
\begin{align}
V_4(x,v) &=
J(x,v)^p -a_4\log|v|
+\la v,\grad U(x)\ra
+c_0 \nonumber\\
&\ge
\bigl(1-2\delta\bigr)J(x,v)^p -C_\delta+c_0. \nonumber
\end{align}
By choosing $\delta>0$ sufficiently small and then choosing $c_0$ large enough, we obtain
\[
V_4(x,v) \ge \frac12 J(x,v)^p-C_\delta +c_0 \ge 1.
\]
On the other hand, we also obtain the bound 
\begin{align*}
    V_4(x,v) \le C(1+J(x,v)^p+|\log|v|| ),
\end{align*}
whence,
\begin{align*}
    V_4(x,v)^{\frac{p-1}{p}} &\le
C\left( 1+J(x,v)^{p-1}
+|\log|v||^{\frac{p-1}{p}} \right) \nonumber \\
&\le
C\left( 1+J(x,v)^{p-1}
+|\log|v|| \right).
\end{align*}
In light of~\eqref{ineq:|v|^2<J^(p-1)|v|^2} and \eqref{eq:LVFour-Remaining-Estimate} and the fact that 
\begin{align*}
    \log|v| \le |v|+\frac{1}{|v|},
\end{align*}
we get
\begin{align}\label{eq:LVFour-Control-V}
V_4(x,v)^{\frac{p-1}{p}} 
&\le C\Big(
1 + J(x,v)^{p-1}|v|^2 + |\grad U(x)+\grad G(v)|^2
+ \frac{1}{|v|^2} \Big).
\end{align}
Combining~\eqref{eq:LVFour-Control-V} with~\eqref{eq:LVFour-Coercive}, we conclude that there exist constants $c,C>0$ such that
\[
LV_4(x,v) \le -cV_4(x,v)^{\frac{p-1}{p}} + C, \qquad (x,v)\in\X.
\]
This proves \eqref{ineq:Lyapunov:singular-damping:polynomial} in the case $d\ge3$, $\beta=1$ with $n=1$. \vspace{.2cm}

(ii) Suppose that $d=2$ and $\beta=1$.
By~\eqref{form:L} with $d=2$, we obtain
\[
L|v|^\theta
= -\theta|v|^\theta -\frac{\theta}{|v|^{2-\theta}} \la v,\grad U(x)+\grad G(v)\ra +\frac{\sigma\theta^2}{2} \frac{1}{|v|^{2-\theta}}.
\]
Applying the generator $L$ to~\eqref{Def:V-Four} in this case and using \eqref{eq:LVFour-d-geq-3} with $d=2$, we obtain
\begin{align}\label{eq:LVFour-dTwo-Simplified}
&LV_4(x,v) \notag \\
& = L\Big(J(x,v)^p-|v|^\theta -a_4\log|v|+\la v, \grad U(x)\ra+c_0\Big)\notag  \\
&=-pJ(x,v)^{p-1}|v|^2
+p(a_4+\sigma )J(x,v)^{p-1} -\la v,R(v)\ra J(x,v)^{p-1}+\frac{\sigma p(p-1)}{2} J(x,v)^{p-2}|v|^2 \nonumber\\
&\qquad +\la v,\nabla^2U(x)v\ra 
+\la R(v)-v,\grad U(x)+\grad G(v)\ra+\la R(v),v\ra \notag \\
&\qquad-|\grad U(x)+\grad G(v)|^2 +\theta|v|^\theta+ \frac{\theta}{|v|^{2-\theta}}\la v,\grad U(x)+\grad G(v)\ra   -\frac{\sigma\theta^2}{2}
\frac{1}{|v|^{2-\theta}}.
\end{align}
Recalling from \eqref{eq:LVFour-Energy-Estimate}, we readily have for all $\delta>0$ sufficiently small
\begin{align*}
&-pJ(x,v)^{p-1}|v|^2- \la v,R(v)\ra J(x,v)^{p-1} + \la v,\nabla^2U(x)v\ra \notag \\
&\qquad-|\grad U(x)+\grad G(v)|^2 +\la R(v)-v,\grad U(x)+\grad G(v)\ra+\la v,R(v )\ra  \\
&\le -\frac{p}{2}J(x,v)^{p-1}|v|^2-\frac{3}{4}|\grad U(x)+\grad G(v)|^2
+  C_{\delta}\Big(\frac{1}{|v|^{2\betat}} +J(x,v)^{p-1}+1\Big).
\end{align*}
Since $\theta <2(1-\betat)$, i.e., $2\betat<2-\theta$, we infer that
\begin{align*}
    \frac{C_{\delta}}{|v|^{2\betat}} \le \frac{1}{100}\frac{\delta}{|v|^{2-\theta}}+C.
\end{align*}
Also, since $\theta <2(1-\betat)<2<2p$, we have
\begin{align*}
    \theta |v|^\theta \le \frac{1}{100}\delta |v|^{2p}+C\le \frac{1}{100}\delta J(x,v)^{p-1} |v|^{2}+C.
\end{align*}
Next, by the Cauchy-Schwarz inequality, it holds that
\begin{align*}
\frac{\theta}{|v|^{2-\theta}} \la v,\grad U(x)+\grad G(v)\ra 
&\le \frac{1}{4}|\grad U(x)+\grad G(v)|^2 + 
\frac{\theta^2}{|v|^{2-2\theta}}.
\end{align*}
We note that since $0<\theta<\frac{\lambda}{\lambda-1}\le 2$, we have for all $\delta$
\begin{align*}
   \theta^2|v|^{2\theta-2} &
\le \frac{1}{100}\frac{\delta}{|v|^{2-\theta}}+\frac{1}{100}\delta|v|^2+C_\delta \\
&\le \frac{1}{100}\frac{\delta}{|v|^{2-\theta}}+\frac{1}{100}\delta J(x,v)^{p-1}|v|^2+C_\delta,
\end{align*}
where the last estimate follows from \eqref{ineq:|v|^2<J^(p-1)|v|^2}. It follows that
\begin{align*}
&
\frac{\theta}{|v|^{2-\theta}} \la v,\grad U(x)+\grad G(v)\ra \nonumber\\
&\le \frac{1}{4}|\grad U(x)+\grad G(v)|^2 + \frac{1}{100}\frac{\delta}{|v|^{2-\theta}}+\frac{1}{100}\delta J(x,v)^{p-1}|v|^2+C_\delta .
\end{align*}

Next, similar to \eqref{eq:LVFour-Remaining-Estimate}, we claim that provided that 
\begin{align*}
    p<1+\frac{(2-\theta)(\lambda-1)}{\lambda},
\end{align*}
for every $\eta>0$,
\begin{align}\label{eq:LVFour-dTwo-Remaining-Estimate}
J(x,v)^{p-1}
\le \eta J(x,v)^{p-1}|v|^2 + \eta|\grad U(x)+\grad G(v)|^2
+ \frac{\eta}{|v|^{2-\theta}} + C_\eta.
\end{align}
Indeed, following the proof of \eqref{eq:LVFour-Remaining-Estimate}, on the one hand, for $|v|\ge R_\eta = \eta^{-1/2}$, we readily have
\[
J(x,v)^{p-1}\le\eta J(x,v)^{p-1}|v|^2.
\]
On the other hand, when $|v|<R_\eta$, we once again invoke estimate \eqref{cond:U:U'>x^(lambda-1)} and condition \eqref{cond:U:U(x)=O(x^lambda+1)} to infer that
\[
J(x,v)^{p-1}
=
\Big(U(x)+\frac12|v|^2 \Big)^{p-1}
\le C\left(1+|\grad U(x)|^{\frac{\lambda(p-1)}{\lambda-1}}\right)
\le\eta|\grad U(x)|^{2-\theta}+C_{\eta}.
\]
In the above, the last inequality follows from Young's inequality thanks to the choice of $p,\theta$ satisfying
$p<1+(2-\theta)(\lambda-1)/\lambda$, which is equivalent to $\lambda(p-1)/(\lambda-1)<2-\theta$. Furthermore, we have
\begin{align*}
    |\grad U(x)|^{2-\theta} &\le C|\grad U(x)+\grad G(v)|^{2-\theta}+ C|\grad G(v)|^{2-\theta}\\
    &\le C|\grad U(x)+\grad G(v)|^{2-\theta}+C \Big(\frac{1}{|v|^{2-\theta}}+1\Big).
\end{align*}
Since the constant $C$ does not depend on $\eta$, we get
\begin{align*}
    J(x,v)^{p-1} \le \frac{\eta}{C}|\grad U(x)|^{2-\theta}+C_{\eta} &\le \eta |\grad U(x)+\grad G(v)|^{2-\theta}+\frac{\eta}{|v|^{2-\theta}}+C_\eta\\
    &\le \eta |\grad U(x)+\grad G(v)|^{2}+\frac{\eta}{|v|^{2-\theta}}+C.
\end{align*}
where the last inequality above follows from the fact that $\theta <2$. Combining two cases proves \eqref{eq:LVFour-dTwo-Remaining-Estimate}.

As consequence, using $|v|^2\le2J(x,v)$ and \eqref{eq:LVFour-dTwo-Remaining-Estimate}, we obtain
\begin{align*}
&p(a_4+\sigma)J(x,v)^{p-1}
+ \frac{\sigma p(p-1)}{2}J(x,v)^{p-2}|v|^2+C_\delta J(x,v)^{p-1} \\ 
&\le C J(x,v)^{p-1} \\
&\le \delta J(x,v)^{p-1}|v|^2
+ \delta|\grad U(x)+\grad G(v)|^2+\frac{\delta}{|v|^{2-\theta}} + C_\delta .
\end{align*}
Now, we collect all of the above estimates with the identity~\eqref{eq:LVFour-dTwo-Simplified} to arrive at
\begin{align*}
   LV_4(x,v)
&\le -\left(\frac{p}{2}-2\delta\right) J(x,v)^{p-1}|v|^2
-\left(\frac{1}{2}-2\delta\right) |\grad U(x)+\grad G(v)|^2\\
&\qquad -\left( \frac{\sigma\theta^2}{4}-2\delta
\right) \frac{1}{|v|^{2-\theta}} + C_\delta. 
\end{align*}
Hence, choosing $\delta>0$ sufficiently small, we obtain constants $c,C>0$ such that
\begin{align}\label{eq:LVFour-dTwo-Coercive}
LV_4(x,v) \le -cJ(x,v)^{p-1}|v|^2 -c|\grad U(x)+\grad G(v)|^2 -\frac{c}{|v|^{2-\theta}} + C.
\end{align}

We next show that $V_4\ge1$. By the same argument below~\eqref{eq:LVFour-Coercive}, we have
\begin{align*}
    J(x,v)^p -a_4\log|v|
+\la v,\grad U(x)\ra
+c_0\ge
\bigl(1-2\delta\bigr)J(x,v)^p -C_\delta+c_0.
\end{align*}
Also, since $0<\theta<2<2p$ and
$J(x,v)\ge\frac12|v|^2$, Young's inequality gives $|v|^\theta
\le \delta J(x,v)^p+C_\delta$. Using the these observations, we obtain 
\begin{align}
V_4(x,v)
&=
J(x,v)^p -|v|^\theta -a_4\log|v| +\la v,\grad U(x)\ra +c_0 \nonumber\\
&\ge \bigl(1-2\delta\bigr)J(x,v)^p -C_\delta+c_0 \nonumber
\end{align}
Choosing $\delta>0$ sufficiently small and taking $c_0$ large enough immediately produce $V_4\ge 1$.

Finally, since for all $a>0$
\begin{align*}
    |\log a|^{\frac{p-1}{p}} \le C\Big(a+\frac{1}{a^{2-\theta}}+1\Big),
\end{align*}
we get
\begin{align} \label{eq:LVFour-dTwo-Control-V}
    V_4(x,v)^{\frac{p-1}{p}}
&\le C\left( J(x,v)^p +a_4|\log v| +c_0\right)^{\frac{p-1}{p}}  \notag \\
&\le C\left( J(x,v)^{p-1} + |v|+\frac{1}{|v|^{2-\theta}} +1\right) \notag \\
&\le C\left( 1 +J(x,v)^{p-1}|v|^2 +|\grad U(x)+\grad G(v)|^2 +\frac{1}{|v|^{2-\theta}} \right).
\end{align}
In the above, the last estimate follows from~\eqref{eq:LVFour-dTwo-Remaining-Estimate} and \eqref{ineq:|v|^2<J^(p-1)|v|^2} allowing us to control $|v|$ and $J(x,v)^{p-1}$ by $J(x,v)^p|v|^2$. Combining~\eqref{eq:LVFour-dTwo-Control-V} with~\eqref{eq:LVFour-dTwo-Coercive}, we conclude that there exist constants $c,C>0$ such that
\[
LV_4(x,v) \le -cV_4(x,v)^{\frac{p-1}{p}} + C,
\qquad (x,v)\in\X.
\]
This establishes \eqref{ineq:Lyapunov:singular-damping:polynomial} for the case $d=2$, $\beta=1$ and $n=1$. \vspace{.2cm}

(iii) Suppose now that $d\ge2$ and $\beta>1$, and let $q,p$ satisfy~\eqref{eq:LVFour-q-range} and~\eqref{eq:LVFour-p-range}. 
Using expression~\eqref{eq:LJ:singular} together with Assumption~\ref{cond:G:singular:beta>1}, we have
\begin{align} \label{eq:LVFour-Jp2}
LJ(x,v)^p&
=-pJ(x,v)^{p-1}|v|^2
+pa_4\frac{J(x,v)^{p-1}}{|v|^{\beta-1}}+pa_6\frac{J(x,v)^{p-1}}{|v|^{\betat-1}}  \notag \\
&\qquad
+\frac{\sigma dp}{2}J(x,v)^{p-1}
+\frac{\sigma p(p-1)}{2}J(x,v)^{p-2}|v|^2.
\end{align}
We next compute the generator of the correction in~\eqref{Def:V-Four}. A routine differentiation gives
\begin{align*}
\grad_v|v|^{-(\beta+q-1)}
&=-(\beta+q-1)|v|^{-(\beta+q+1)}v,\\
\triangle_v|v|^{-(\beta+q-1)}
&=-(\beta+q-1)(d-\beta-q-1)|v|^{-(\beta+q+1)},\\
\grad_v\bigl(|v|^{-q}\la v,\grad U(x)\ra\bigr)
&=|v|^{-q}\grad U(x)-q|v|^{-q-2}\la v,\grad U(x)\ra v,\\
\triangle_v\bigl(|v|^{-q}\la v,\grad U(x)\ra\bigr)
&=-q(d-q)|v|^{-q-2}\la v,\grad U(x)\ra.
\end{align*}
Thus, by~\eqref{form:L}, writing $\grad_v|v|^{-(\beta+q-1)} = -(\beta+q-1) |v|^{-q} \frac{v}{|v|^{\beta+1}}$, we obtain
\begin{align*}
    &\frac{a_4}{\beta+q-1}L|v|^{-(\beta+q-1)} \\
    &=|v|^{-q}\Big\la -a_4 \frac{v}{|v|^{\beta+1}},-v-\grad U(x)-\grad G(v)\Big\ra - a_4\frac{\sigma}{2}(d-\beta-q-1)|v|^{-(\beta+q+1)}\\
    &=a_4|v|^{-(\beta+q-1)}+ |v|^{-q}\Big\la a_4 \frac{v}{|v|^{\beta+1}},\grad U(x)+\grad G(v)\Big\ra - a_4\frac{\sigma}{2}(d-\beta-q-1)|v|^{-(\beta+q+1)},
\end{align*}
and similarly,
\begin{align*}
    &\frac{a_6}{\betat+q-1}L|v|^{-(\betat+q-1)}\\
    &=a_6|v|^{-(\betat+q-1)}+|v|^{-q} \Big \la a_6 \frac{v}{|v|^{\betat+1}},\grad U(x)+\grad G(v) \Big \ra  - a_6\frac{\sigma}{2}(d-\betat-q-1)|v|^{-(\betat+q+1)},
\end{align*}
whence, by Assumption~\ref{cond:G:singular:beta>1},
\begin{align}\label{eq:LVFour-Weighted-Singular-Generator}
&L\bigg(\frac{a_4}{\beta+q-1}|v|^{-(\beta+q-1)}+\frac{a_6}{\betat+q-1}|v|^{-(\betat+q-1)}\bigg) \notag \\
&= a_4|v|^{-(\beta+q-1)}+a_6|v|^{-(\betat+q-1)}-|v|^{-q}\la\grad G(v),\grad U(x)+\grad G(v)\ra \nonumber\\
&\qquad-a_4\frac{\sigma}{2}(d-\beta-q-1)|v|^{-(\beta+q+1)}-a_6\frac{\sigma}{2}(d-\betat-q-1)|v|^{-(\betat+q+1)}.
\end{align}
Also,
\begin{align}\label{eq:LVFour-Weighted-Cross-Generator}
&L\bigl(|v|^{-q}\la v,\grad U(x)\ra\bigr) \notag \\
&=|v|^{-q}\la v,\nabla^2U(x)v\ra
-(1-q)|v|^{-q}\la v,\grad U(x)\ra\nonumber\\
&\quad-|v|^{-q}\la\grad U(x),\grad U(x)+\grad G(v)\ra
+q|v|^{-q-2}\la v,\grad U(x)\ra\la v,\grad U(x)+\grad G(v)\ra\nonumber\\
&\quad-\frac{\sigma}{2}q(d-q)|v|^{-q-2}\la v,\grad U(x)\ra.
\end{align}
Combining~\eqref{eq:LVFour-Weighted-Singular-Generator} and~\eqref{eq:LVFour-Weighted-Cross-Generator}, we obtain the following identity; see Appendix~\ref{sec:LVFour-Weighted-Correction-Derivation} for the detailed derivation. 
\begin{align}\label{eq:LVFour-Weighted-Correction}
&L\left(\frac{a_4}{\beta+q-1}|v|^{-(\beta+q-1)}
+\frac{a_6}{\betat+q-1}|v|^{-(\betat+q-1)}
+\frac{1}{1-q}|v|^{-q}\la v,\grad U(x)\ra\right)\nonumber\\
&=-\frac{1}{1-q}|v|^{-q}
\left(|\grad U(x)+\grad G(v)|^2
-q\frac{\la v,\grad U(x)+\grad G(v)\ra^2}{|v|^2}\right)\nonumber\\
&\quad+\frac{1}{1-q}|v|^{-q}\la v,\nabla^2U(x)v\ra
-|v|^{-q}\la v,\grad U(x)+\grad G(v)\ra\nonumber\\
&\quad-\frac{\sigma a_4}{2(1-q)}
\bigl(d-1-\beta(1-q)\bigr)|v|^{-(\beta+q+1)}\nonumber\\
&\quad-\frac{\sigma a_6}{2(1-q)}
\bigl(d-1-\betat(1-q)\bigr)|v|^{-(\betat+q+1)}\nonumber\\
&\quad-\frac{\sigma q(d-q)}{2(1-q)}
|v|^{-q-2}\la v,\grad U(x)+\grad G(v)\ra.
\end{align}
In particular, the coefficient of the leading singular term $|v|^{-(\beta+q+1)}$ is strictly negative precisely when $\beta(1-q)<d-1$, as ensured by~\eqref{eq:LVFour-q-range}.

By Cauchy-Schwarz, the first term on the right-hand side of~\eqref{eq:LVFour-Weighted-Correction} is at most $$-|v|^{-q}|\grad U(x)+\grad G(v)|^2.$$
For the two mixed terms, Young's inequality gives
\begin{align*}
|v|^{-q}|\la v,\grad U(x)+\grad G(v)\ra|
&\le\frac14|v|^{-q}|\grad U(x)+\grad G(v)|^2+|v|^{2-q},\\
\frac{\sigma q(d-q)}{2(1-q)}|v|^{-q-2}|\la v,\grad U(x)+\grad G(v)\ra|
&\le\frac14|v|^{-q}|\grad U(x)+\grad G(v)|^2+C|v|^{-q-2},
\end{align*}
where $C>0$ is a constant. Moreover, since $\beta>1$ and $0\le\betat<\beta$, both $q+2$ and $\betat+q+1$ are positive and strictly smaller than $\beta+q+1$. Thus, Young's inequality yields
\[
C|v|^{-q-2}
+\frac{\sigma|a_6|}{2(1-q)}|d-1-\betat(1-q)|\,|v|^{-(\betat+q+1)}
\le\frac{\sigma a_4}{4(1-q)}\bigl(d-1-\beta(1-q)\bigr)|v|^{-(\beta+q+1)}+C.
\]
Consequently, regardless of the sign of $a_6$,
\begin{align}\label{eq:LVFour-Weighted-Correction-Bound}
&L\left(\frac{a_4}{\beta+q-1}|v|^{-(\beta+q-1)}
+\frac{a_6}{\betat+q-1}|v|^{-(\betat+q-1)}
+\frac{1}{1-q}|v|^{-q}\la v,\grad U(x)\ra\right)\nonumber\\
&\le-\frac12|v|^{-q}|\grad U(x)+\grad G(v)|^2
-\frac{\sigma a_4}{4(1-q)}\bigl(d-1-\beta(1-q)\bigr)|v|^{-(\beta+q+1)}\nonumber\\
&\quad+\frac{1}{1-q}|v|^{-q}|\la v,\nabla^2U(x)v\ra|
+|v|^{2-q}+C.
\end{align}

We next claim that, for every $\delta>0$,
\begin{align}\label{eq:LVFour-Weighted-Hessian}
|v|^{-q}|\la v,\nabla^2U(x)v\ra|
\le\delta J(x,v)^{p-1}|v|^2
+\delta|v|^{-q}|\grad U(x)+\grad G(v)|^2
+\delta|v|^{-(\beta+q+1)}+C_\delta.
\end{align}
Suppose first that $|v|<1$. If $\lambda=2$,~\eqref{cond:U:U''(x)=O(x^lambda)} implies that the Hessian $\nabla^2U(x)$ is bounded, and hence the left-hand side of~\eqref{eq:LVFour-Weighted-Hessian} is bounded by a constant in this case. If $\lambda>2$, the Hessian assumption~\eqref{cond:U:U''(x)=O(x^lambda)} and the bound
$|x|^{\lambda-1}\le C(1+|\grad U(x)|)$, which follows from~\eqref{cond:U:U'>x^(lambda-1)}, give the first inequality below:
\begin{align*}
|v|^{-q}|\la v,\nabla^2U(x)v\ra|
&\le C|v|^{2-q}\left(1+|\grad U(x)|^{\frac{\lambda-2}{\lambda-1}}\right)\\
&\le C|v|^{2-q}
+C|v|^{2-q}|\grad U(x)+\grad G(v)|^{\frac{\lambda-2}{\lambda-1}}
+C|v|^{2-q-\frac{\beta(\lambda-2)}{\lambda-1}}.
\end{align*}
Here the second inequality uses the triangle inequality and the estimate $|\grad G(v)|\le C|v|^{-\beta}$ for $|v|<1$, which follows from Assumption~\ref{cond:G:singular:beta>1}. Since $q<1$ and $|v|<1$, we have $|v|^{2-q}\le1$. Also, Young's inequality with conjugate exponents $2(\lambda-1)/(\lambda-2)$ and $2(\lambda-1)/\lambda$ gives, for every $\eta>0$,
\begin{align*}
C|v|^{2-q}|\grad U(x)+\grad G(v)|^{\frac{\lambda-2}{\lambda-1}}
&\le\eta|v|^{-q}|\grad U(x)+\grad G(v)|^2
+C_\eta|v|^{\frac{4(\lambda-1)}{\lambda}-q}\\
&\le\eta|v|^{-q}|\grad U(x)+\grad G(v)|^2+C_\eta.
\end{align*}
The last inequality holds because $4(\lambda-1)/\lambda-q>0$ and $|v|<1$. If $2-q-\beta(\lambda-2)/(\lambda-1)\ge0$, the corresponding power of $|v|$ is bounded by $1$. Otherwise, since
\[
2-q-\frac{\beta(\lambda-2)}{\lambda-1}
+(\beta+q+1)=3+\frac{\beta}{\lambda-1}>0,
\]
we have $2-q-\beta(\lambda-2)/(\lambda-1)>-(\beta+q+1)$. Therefore, Young's inequality gives, for every $\eta>0$,
\[
C|v|^{2-q-\frac{\beta(\lambda-2)}{\lambda-1}}
\le
\eta|v|^{-(\beta+q+1)}+C_\eta.
\]
Thus, from the above observations, choosing $\eta>0$ sufficiently small in terms of $\delta$ proves~\eqref{eq:LVFour-Weighted-Hessian} when $|v|<1$.

Suppose next that $|v|\ge1$. Fix $\eta>0$ and first consider $|x|\ge1$ and $|v|^2\le\eta|x|^\lambda$. By~\eqref{cond:U:U''(x)=O(x^lambda)} and the bound
$|x|^{\lambda-1}\le C(1+|\grad U(x)|)$ from~\eqref{cond:U:U'>x^(lambda-1)}, we obtain
\begin{align*}
|v|^{-q}|\la v,\nabla^2U(x)v\ra|
\le2a_1|v|^{2-q}|x|^{\lambda-2}
&\le C\eta|v|^{-q}\left(1+|\grad U(x)|^2\right)\\
&\le C\eta|v|^{-q}|\grad U(x)+\grad G(v)|^2+C\eta,
\end{align*}
where the final inequality uses the triangle inequality, $|v|^{-q}\le1$, and $|\grad G(v)|\le C$ since $|v|\ge1$. Thus, choosing $\eta>0$ sufficiently small in terms of $\delta$ yields~\eqref{eq:LVFour-Weighted-Hessian} in the current region. Next, consider the remaining region, where either $|x|<1$ or $|v|^2>\eta|x|^\lambda$. Then, using~\eqref{cond:U:U''(x)=O(x^lambda)}, we obtain
\[
|v|^{-q}|\la v,\nabla^2U(x)v\ra|
\le C_\eta\left(|v|^{2-q}+|v|^{4-q-\frac4\lambda}\right).
\]
Both exponents are positive and strictly less than $2p$: the first follows from $0<q<1$ and $p>1$, while the second follows from~\eqref{eq:LVFour-p-range}, since $p>2-2/\lambda-q/2$. Applying Young's inequality with conjugate exponents $2p/(2-q)$ and $2p/(2p-2+q)$ to the first term, and $2p/(4-q-4/\lambda)$ and $2p/(2p-4+q+4/\lambda)$ to the second, we obtain
\[
C_\eta\left(|v|^{2-q}+|v|^{4-q-\frac4\lambda}\right)
\le\delta\,2^{1-p}|v|^{2p}+C_{\delta,\eta}
\le\delta J(x,v)^{p-1}|v|^2+C_{\delta,\eta}.
\]
The last inequality uses $J(x,v)\ge|v|^2/2$ and $p>1$. This completes the proof of~\eqref{eq:LVFour-Weighted-Hessian}.

Next, since $2-q<2p$, Young's inequality gives
$|v|^{2-q}\le\delta J(x,v)^{p-1}|v|^2+C_\delta$.
Using this together with~\eqref{eq:LVFour-Weighted-Hessian}, with its small parameter chosen after the fixed coefficient $1/(1-q)$, in~\eqref{eq:LVFour-Weighted-Correction-Bound}, we obtain, for every sufficiently small $\delta>0$,
\begin{align}\label{eq:LVFour-Weighted-Correction-Coercive}
&L\left(\frac{a_4}{\beta+q-1}|v|^{-(\beta+q-1)}
+\frac{a_6}{\betat+q-1}|v|^{-(\betat+q-1)}
+\frac{1}{1-q}|v|^{-q}\la v,\grad U(x)\ra\right)\nonumber\\
&\le
2\delta J(x,v)^{p-1}|v|^2
-\left(\frac12-\delta\right)|v|^{-q}|\grad U(x)+\grad G(v)|^2\notag\\
&\quad
-\left(
\frac{\sigma a_4}{4(1-q)}
\bigl(d-1-\beta(1-q)\bigr)-\delta
\right)|v|^{-(\beta+q+1)}
+C_\delta.
\end{align}

We now estimate the positive terms in~\eqref{eq:LVFour-Jp2}. We claim that, for every $\delta>0$,
\begin{align}\label{eq:LVFour-Weighted-Remaining}
\frac{J(x,v)^{p-1}}{|v|^{\beta-1}}+J(x,v)^{p-1}
&\le\delta J(x,v)^{p-1}|v|^2
+\delta|v|^{-q}|\grad U(x)+\grad G(v)|^2
+\delta|v|^{-(\beta+q+1)}+C_\delta.
\end{align}
First observe that the upper bounds on $p$ in~\eqref{eq:LVFour-p-range} imply
\begin{align}\label{eq:LVFour-Weighted-Powers}
0<\frac{\lambda(p-1)}{\lambda-1}<2,
\qquad
\frac{\beta\lambda(p-1)}{\lambda-1}<q+2.
\end{align}
Now, fix $\delta>0$. For $|v|\ge\max\{1,\sqrt{2/\delta}\}$, we have
\[
\frac{J(x,v)^{p-1}}{|v|^{\beta-1}}+J(x,v)^{p-1}
\le2J(x,v)^{p-1}\le\delta J(x,v)^{p-1}|v|^2.
\]
On the other hand, for $1\le|v|<\max\{1,\sqrt{2/\delta}\}$, the upper bound on $U$ in~\eqref{cond:U:U(x)=O(x^lambda+1)} and the bound
$|x|^{\lambda-1}\le C(1+|\grad U(x)|)$, together with $|\grad G(v)|\le C$, give
\[
J(x,v)^{p-1}
\le C_\delta\left(1+|\grad U(x)|^{\frac{\lambda(p-1)}{\lambda-1}}\right)
\le C_\delta\left(1+|\grad U(x)+\grad G(v)|^{\frac{\lambda(p-1)}{\lambda-1}}\right).
\]
Using Young's inequality with conjugate exponents $2(\lambda-1)/ (\lambda(p-1))$ and $2(\lambda-1)/(2(\lambda-1)-\lambda(p-1))$, which are greater than one by the first inequality in~\eqref{eq:LVFour-Weighted-Powers}, we obtain
\[
2J(x,v)^{p-1}\le\eta|\grad U(x)+\grad G(v)|^2+C_{\delta,\eta}
\le\delta|v|^{-q}|\grad U(x)+\grad G(v)|^2+C_{\delta},
\]
where the last inequality follows by choosing
$0<\eta\le\delta\max\{1,\sqrt{2/\delta}\}^{-q}$.
This proves~\eqref{eq:LVFour-Weighted-Remaining} for $|v|\ge1$.

For $|v|<1$, as above, we have
\begin{align*}
J(x,v)^{p-1}
\le C\left(1+|\grad U(x)|^{\frac{\lambda(p-1)}{\lambda-1}}\right)
\le C\left(1+|\grad U(x)+\grad G(v)|^{\frac{\lambda(p-1)}{\lambda-1}}
+|v|^{-\frac{\beta\lambda(p-1)}{\lambda-1}}\right),
\end{align*}
where the final inequality uses the triangle inequality and
$|\grad G(v)|\le C|v|^{-\beta}$. Since $\beta>1$ and $|v|<1$, and using the preceding bound on $J(x,v)^{p-1}$, it follows that
\begin{align*}
&\frac{J(x,v)^{p-1}}{|v|^{\beta-1}}+J(x,v)^{p-1}\\
&\le\frac{2J(x,v)^{p-1}}{|v|^{\beta-1}}\\
&\le C|v|^{-(\beta-1)}
+C|v|^{-(\beta-1)}|\grad U(x)+\grad G(v)|^{\frac{\lambda(p-1)}{\lambda-1}}
+C|v|^{-\left(\beta-1+\frac{\beta\lambda(p-1)}{\lambda-1}\right)}.
\end{align*}
The first and third powers of $|v|^{-1}$ are strictly smaller than $\beta+q+1$, by~\eqref{eq:LVFour-Weighted-Powers}. Therefore, Young's inequality for each of these powers gives, for every $\eta>0$,
\begin{align*}
C|v|^{-(\beta-1)}
+C|v|^{-\left(\beta-1+\frac{\beta\lambda(p-1)}{\lambda-1}\right)}
\le\eta|v|^{-(\beta+q+1)}+C_\eta.
\end{align*}
For the middle term, Young's inequality gives
\begin{align*}
C|v|^{-(\beta-1)}|\grad U(x)+\grad G(v)|^{\frac{\lambda(p-1)}{\lambda-1}}
\le\eta|v|^{-q}|\grad U(x)+\grad G(v)|^2
+C_\eta|v|^{-\frac{2(\beta-1)(\lambda-1)-q\lambda(p-1)}{2(\lambda-1)-\lambda(p-1)}}.
\end{align*}
Using~\eqref{eq:LVFour-Weighted-Powers}, we obtain the first inequality below:
\[
\frac{\lambda(p-1)}{\lambda-1}(\beta+1)
<\frac{(q+2)(\beta+1)}{\beta}<2(q+2),
\]
where the second uses $\beta>1$. Rearranging the above inequality gives
\[
\frac{2(\beta-1)(\lambda-1)-q\lambda(p-1)}{2(\lambda-1)-\lambda(p-1)}
<\beta+q+1.
\]
The denominator on the left is positive by the first inequality in~\eqref{eq:LVFour-Weighted-Powers}. If the numerator is nonpositive, the corresponding power of $|v|$ is bounded by $1$ on $|v|<1$. Otherwise, another application of Young's inequality gives, for every $\eta_1>0$,
\[
C_\eta|v|^{-\frac{2(\beta-1)(\lambda-1)-q\lambda(p-1)}{2(\lambda-1)-\lambda(p-1)}}
\le\eta_1|v|^{-(\beta+q+1)}+C_{\eta,\eta_1}.
\]
Choosing first $\eta>0$ and then $\eta_1>0$ sufficiently small proves~\eqref{eq:LVFour-Weighted-Remaining} also for $|v|<1$.

Next observe that, for $|v|<1$, we have $|v|^{1-\betat}\le|v|^{1-\beta}$, while for $|v|\ge1$, we have $|v|^{1-\betat}\le|v|\le\eta|v|^2+C_\eta$ for every $\eta>0$. Therefore,
\[
\frac{J(x,v)^{p-1}}{|v|^{\betat-1}}
\le
\frac{J(x,v)^{p-1}}{|v|^{\beta-1}}
+\eta J(x,v)^{p-1}|v|^2
+C_\eta J(x,v)^{p-1}.
\]
Using this together with~\eqref{eq:LVFour-Weighted-Remaining} and the bound $|v|^2\le2J(x,v)$ in~\eqref{eq:LVFour-Jp2}, and choosing the corresponding small parameters sufficiently small in terms of $\delta$, we obtain, for every $\delta>0$,
\begin{align} \label{eq:LVFour-Jp2-bound}
LJ(x,v)^p
&\le
-(p-\delta)J(x,v)^{p-1}|v|^2
+\delta|v|^{-q}|\grad U(x)+\grad G(v)|^2
+\delta|v|^{-(\beta+q+1)}
+C_\delta.
\end{align}
Combining~\eqref{eq:LVFour-Weighted-Correction-Coercive} and~\eqref{eq:LVFour-Jp2-bound}, and choosing $\delta>0$ sufficiently small, there exist constants $c,C>0$ such that
\begin{align}\label{eq:LVFour-Weighted-Coercive}
LV_4(x,v)
\le
-c\left(
J(x,v)^{p-1}|v|^2
+|v|^{-q}|\grad U(x)+\grad G(v)|^2
+|v|^{-(\beta+q+1)}
\right)+C.
\end{align}

It remains to prove positivity of $V_4$ and compare the right-hand side of~\eqref{eq:LVFour-Weighted-Coercive} with $V_4^{(p-1)/p}$. By Assumption~\ref{cond:U}, $|v|^2\le2J(x,v)$, and $|x|^\lambda\le a_1J(x,v)+a_1^2$, we have
\begin{align*}
\frac{1}{1-q}|v|^{-q}|\la v,\grad U(x)\ra|
\le C\left(J(x,v)^{\frac{1-q}{2}}
+J(x,v)^{\frac32-\frac1\lambda-\frac q2}\right)
\le\frac14J(x,v)^p+C,
\end{align*}
where the last inequality follows from Young's inequality, since both exponents on the right-hand side are strictly less than $p$. Moreover, for every $\eta>0$, we have
\[
\frac{|a_6|}{|\betat+q-1|}|v|^{-(\betat+q-1)}
\le
\eta J(x,v)^p+\eta|v|^{-(\beta+q-1)}+C_\eta.
\]
Combining the preceding two estimates with the definition of $V_4$ in~\eqref{Def:V-Four}, then choosing $\eta>0$ sufficiently small and $c_0>0$ sufficiently large, we obtain
\begin{align}
V_4(x,v)&\ge1+\frac12J(x,v)^p+\frac{a_4}{2(\beta+q-1)}|v|^{-(\beta+q-1)}, \label{eq:LVFour-Weighted-Comparison:lower}\\
V_4(x,v)&\le C\left(1+J(x,v)^p+|v|^{-(\beta+q-1)}\right). \label{eq:LVFour-Weighted-Comparison}
\end{align}
In particular,~\eqref{eq:LVFour-Weighted-Comparison:lower} implies that $V_4(x,v)\ge1$. Also, since $0<(p-1)/p<1$, the upper bound in~\eqref{eq:LVFour-Weighted-Comparison} implies
\begin{align*}
V_4(x,v)^{\frac{p-1}{p}}
&\le C\left(1+J(x,v)^{p-1}+|v|^{-\frac{(\beta+q-1)(p-1)}p}\right)\\
&\le C\left(1+J(x,v)^{p-1}+|v|^{-(\beta+q+1)}\right)\\
&\le C\left(1+J(x,v)^{p-1}|v|^2
+|v|^{-q}|\grad U(x)+\grad G(v)|^2
+|v|^{-(\beta+q+1)}\right).
\end{align*}
The second inequality follows by considering $|v|\le1$ and $|v|\ge1$ separately and using $(\beta+q-1)(p-1)/p<\beta+q+1$. The third inequality follows from~\eqref{eq:LVFour-Weighted-Remaining} with $\delta=1$. Combining this comparison with~\eqref{eq:LVFour-Weighted-Coercive} yields
\[
LV_4(x,v)\le-cV_4(x,v)^{\frac{p-1}{p}}+C,
\qquad (x,v)\in\X.
\]

Finally, we verify the parameter ranges. The singular diffusion term in~\eqref{eq:LVFour-Weighted-Correction} requires
$q>1-(d-1)/\beta$, which is one of the lower bounds in~\eqref{eq:LVFour-q-range}. The lower bound $2-2/\lambda-q/2$ on $p$ is always smaller than $3-2/\lambda$, and both upper bounds in~\eqref{eq:LVFour-p-range} exceed $1$. Thus the interval for $p$ is nonempty precisely when
\[
2-\frac2\lambda-\frac q2
<1+\frac{(q+2)(\lambda-1)}{\beta\lambda},
\]
which is equivalent to
\[
q>\frac{2\beta(\lambda-2)-4(\lambda-1)}{\beta\lambda+2(\lambda-1)}.
\]
Together with $0<q<1$, this gives~\eqref{eq:LVFour-q-range}. Since $1-(d-1)/\beta<1$, such a $q$ exists precisely when
\[
\beta(\lambda-4)<6(\lambda-1).
\]
This holds for every $\beta>1$ if $2\le\lambda\le4$, and is equivalent to the stated upper bound on $\beta$ if $\lambda>4$. \vspace{.2cm}

2. We now consider the general case $n>1$ to establish~\eqref{ineq:Lyapunov:singular-damping:polynomial}. Using the chain rule, we obtain
\begin{align}\label{eq:LVFour-Powers-Chain-Rule}
L(V_4^n)=nV_4^{n-1}LV_4+\frac{\sigma}{2}n(n-1)V_4^{n-2}|\grad_vV_4|^2.
\end{align}
Next, define
\[
D(x,v):=
\begin{cases}
J(x,v)^{p-1}|v|^2+|\grad U(x)+\grad G(v)|^2+|v|^{-2}, & d\ge3,\ \beta=1,\\[1mm]
J(x,v)^{p-1}|v|^2+|\grad U(x)+\grad G(v)|^2+|v|^{-(2-\theta)}, & d=2,\ \beta=1,\\[1mm]
J(x,v)^{p-1}|v|^2+|v|^{-q}|\grad U(x)+\grad G(v)|^2+|v|^{-(\beta+q+1)}, & \beta>1.
\end{cases}
\]
Then, by~\eqref{eq:LVFour-Coercive},~\eqref{eq:LVFour-dTwo-Coercive}, and~\eqref{eq:LVFour-Weighted-Coercive}, respectively, there exist constants $c,C>0$, independent of $n$, such that
\begin{align}\label{eq:LVFour-Powers-Drift}
LV_4\le-cD+C.
\end{align}
Moreover,~\eqref{eq:LVFour-Control-V},~\eqref{eq:LVFour-dTwo-Control-V}, and the display following~\eqref{eq:LVFour-Weighted-Comparison} give, in the corresponding cases,
\begin{align}\label{eq:LVFour-Powers-Comparison}
V_4^{\frac{p-1}{p}}\le C(1+D).
\end{align}

We first record lower bounds on $V_4$. When $\beta=1$, we write $-\log|v|=(-\log|v|)_+-(\log|v|)_+$. The estimates following~\eqref{eq:LVFour-Coercive} and~\eqref{eq:LVFour-dTwo-Coercive}
give, for every $\delta>0$,
\[
\big|\la v,\grad U(x)\ra\big|
\le \delta J(x,v)^p+C_\delta,
\qquad
(\log|v|)_+
\le \delta J(x,v)^p+C_\delta,
\qquad
|v|^\theta
\le \delta J(x,v)^p+C_\delta,
\]
where the last inequality is used only in the case $d=2$ and $\beta=1$. Substituting these estimates into the definition~\eqref{Def:V-Four} of $V_4$ and choosing
$\delta>0$ sufficiently small, we obtain
\[
V_4(x,v)
\ge
\frac12J(x,v)^p
+a_4(-\log|v|)_+
+c_0-C.
\]
Thus, choosing $c_0$ sufficiently large, independently of $n$, we obtain
\begin{align}\label{eq:LVFour-Powers-Lower-Bound}
V_4(x,v)\ge1+\frac12J(x,v)^p+
\begin{cases}
a_4(-\log|v|)_+, & \beta=1,\\
\dfrac{a_4}{2(\beta+q-1)}|v|^{-(\beta+q-1)}, & \beta>1,
\end{cases}
\end{align}
where the estimate for $\beta>1$ follows from~\eqref{eq:LVFour-Weighted-Comparison:lower}. Suppose that, in all three cases, there exists a function $h:[1,\infty)\to(0,\infty)$ such that
\begin{align}\label{eq:LVFour-Powers-Gradient-Control}
\frac{|\grad_vV_4|^2}{V_4}\le h(V_4)D, \qquad h(t)\longrightarrow0 \quad\text{as }t\longrightarrow\infty,
\end{align}
where $h$ is independent of $n$. Since $h(t)\to0$, for each $n\ge1$ we may choose $R_n\ge1$ sufficiently large so that
\[
\frac{\sigma}{2}(n-1)h(t)\le \frac{c}{2},
\qquad t\ge R_n,
\]
where $c$ is the constant in~\eqref{eq:LVFour-Powers-Drift}. Consequently, using~\eqref{eq:LVFour-Powers-Drift} and~\eqref{eq:LVFour-Powers-Gradient-Control}, we obtain the inequality below, on $\{V_4\ge R_n\}$
\begin{align*}
L(V_4^n)
&\overset{\eqref{eq:LVFour-Powers-Chain-Rule}}{=}nV_4^{n-1}LV_4+nV_4^{n-1}
\left[
\frac{\sigma}{2}(n-1)
\frac{|\nabla_vV_4|^2}{V_4}
\right] \\
&\le
-ncV_4^{n-1}D+nCV_4^{n-1}
+\frac{nc}{2}V_4^{n-1}D\\
&=
-\frac{nc}{2}V_4^{n-1}D+nCV_4^{n-1}.
\end{align*}
On the other hand, on the compact set $\{V_4\le R_n\}$, since $L(V_4^n)$, $V_4^{n-1}$, and $D$ are continuous, there exists a constant $C_n>0$ such that
\[
L(V_4^n)+\frac{nc}{2}V_4^{n-1}D\le C_n \le C_nV_4^{n-1},
\]
where the second inequality follows from $V_4^{n-1}\ge1$. Combining this with the estimate on $\{V_4\ge R_n\}$ and enlarging $C_n$ if necessary, we obtain, for $(x,v)\in\X$,
\[
L(V_4^n)
\le
-\frac{nc}{2}V_4^{n-1}D
+C_nV_4^{n-1}
\le
-a_nV_4^{n-\frac1p}
+b_nV_4^{n-1},
\]
where $a_n=\frac{nc}{2C}$ and $b_n=\frac{nc}{2}+C_n$, and the second inequality uses that $D\ge C^{-1}V_4^{(p-1)/p}-1$, which follows from~\eqref{eq:LVFour-Powers-Comparison}. Finally, $p>1$ implies $n-1<n-1/p$, and therefore Young's inequality gives
\[
b_nV_4^{n-1}\le\frac{a_n}{2}V_4^{n-\frac1p}+C_n.
\]
Substituting this into the preceding estimate proves the assertion. It only remains to prove~\eqref{eq:LVFour-Powers-Gradient-Control}, which we do by considering the following three cases. \vspace{.2cm}

(i) Suppose that $d\ge3$ and $\beta=1$. Differentiating~\eqref{Def:V-Four}, we obtain
\[
\grad_vV_4=pJ^{p-1}v+\grad U(x)+\grad G(v)
-\left(\grad G(v)+a_4\frac{v}{|v|^2}\right).
\]
By Assumption~\ref{cond:G:singular}, the squared norm of the last term is bounded above by $CD$. Consequently,
\begin{align*}
|\grad_vV_4|^2
&\le3p^2J^{2p-2}|v|^2+3|\grad U(x)+\grad G(v)|^2
+3\left|\grad G(v)+a_4\frac{v}{|v|^2}\right|^2\\
&\le C(J^{p-1}+1)D.
\end{align*}
Using $J^p\le2V_4$ from~\eqref{eq:LVFour-Powers-Lower-Bound}, we obtain~\eqref{eq:LVFour-Powers-Gradient-Control} with $h(t)=C(t^{-1/p}+t^{-1})$. \vspace{.2cm}

(ii) Suppose that $d=2$ and $\beta=1$. In this case,~\eqref{Def:V-Four} yields
\begin{align*}
|\grad_vV_4|^2
&\le4p^2J^{2p-2}|v|^2+4\theta^2|v|^{2\theta-2}
+4|\grad U(x)+\grad G(v)|^2\\
&\quad+4\left|\grad G(v)+a_4\frac{v}{|v|^2}\right|^2.
\end{align*}
The right-hand side is bounded by $C(J^{p-1}+1)D$, where for the last term we use Assumption~\ref{cond:G:singular} and $\theta<2(1-\betat)$. Hence, applying~\eqref{eq:LVFour-Powers-Lower-Bound}, we obtain~\eqref{eq:LVFour-Powers-Gradient-Control} with $h(t)=C(t^{-1/p}+t^{-1})$. \vspace{.2cm}

(iii) Suppose that $d\ge2$ and $\beta>1$. Differentiating the two radial terms in~\eqref{Def:V-Four} and using Assumption~\ref{cond:G:singular:beta>1}, we obtain $|v|^{-q}\grad G(v)$. Using this together with $v\,\la v,\grad U(x)\ra=(v\otimes v)\grad U(x)$, we obtain
\begin{align*}
\grad_vV_4
&=pJ^{p-1}v+|v|^{-q}\grad G(v)
+\frac{1}{1-q}|v|^{-q}
\left(I_d-q\frac{v\otimes v}{|v|^2}\right)\grad U(x)\\
&=pJ^{p-1}v+\frac{1}{1-q}|v|^{-q}
\left(I_d-q\frac{v\otimes v}{|v|^2}\right)\bigl(\grad U(x)+\grad G(v)\bigr).
\end{align*}
The second equality follows because $\grad G(v)$ is parallel to $v$. Consequently,
\[
|\grad_vV_4|^2
\le2p^2J^{2p-2}|v|^2+\frac{2}{(1-q)^2}|v|^{-2q}|\grad U(x)+\grad G(v)|^2
\le C\left(J^{p-1}+|v|^{-q}\right)D.
\]
The above together with~\eqref{eq:LVFour-Weighted-Comparison:lower} implies~\eqref{eq:LVFour-Powers-Gradient-Control} with $h(t)=C(t^{-1/p}+t^{q/(\beta+q-1)-1})$.
\end{proof}

\section{Ergodicity and Mixing} \label{sec:ergodicity-mixing}

In this section, we proceed to establish Theorem \ref{thm:mixing}. The first two auxiliary results in Lemma \ref{lem:Hormander} and Lemma \ref{lem:control} provide the hypoellipticity and the solvability of the associated control problem for \eqref{eq:Langevin-SDE}. Together with the Lyapunov functions from Section \ref{sec:ergodicity}, they furnish the proof for Theorem \ref{thm:mixing} part 1 and part 3. Section \ref{sec:ergodicity-mixing} concludes with Lemma \ref{lem:dimension-1:no-IPM} establishing the non-existence of invariant measure in the case of singular damping in dimension $d=1$, as stated in Theorem \ref{thm:mixing}, part 2. 

We start with the verification of hypoellipticity through the Hormander's condition and note that the operator $L$ from \eqref{form:L} can be recast as
\begin{align*}
\L_N =  X_0 + (X_1)^2,
\end{align*}
where $X_0,X_1$ denote the vector fields given by
\begin{align*}
X_0& = v\partial_{x}+ \big( - v -\grad U(x)  -\grad G(v)\big)\partial_{v}, \notag \\
X_1 &= \partial_{v}.
\end{align*}
Following \cite[Theorem 1.1]{hormander1967hypoelliptic} and \cite[Condition (159)]{bellet2006ergodic}, we have the following auxiliary result asserting that the Lie algebra generated by $X_0,X_1$ has maximal rank in $\X$.

\begin{lemma} \label{lem:Hormander}
     Suppose that $U$ satisfies Assumption~\ref{cond:U} and that $G$ satisfies either Assumption \ref{cond:G:poly} or Assumption \ref{cond:G:singular}. Then, the vector fields $X_0,X_1$ satisfy the Hormander's condition in $\X$. That is, the Lie algebra generated by the family
\begin{align*}
X_1,\quad\{[X_i,X_j]\}_{i,j=0}^1,\quad \{[[X_i,X_j],X_k]\}_{i,j,k=0}^1,\quad\dots,
\end{align*}
has maximal rank at every point $(x,v)\in\X$.
\end{lemma}
The proof of Lemma \ref{lem:Hormander} is standard and can be found in literature, e.g, \cite{mattingly2002ergodicity,pavliotis2014stochastic}. Since the argument is relatively short, we include it here for the sake of completeness.

\begin{proof}[Proof of Lemma \ref{lem:Hormander}]
    Observe that 
    \begin{align*}
        [X_0,X_1]  = -\partial_x + \partial_v+\grad^2G(v)\partial_v,
    \end{align*}
    whence
    \begin{align*}
        \partial_x = \partial_v+\grad^2G(v)\partial_v-[X_0,X_1].
    \end{align*}
    It follows that the pair $\{x_1,[X_0,X_1]\}$ has maximal rank at every point $(x,v)\in \X$, as claimed.
\end{proof}

Next, we turn to the second ingredient to establish Theorem \ref{thm:mixing} concerning the solvability of the control problem associated with \eqref{eq:Langevin-SDE}. Namely, consider the following system of ODEs:
\begin{align} \label{eq:Langevin-SDE:control-problem}
\d x &= v\d t,\notag\\
\d v(t) & = -v\d t -\grad U(x)\d t-  \grad G\big(v\big) \d t +\sqrt{\sigma} \,\d \Gamma,
\end{align}
where $\Gamma(\cdot)\in C^1(\X;\rbb)$ is a control path. We then have the following condition requiring that the element $X_1\in \X=\rbb^d \times\D$ be reachable through \eqref{eq:Langevin-SDE:control-problem} from $X_0\in \X$.

\begin{definition} \label{def:control-problem}
Given any $X_0,X_1\in\X$, the control problem \eqref{eq:Langevin-SDE:control-problem} is said to be solvable if there exist a time $T>0$ and a path $\Gamma$ such that the solution $Z(t)$ of \eqref{eq:Langevin-SDE:control-problem} satisfies
\begin{align*}
Z(0)= X_0,\quad \textup{and}\quad Z(T)=X_1.
\end{align*}
\end{definition}

Following standard arguments in literature \cite{hairer2009hot,mattingly2002ergodicity,
ottobre2011asymptotic,pavliotis2014stochastic,bellet2006ergodic}, verifying the solvability of \eqref{eq:Langevin-SDE:control-problem} as in Definition \ref{def:control-problem} is relatively not difficult. This result is summarized below through Lemma \ref{lem:control}.

\begin{lemma} \label{lem:control}
   Suppose that $U$ satisfies Assumption~\ref{cond:U} and that 
   \begin{itemize}
       \item[\textup{(i)}] in dimension $d=1$, $G$ satisfies Assumption \ref{cond:G:poly}; and
       \item[\textup{(ii)}] in dimension $d\ge 2$, $G$ satisfies either Assumption \ref{cond:G:poly} or Assumption \ref{cond:G:singular}.
   \end{itemize}
Then, the associated control problem \eqref{eq:Langevin-SDE:control-problem} is solvable in the sense of Definition \ref{def:control-problem}.
\end{lemma}
\begin{proof} In order to construct a triple $(x,v,\Gamma)$ solving \eqref{eq:Langevin-SDE:control-problem}, there are two cases to be considered depending on the choice of $G$. We start with the simpler case where $G$ is a polynomial.

Case 1: $G$ satisfies Assumption \ref{cond:G:poly}. In this case, we simply pick an arbitrary function $\f\in C^2(\rbb;\rbb^d)$ such that
\begin{align*}
    \f(0)=x_0,\quad \f(1)=x_1,\quad \f'(0)=v_),\quad \f'(1)=v_1.
\end{align*}
Then, setting $T=1$ and
\begin{align} \label{form:Gamma}
    \Gamma(t) = \frac{1}{\sqrt{\sigma}}\int_0^t \f'(s)+\grad U(\f(s))+\grad G(\f'(s))+\f''(s)\d s,
\end{align}
observe that the triple $(\f,\f',\Gamma)$ solves the control problem \eqref{eq:Langevin-SDE:control-problem} on the interval $[0,1]$.

Case 2: $G$ satisfies Assumption \ref{cond:G:poly} in dimension $d\ge 2$. We note that the main difference of this situation from the previous case is to ensure that the choice of $\f$ satisfies $\f'(t)\neq \vec{0}$, $t\in[0,1]$, owing to the fact that the phase space in the singular damping setting is $\X=\rbb^d\times (\rbb^d\setminus\{0\})$. For this purpose, let 
\begin{align*}
    a = (0,\dots,0,a_d),\quad b = (0,\dots,0,b_d)
\end{align*}
be given elements in $\rbb^d$ and be specified later. Observe that the function $\f$ defined as
\begin{align*}
    \f(t) &= at^5+bt^4+ 3\big[v_1+v_0-2(x_1-x_0)-3a-2b\big]t^2\\
    &\qquad+\big[3(x_1-x_0)-(v_1+2v_0)+2a+b\big]t^2+v_0t+x_0
\end{align*}
satisfies
\begin{align*}
    \f(0)=x_0,\quad \f(1)=x_1,\quad \f'(0)=v_0,\quad \f'(1)=v_1.
\end{align*}
It remains to carefully tune $a_d,b_d\in\rbb$ so as to guarantee $\f'(t)\neq 0$, $t\in[0,1]$. To this end, we recast $\f'(t)$ as
\begin{align*}
    \f'(t) &= at(5t^3-9t+4)+2bt(2t^2-3t+1)+3\big[v_1+v_0-2(x_1-x_0)  \big]t^2\\
    &\qquad+2\big[ 3(x_1-x_0)-(v_1+2v_0)  \big]t+v_0.
\end{align*}
In particular, denoting by $[a]_i$ the $i$-component of vector $a$, for $i=1,\dots,d-1$, 
\begin{align*}
    [\f'(t)]_i &= 3\big[v_1+v_0-2(x_1-x_0)  \big]_it^2+2\big[ 3(x_1-x_0)-(v_1+2v_0)  \big]_it+[v_0]_i,
\end{align*}
which depends on $x_0,x_1,v_0,v_1$ but is independent of the choice of $a_d,b_d$. So, on the one hand, if there does not exist $t^*\in(0,1)$ such that 
\begin{align*}
       [\f'(t^*)]_i = 0,\quad i=1,\dots,d-1,
\end{align*}
we deduce immediately that $(\f,\f',\Gamma)$ where $\Gamma$ is defined in \eqref{form:Gamma} is a solution for the control problem \eqref{eq:Langevin-SDE:control-problem}. On the other hand, if such a $t^*$ exists (depending only on $x_0,x_1,v_0,v_1$), we note that the number of such $t^*$ is finite. Also, it is clear that at least either $5(t^*)^3-9t^*+4\neq 0$ or $2(t^*)^2-3t^*+1\neq 0$. It follows that we may pick $a_d,b_d$ such that for all such $t^*$, the $d$-component of $\f'(t^*)$ is non zero, i.e.,
\begin{align*}
    [\f'(t^*)]_d &= a_dt^*(5(t^*)^3-9t^*+4)+2b_dt^*(2(t^*)^2-3t^*+1)+3\big[v_1+v_0-2(x_1-x_0)  \big]_d(t^*)^2\\
    &\qquad+2\big[ 3(x_1-x_0)-(v_1+2v_0)  \big]_dt^*+[v_0]_d\\
    &\neq 0.
\end{align*}
Once again, it allows for choosing $(\f,\f',\Gamma)$ as a solution for the control problem \eqref{eq:Langevin-SDE:control-problem}. The proof is thus complete.

\end{proof}

The last ingredient to prove Theorem \ref{thm:mixing} is the non-existence of invariant measure when $G$ is repulsive in dimension $d=1$. More precisely, we have the following result.

\begin{lemma} \label{lem:dimension-1:no-IPM}
    Suppose that $U$ satisfies Assumption~\ref{cond:U} and that $G$ satisfies Assumption \ref{cond:G:singular}. Then, in dimension $d=1$, there does not exist an invariant probability measure $\nu$ for $P_t$.
\end{lemma}

For the sake of clarity, the proof of Lemma \ref{lem:dimension-1:no-IPM} is deferred by a few lines. Assuming the above result, let us conclude Theorem \ref{thm:mixing} by verifying the hypothesis of \cite[Theorem 3.5]{hairer2009hot}. Since the argument is short, we include it here for completeness.

\begin{proof}[Proof of Theorem \ref{thm:mixing}]
First of all, the functions $V_3$ and $V_4$ respectively constructed in Lemma \ref{lem:ergodicity:poly-damping} and Lemma \ref{lem:ergodicity:singular-damping} corresponding to polynomial damping and singular damping provide the Lyapunov condition required in \cite[Theorem 3.5]{hairer2009hot}. Also, the results of Lemma \ref{lem:Hormander} and Lemma \ref{lem:control} respectively verify \cite[Assumption 2]{hairer2009hot} and \cite[Assumption 3]{hairer2009hot}. Altogether, we invoke \cite[Theorem 3.5]{hairer2009hot} to deduce the uniqueness of the invariant measure $\nu$ and the mixing rates as stated in Theorem \ref{thm:mixing}, parts 1 and 3. Lastly, Lemma \ref{lem:dimension-1:no-IPM} furnishes the argument for Theorem \ref{thm:mixing}, part 2.

\end{proof}

Finally, we provide the proof of Lemma \ref{lem:dimension-1:no-IPM}, which is the main ingredient to establish Theorem \ref{thm:mixing}, part 2.

\begin{proof}[Proof of Lemma \ref{lem:dimension-1:no-IPM}]
    By contradiction, suppose that there exists at least one invariant probability measure $\nu$ for $P_t$. Without loss of
generality, we may assume that $\nu$ is ergodic by ergodic decomposition,  (see e.g., \cite[Chapter
1.3]{cornfeld2012ergodic} and  \cite[Theorem 5.1.3]{viana2016foundations}). Then, there exists a set $A\subset \X$ such that $\nu(A)=1$ and that for all initial condition $(x_0,v_0)\in A$ and bounded function $\f:\X\to\rbb$, it holds that
\begin{align*}
    \E \f(x(t),v(t))\to \int_\X f(x,v)\nu(\d x,\d v).
\end{align*}
Denote by $A_1$ the projection of $A$ on the $x$-marginal, i.e., 
\begin{align*}
    A_1= \{x\in\rbb:\exists\, v , (x,v)\in A\}.
\end{align*}
There are two cases to be considered depending on $\sup A_1$. 

Case 1: $\sup A_1 =\infty$. In this case, we may extract a sequence of initial conditions $\{(x_n,v_n)\}_{n\ge 1}\subset A$ such that $\{x_n\}_{n\ge 1}\subset A_1$ diverging to infinity. Observe that for each such initial condition $(x_n,v_n)$,  
\begin{align*}
   \E\big[ \mathbf{1}_{(-\infty,x_n]\times \D}\big(x(t),v(t)\big)\big] = \P(x(t)\le x_n) \to \nu\big((-\infty,x_n]\times \D\big).
\end{align*}
Since $v(t)>0$ a.s., we get $x(t)> x(0)=x_n$ a.s., implying that $\P(x(t)\le x_n)=0$. It follows immediately from the above limiting behavior that
$$\nu((-\infty,x_n]\times \D)=0.$$ 
Sending $n\to\infty$, we obtain by virtue of the Monotone Convergence Theorem that $\nu(\X)=0$, a contradiction.

Case 2: $\sup A_1 = c<\infty$. In this case, on the one hand, we get $$\nu((c,\infty)\times\D)=0.$$ 
On the other hand, there exists a sequence of initial conditions $\{(x_n,v_n)\}_{n\ge 1}\subset A$ such that $\{x_n\}_{n\ge 1}\subset A_1$ is increasing and that $x_n\to c$, as $n\to\infty$. Similar to Case 1, we also deduce that $\nu((-\infty,x_n]\times \D)=0$, implying that $$\nu((-\infty,c)\times \D)=0.$$
As consequence, we obtain
\begin{align*}
    \nu(\{c\}\times \D)=0.
\end{align*}
By invariance property, if $(x(0),v(0))\sim \nu$, it holds that $(x(t),v(t))\sim \nu$. In particular, we have
\begin{align*}
    \P(x(0)=c)=1 =\P(x(t)=c),
\end{align*}
contradicting the fact that in dimension $d=1$, $x(t)>x(0)$ a.s.

From both cases, we conclude that in dimension $d=1$, there does not exist any invariant probability measure for $P_t$, as claimed.

\end{proof}

\appendix

\section{Existence and finite-time explosion of the comparison diffusion} \label{sec:illposed:auxiliary-equation}

Fix $k>0$, $\sigma>2k$, $C\ge0$, and $x>0$. Let $W$ be a standard
one-dimensional Brownian motion on a filtered probability space
$(\Omega,\mathcal F,\{\mathcal F_t\}_{t\ge0},\mathbb P)$ satisfying the usual conditions. Consider the stochastic differential equation on the state space $(0,\infty)$
\begin{align}\label{eq:Comparison-Diffusion}
\d u(t)
=
\left( -u(t)+C+\frac{k}{u(t)} \right)\d t + \sqrt{\sigma}\,\d W(t), \qquad u(0)=x.
\end{align}
We use the notion of a solution up to an explosion time in the sense of~\cite[Def.~5.5.20, pp.~343--344]{karatzas2012brownian}. Since the state space in the present setting is $(0,\infty)$, we use the localization intervals $(1/n,n)$. More precisely, whenever $u$ is a solution of~\eqref{eq:Comparison-Diffusion} up to an explosion time, choose $n_0\ge1$ such that $1/n_0<x<n_0$ and define the successive exit times
\begin{align}\label{eq:Comparison-Exit-Times}
S_n
:=
\inf\left\{ t\ge0: u(t)\notin\left(\frac1n,n\right) \right\},
\qquad n\ge n_0.
\end{align}
Since the intervals $(1/n,n)$ increase with $n$, the sequence
$\{S_n\}_{n\ge n_0}$ is increasing. We define
\begin{align}\label{eq:Comparison-Explosion-Time}
S
:=
\lim_{n\to\infty}S_n \in(0,\infty]
\end{align}
and refer to $S$ as the \textit{explosion time}.

\begin{lemma} \label{Lem:Comparison-Diffusion}
For $k>0$ and $\sigma>2k$, the following statements hold.
\begin{enumerate}
\item[\textup{(i)}] There exists a unique strong solution $\{u(t)\}_{ t\in[0,S)}$ to~\eqref{eq:Comparison-Diffusion} taking values in $(0,\infty)$ up to the explosion time $S$ defined in~\eqref{eq:Comparison-Explosion-Time}.

\item[\textup{(ii)}] The explosion time is finite almost surely and the explosion occurs at the left endpoint $0$; more precisely,
\[
S<\infty,
\qquad
\lim_{t\uparrow S}u(t)=0,
\qquad
\sup_{0\le t<S}u(t)<\infty
\qquad
\textup{a.s.}
\]
\end{enumerate}
\end{lemma}

\begin{remark}
In view of Lemma~\ref{Lem:Comparison-Diffusion}\textup{(ii)}, the solution admits a continuous extension to the explosion time by setting
\[
u(S):=0.
\]
If we further stipulate that $u(t)=u(S)=0$ for $t\ge S$, then, in analogy with~\cite[Eq.~(5.6), p.~330]{karatzas2012brownian},
\[
S
=
\inf\left\{ t\ge0: u(t)\notin(0,\infty) \right\}
=
\inf\left\{ t\ge0: u(t)\in\{0,+\infty\} \right\}
=
\inf\left\{ t\ge0: u(t)=0 \right\}
\quad
\textup{a.s.}
\]
\end{remark}

\begin{proof}[Proof of Lemma~\ref{Lem:Comparison-Diffusion}]
 
Part (i). Recall that $n_0\ge1$ such that $1/n_0<x<n_0$. For each $n\ge n_0$, set $I_n:=(1/n,n)$
and define the projection map $\pi_n:\mathbb R\longrightarrow [1/n,n]$ by
\[
\pi_n(r)
=
\begin{cases}
\frac{1}{n}, & r<\frac{1}{n},\\[1mm]
r, & \frac{1}{n}\le r\le n,\\[1mm]
n, & r>n.
\end{cases}
\]
We extend the drift coefficient $b(r):=-r+C+k/r$ from $I_n$ to a
globally Lipschitz function $b_n:\mathbb R\to\mathbb R$ by setting
\[
b_n(r) := b\bigl(\pi_n(r)\bigr) = -\pi_n(r)+C+\frac{k}{\pi_n(r)}, \qquad r\in\mathbb R,
\]
and keep the diffusion coefficient equal to $\sqrt{\sigma}$. The
resulting coefficients are bounded and globally Lipschitz, so~\cite[Prop.~5.5.17, p.~341]{karatzas2012brownian} yields a unique strong solution $\{u^{(n)}(t)\}_{t\in[0,\infty)}$ of
\[
\d u^{(n)}(t)
=
b_n\bigl(u^{(n)}(t)\bigr)\,\d t
+
\sqrt{\sigma}\,\d W(t),
\qquad
u^{(n)}(0)=x.
\]
Since $b_n=b$ on $I_n$, the process $u^{(n)}$ satisfies the original equation~\eqref{eq:Comparison-Diffusion} on $[0,\eta_n)$, where
\[
\eta_n
:=
\inf\left\{ t\in[0,\infty): u^{(n)}(t)\notin I_n \right\}.
\]

If $m>n$, define
\[
\eta_{m,n}
:=
\inf\left\{ t\in[0,\infty): u^{(m)}(t)\notin I_n \right\}.
\]
Then the processes $u^{(n)}$ and $u^{(m)}$ satisfy the same
equation~\eqref{eq:Comparison-Diffusion}, with the same initial
condition and the same Brownian motion, on $[0,\eta_n\wedge\eta_{m,n})$. By strong uniqueness,
\[
u^{(m)}(t)=u^{(n)}(t),
\qquad
0\le t<\eta_n\wedge\eta_{m,n},
\]
almost surely. By continuity, the two processes therefore exit $I_n$ at the same time, so
\[
\eta_{m,n}=\eta_n
\qquad
\textup{a.s.}
\]
Since $I_n\subset I_m$, we have $\eta_{m,n}\le\eta_m$, and hence
$\eta_n\le\eta_m$. Thus,
\[
\eta:=\lim_{n\to\infty}\eta_n\in(0,\infty]
\]
is well defined, and the localized solutions can be patched together by setting
\[
u(t):=u^{(n)}(t)
\qquad
\text{whenever }t<\eta_n.
\]
The resulting process $\{u(t)\}_{t\in[0,\eta)}$ is therefore a unique strong solution of~\eqref{eq:Comparison-Diffusion} taking values in $(0,\infty)$. Moreover, for this process,
\[
\eta_n
=
\inf\left\{ t\ge0: u(t)\notin\left(\frac{1}{n},n\right) \right\}
=
S_n,
\]
where $S_n$ is defined in~\eqref{eq:Comparison-Exit-Times}. Hence $S=\eta$ is precisely the explosion time defined
in~\eqref{eq:Comparison-Explosion-Time}. \vspace{.2cm}

Part (ii). The scale function (see~\cite[Eq.~(5.42), p.~339]{karatzas2012brownian} with reference point $c=1$) associated with the diffusion~\eqref{eq:Comparison-Diffusion} is given by
\begin{align}\label{eq:Comparison-Scale-Function}
p(r)
:=
\int_1^r \exp\left( -\int_1^y\frac{2b(z)}{\sigma}\,\d z \right)
\d y, \qquad r>0,
\end{align}
where $b(r):=-r+C+k/r$ for $r>0$ is the drift coefficient. Set $c:=\exp\left(\frac{2C-1}{\sigma}\right)>0$. Then, by the Fundamental Theorem of Calculus, the derivative of the scale function is
\begin{align}\label{eq:Comparison-Scale-Density}
p'(r)
=
\exp\left( -\int_1^r\frac{2b(z)}{\sigma}\,\d z \right)
=
c\,r^{-\frac{2k}{\sigma}} \exp\left( \frac{r^2-2Cr}{\sigma} \right), \qquad r>0,
\end{align}
where the second equality holds since
\begin{align*}
\int_1^r \frac{2b(z)}{\sigma}\,\d z
&=
\frac{2}{\sigma} \int_1^r \left( -z+C+\frac{k}{z} \right)\d z \\
&=
\frac{2}{\sigma} \left[ -\frac{z^2}{2}+Cz+k\log z \right]_{z=1}^{z=r}
=
-\frac{r^2}{\sigma} +\frac{2Cr}{\sigma} +\frac{2k}{\sigma}\log r -\frac{2C-1}{\sigma}.
\end{align*}
Since $C\ge0$, for $0<z<1$ we have $z^2-2Cz\le1$. This yields the first inequality below.
\begin{align} \label{eq:finiteness-Integral-pDensity}
\int_0^1p'(z)\,\d z
&\overset{\eqref{eq:Comparison-Scale-Density}}{=}
c\int_0^1 z^{-\frac{2k}{\sigma}} \exp\left( \frac{z^2-2Cz}{\sigma} \right) \d z \nonumber \\
&\le
c\exp\left( \frac1{\sigma} \right) \int_0^1 z^{-\frac{2k}{\sigma}} \,\d z
=
c\exp\left( \frac1{\sigma} \right) \frac{1}{1-\frac{2k}{\sigma}}
<\infty, 
\end{align}
where the last inequality follows from $\sigma>2k$. Hence,
\begin{align}\label{eq:Comparison-Scale-Zero-Limit}
p(0+)
:=
\lim_{r\downarrow0}p(r)
=
\lim_{r\downarrow0}\int_1^r p'(z)\,\d z
=
-\lim_{r\downarrow0}\int_r^1p'(z)\,\d z
=
-\int_0^1p'(z)\,\d z
>
-\infty,
\end{align}
where the second last equality uses the definition of the improper integral, and the final inequality follows from the finiteness established in~\eqref{eq:finiteness-Integral-pDensity}.

On the other hand, since $r^{-\frac{2k}{\sigma}}
\exp\big(\frac{r^2}{4\sigma}\big)\longrightarrow\infty$ as $r\to\infty$, there exists $R_0>0$ such that $r^{-\frac{2k}{\sigma}} \exp\big(\frac{r^2}{4\sigma}\big)\ge1$ for all $r\ge R_0$. Fix $R>\max\{1,2C,R_0\}$. Then, for every $r\ge R$,
\[
r^{-\frac{2k}{\sigma}} \exp\left( \frac{r^2}{4\sigma} \right) \ge 1.
\]
For $r>R$, the above implies the second inequality below.
\begin{align}\label{eq:Comparison-Scale-Infinity-Limit}
p(r)
=&
\int_1^r p'(y)\,\d y \nonumber\\
\overset{\eqref{eq:Comparison-Scale-Density}}{=}&
p(R)
+
c\int_R^r y^{-\frac{2k}{\sigma}} \exp\left( \frac{(y-C)^2-C^2}{\sigma} \right) \,\d y \nonumber\\
\ge&
c e^{-\frac{C^2}{\sigma}} \int_R^r y^{-\frac{2k}{\sigma}} \exp\left( \frac{y^2}{4\sigma} \right) \,\d y
\ge c e^{-\frac{C^2}{\sigma}} (r-R) \longrightarrow\infty \qquad\text{as }r\to\infty.
\end{align}
The first inequality uses that $p(R)>0$ and
$(y-C)^2\ge y^2/4$ for $y\ge R>2C$.

Since the diffusion coefficient $\sqrt{\sigma}$ is strictly positive and the drift $b(r)=-r+C+k/r$ is smooth on $I:=(0,\infty)$, the conditions {\rm(ND)}$'$ and {\rm(LI)}$'$ in~\cite[p.~343]{karatzas2012brownian} are satisfied. Moreover, by part~\textup{(i)}, $u$ is a weak solution on
$I$ with the nonrandom initial condition $u(0)=x\in I$ and explosion time $S$. Also, by~\eqref{eq:Comparison-Scale-Zero-Limit} and~\eqref{eq:Comparison-Scale-Infinity-Limit}, respectively, $p(0+)>-\infty$ and $p(\infty)=\infty$. Hence, with $\ell=0$ and $r=\infty$,~\cite[Prop.~5.5.22\textup{(b)}, p.~345]{karatzas2012brownian} yields
\[
\mathbb P\left[ \lim_{t\uparrow S}u(t)=0 \right]
=
\mathbb P\left[ \sup_{0\le t<S}u(t)<\infty \right]
=
1.
\]

To prove that the explosion occurs in finite time, let $v$ be the function defined in~\cite[Eq.~(5.65), p.~347]{karatzas2012brownian}, which in the present setting (with reference point $c=1$) is given, for $0<r<1$, by
\begin{align} 
v(r)
:=
\int_1^r p'(y) \left( \int_1^y \frac{2}{\sigma p'(z)} \,\d z \right) \d y
=
\int_r^1 p'(y) \left( \int_y^1 \frac{2}{\sigma p'(z)} \,\d z \right) \d y. \nonumber
\end{align}
Here, the second equality follows by reversing the limits of both integrals. Using~\eqref{eq:Comparison-Scale-Density}, we have the equality below.
\begin{align*}
\int_0^1 \frac{2}{\sigma p'(z)} \,\d z 
=
\frac{2}{\sigma c} \int_0^1 z^{\frac{2k}{\sigma}} \exp\left( \frac{-z^2+2Cz}{\sigma} \right) \d z
\le
\frac{2e^{2C/\sigma}}{\sigma c} \int_0^1 z^{\frac{2k}{\sigma}} \,\d z
<\infty
\end{align*}
The first inequality follows from $-z^2+2Cz\le2C$ for $0<z<1$. The above, together with~\eqref{eq:finiteness-Integral-pDensity}, implies the finiteness below.
\begin{align*}
v(0+)
:=
\lim_{r\downarrow0}v(r)
&=
\int_0^1 p'(y) \left( \int_y^1 \frac{2}{\sigma p'(z)} \,\d z \right) \d y 
\le
\left( \int_0^1 \frac{2}{\sigma p'(z)} \,\d z \right) \left( \int_0^1 p'(y) \,\d y \right) <\infty,
\end{align*}
where the inequality follows since $p'(y)>0$ and $0<y<1$. Since $v(\ell+)=v(0+)<\infty$ and $p(r-)=p(\infty)=\infty$ by~\eqref{eq:Comparison-Scale-Infinity-Limit},~\cite[Prop.~5.5.32\textup{(iii)}, p.~350]{karatzas2012brownian} yields
\[
\mathbb P[S<\infty]=1.
\]
\end{proof}

\section{Derivation of Identity~(\ref{eq:LVFour-Weighted-Correction})} \label{sec:LVFour-Weighted-Correction-Derivation}

By Assumption~\ref{cond:G:singular:beta>1},~\eqref{eq:LVFour-Weighted-Singular-Generator} becomes
\begin{align*}
&L\bigg(\frac{a_4}{\beta+q-1}|v|^{-(\beta+q-1)}
+\frac{a_6}{\betat+q-1}|v|^{-(\betat+q-1)}\bigg)\\
&=-\underbrace{|v|^{-q}\la v,\grad G(v)\ra}_{\textup{(A)}}
-\underbrace{|v|^{-q}\la\grad G(v),\grad U(x)+\grad G(v)\ra}_{\textup{(B)}}\\
&\quad-\underbrace{\frac{\sigma}{2}\Big[
a_4(d-\beta-q-1)|v|^{-(\beta+q+1)}
+a_6(d-\betat-q-1)|v|^{-(\betat+q+1)}
\Big]}_{\textup{(C)}}.
\end{align*}
Also, multiplying~\eqref{eq:LVFour-Weighted-Cross-Generator} by $1/(1-q)$, we obtain
\begin{align}\label{eq:LVFour-Weighted-Cross-Generator2}
&\frac{1}{1-q}L\bigl(|v|^{-q}\la v,\grad U(x)\ra\bigr) \nonumber\\
&=\frac{|v|^{-q}}{1-q}\la v,\nabla^2U(x)v\ra
-\underbrace{|v|^{-q}\la v,\grad U(x)\ra}_{\textup{(A)}'} -\underbrace{\frac{|v|^{-q}}{1-q}
\la\grad U(x),\grad U(x)+\grad G(v)\ra}_{\textup{(B)}'}\nonumber\\
&\quad
+\underbrace{\frac{q}{1-q}|v|^{-q-2}\la v,\grad U(x)\ra
\la v,\grad U(x)+\grad G(v)\ra}_{\textup{(B)}''}
-\underbrace{\frac{\sigma q(d-q)}{2(1-q)}|v|^{-q-2}
\la v,\grad U(x)\ra}_{\textup{(C)}'}.
\end{align}

We first combine the terms marked $\textup{(A)}$ and $\textup{(A)}'$. We have
\begin{align}\label{eq:LVFour-Weighted-A}
-\textup{(A)}-\textup{(A)}'
=-|v|^{-q}\la v,\grad U(x)+\grad G(v)\ra.
\end{align}

We next combine the terms marked $\textup{(B)}$, $\textup{(B)}'$, and $\textup{(B)}''$. Their contribution is
\begin{align*}
-\textup{(B)}-\textup{(B)}'+\textup{(B)}''
&=-\frac{|v|^{-q}}{1-q}\bigg[
(1-q)\la\grad G(v),\grad U(x)+\grad G(v)\ra
+\la\grad U(x),\grad U(x)+\grad G(v)\ra\\
&\hspace{2cm}-q\frac{\la v,\grad U(x)\ra
\la v,\grad U(x)+\grad G(v)\ra}{|v|^2}
\bigg].
\end{align*}
Combining the first two terms and using
$\la v,\grad U(x)\ra
=\la v,\grad U(x)+\grad G(v)\ra-\la v,\grad G(v)\ra$,
we obtain
\begin{align*}
-\textup{(B)}-\textup{(B)}'+\textup{(B)}''
&=-\frac{|v|^{-q}}{1-q}\bigg[
|\grad U(x)+\grad G(v)|^2
-q\frac{\la v,\grad U(x)+\grad G(v)\ra^2}{|v|^2}\\
&\qquad
-q\la\grad G(v),\grad U(x)+\grad G(v)\ra
+q\frac{\la v,\grad G(v)\ra
\la v,\grad U(x)+\grad G(v)\ra}{|v|^2}
\bigg].
\end{align*}
Since, by Assumption~\ref{cond:G:singular:beta>1}, $\grad G(v)$ is parallel to $v$, we have
\[
\la\grad G(v),\grad U(x)+\grad G(v)\ra
=\frac{\la v,\grad G(v)\ra
\la v,\grad U(x)+\grad G(v)\ra}{|v|^2}.
\]
Hence, the last two terms in brackets cancel, and therefore
\begin{align}\label{eq:LVFour-Weighted-B}
-\textup{(B)}-\textup{(B)}'+\textup{(B)}''
&=-\frac{|v|^{-q}}{1-q}\left(
|\grad U(x)+\grad G(v)|^2
-q\frac{\la v,\grad U(x)+\grad G(v)\ra^2}{|v|^2}
\right).
\end{align}

Finally, using Assumption~\ref{cond:G:singular:beta>1}, we have
\[
\la v,\grad U(x)\ra
=\la v,\grad U(x)+\grad G(v)\ra
+a_4|v|^{1-\beta}+a_6|v|^{1-\betat}.
\]
Substituting this into the term marked $\textup{(C)}'$ and combining it with the term marked $\textup{(C)}$, we obtain
\begin{align*}
-\textup{(C)}-\textup{(C)}'
&=-\frac{\sigma a_4}{2(1-q)}
\Big[(1-q)(d-\beta-q-1)+q(d-q)\Big]
|v|^{-(\beta+q+1)}\\
&\quad-\frac{\sigma a_6}{2(1-q)}
\Big[(1-q)(d-\betat-q-1)+q(d-q)\Big]
|v|^{-(\betat+q+1)}\\
&\quad-\frac{\sigma q(d-q)}{2(1-q)}
|v|^{-q-2}\la v,\grad U(x)+\grad G(v)\ra.
\end{align*}
Since the expression in the first brackets simplifies to $d-1-\beta(1-q)$, while the second is the same expression with $\beta$ replaced by $\betat$, it follows that
\begin{align}\label{eq:LVFour-Weighted-C}
-\textup{(C)}-\textup{(C)}'
&=-\frac{\sigma a_4}{2(1-q)}
\bigl(d-1-\beta(1-q)\bigr)
|v|^{-(\beta+q+1)}\nonumber\\
&\quad-\frac{\sigma a_6}{2(1-q)}
\bigl(d-1-\betat(1-q)\bigr)
|v|^{-(\betat+q+1)}\nonumber\\
&\quad-\frac{\sigma q(d-q)}{2(1-q)}
|v|^{-q-2}\la v,\grad U(x)+\grad G(v)\ra.
\end{align}

Combining~\eqref{eq:LVFour-Weighted-A},~\eqref{eq:LVFour-Weighted-B}, and~\eqref{eq:LVFour-Weighted-C}, together with the Hessian term in~\eqref{eq:LVFour-Weighted-Cross-Generator2}, yields~\eqref{eq:LVFour-Weighted-Correction}.

\bibliographystyle{abbrv}
{\footnotesize\bibliography{GLE-bib}}

\end{document}